\documentclass[12pt]{ucthesis}
\usepackage{amsthm,amssymb,amsmath,amscd}
\usepackage{fancyhdr}

\author{Jae-Wook Chung}
\title{The Algebraic Structure of $n$-Punctured Ball Tangles}

\degreeyear{2005}

\degree{Doctor of Philosophy}

\chair{Professor Xiao-Song Lin}

\othermembers{Professor Gerhard Gierz\\Professor Bun Wong}

\prevdegrees{B.S. (Yeungnam University, Korea) 1990\\
M.S. (Yeungnam University, Korea) 1992}

\field{Mathematics}

\campus{Riverside}

\newtheorem{thm}{Theorem}[chapter]
\newtheorem{lm}[thm]{Lemma}

\newtheorem{pro}[thm]{Proposition}
\newtheorem{pro-def}[thm]{Proposition and Definition}

\theoremstyle{definition}
\newtheorem{df}[thm]{Definition}

\input epsf

\def\dsp{\def\baselinestretch{2.0}\large\normalsize}
\dsp

\fancypagestyle{plain}{
\lhead{} \chead{} \rhead{} \lfoot{}
\fancyfoot[C]{\vspace{-1pt}\thepage}

}

\begin{document}

\thispagestyle{empty}
\
\ssp

\begin{center}\

UNIVERSITY OF CALIFORNIA\\
RIVERSIDE\\
\vspace{.7in} The Algebraic Structure of $n$-Punctured Ball Tangles\\
\vspace{.7in} A Dissertation submitted in partial satisfaction\\
of the requirements for the degree of\\
\vspace{.2in} Doctor of Philosophy\\[12pt] in\\[12pt] Mathematics\\[12pt] by\\
[12pt] Jae-Wook Chung\\[12pt] June 2005

\end{center}

\vspace{1.5in}
\ssp

\noindent Dissertation Committee:\\
\indent Dr. Xiao-Song Lin, Chairperson\\
\indent Dr. Gerhard Gierz\\
\indent Dr. Bun Wong

\newpage

\ssp
\thispagestyle{empty}
\
\vspace{7.38in}
\begin{center}
Copyright by\\Jae-Wook Chung\\
2005
\end{center}

\newpage

\renewcommand{\thepage}{\roman{page}}

\setcounter{page}{4}

\begin{center}
\textbf{Acknowledgements}\\
\end{center}

First of all, I would like to express my gratitude to my thesis
adviser Professor Xiao-Song Lin. I sincerely appreciate his
continuous and devoted advice to me not only on my research, but
also in learning mathematics. Indeed, his sharp ideas and wide
knowledge on research fields have encouraged me to do research.

I would also like to thank Professor Gerhard Gierz and Professor
Bun Wong. I took several courses from them. Both of the professors
gave me very impressive lectures which have helped me understand
many important things that support my research work.

Finally, I would like to thank all other people who helped me in
the Department of Mathematics, University of California,
Riverside.

\newpage

\renewcommand{\thepage}{\roman{page}}

\setcounter{page}{5}
\
\vspace{2.0in}

\begin{center}
Dedicated to my wife Chang-Hee Lee and my daughter Eun-Ah Chung
\end{center}

\newpage

\thispagestyle{fancy}

\begin{center}
ABSTRACT OF THE DISSERTATION\\
\vspace{.3in} The Algebraic Structure of $n$-Punctured Ball Tangles\\
[12pt] by\\[12pt] Jae-Wook Chung\\
\vspace{.3in} Doctor of Philosophy, Graduate Program in Mathematics\\
University of California, Riverside, June 2005\\
Professor Xiao-Song Lin, Chairperson
\end{center}

\dsp
\vspace{.7in}

\indent We consider a class of topological objects in the 3-sphere
$S^3$ which will be called {\it $n$-punctured ball tangles}. Using
the Kauffman bracket at $A=e^{i \pi/4}$, an invariant for a
special type of $n$-punctured ball tangles is defined. The
invariant $F^n$ takes values in $PM_{2\times2^n}(\mathbb Z)$, that
is the set of $2\times 2^n$ matrices over $\mathbb Z$ modulo the
scalar multiplication of $\pm1$. This invariant leads to a
generalization of a theorem of D. Krebes which gives a necessary
condition for a given collection of tangles to be embedded in a
link in $S^3$ disjointly. Furthermore, we provide the general
formula to compute the invariant of $k_1 + \cdots + k_n$-punctured
ball tangle $T^n(T^{k_1(1)},\dots,T^{k_n(n)})$ determined by
$n,k_1,\dots,k_n$-punctured ball tangles
$T^n,T^{k_1(1)},\dots,T^{k_n(n)}$, respectively, when those of
$T^n, T^{k_1(1)},\dots,T^{k_n(n)}$ are given. Also, we consider
various connect sums among punctured ball tangles and provide the
formulas for the invariant of $T^{k_1(1)} +_h T^{k_2(2)}$ and that
of $T^{k_1(1)} +_v T^{k_2(2)}$ when those of $T^{k_1(1)}$ and
$T^{k_2(2)}$ are given, where $+_h$ and $+_v$ mean the outer
horizontal and the outer vertical connect sums of punctured ball
tangles, respectively. We also address the question of whether the
invariant $F^n$ is surjective onto $PM_{2\times2^n}(\mathbb Z)$.
We will show that the invariant $F^n$ is surjective when $n=0$.
When $n=1$, $n$-punctured ball tangles will also be called {\it
spherical tangles}. We show that ${\rm det}\,F^1(S)\equiv 0$ or 1
{\rm mod} 4 for every spherical tangle $S$. Thus, $F^n$ is not
surjective when $n=1$. In addition, we introduce monoid structures
on the class of $0$-punctured ball tangles and the class of
spherical tangles and show that the group generated by the
elementary operations on $PM_{2\times2}(\mathbb Z)$ induced by
those on the spherical tangles is isomorphic to a Coxeter group.

\newpage

\thispagestyle{fancy}
\tableofcontents

\newpage

\renewcommand{\thepage}{\arabic{page}}

\setcounter{page}{1}

\chapter{Introduction}

In this thesis, we will work in either the smooth or the piecewise
linear category. For basic terminologies of knot theory, see
\cite{A,BZ}.

The notion of {\it tangles} was introduced by J. Conway
\cite{Conway} as the basic building blocks of links in the
3-dimensional sphere $S^3$. Slightly abusing the notation, a
tangle $T$ is a pair $(B^3,T)$, where $B^3$ is a 3-dimensional
ball and $T$ is a proper 1-dimensional submanifold of $B^3$ with 2
non-circular components. The points in $\partial T\subset\partial
B^3$ will be fixed once and for all. Recall that a link $L$ is a
submanifold of $S^3$ homeomorphic to a disjoint union of several
copies of the circle $S^1$. A tangle $T=(B^3,T)$ can be embedded
in a link $L$ in $S^3$ if there is an embedding
$\phi:B^3\longrightarrow S^3$ such that $\phi(B^3)\cap L=\phi(T)$.
Using the Kauffman bracket at $A=e^{i \pi/4}$, a necessary
condition that one can embed a tangle $T$ in a link $L$ is given
by D. Krebes in \cite{K}.

One of the purposes of this thesis is to give a generalization to
Krebes' theorem. Suppose that $k$ tangles $T_i=(B^3,T_i)$,
$i=1,\dots,k$, are given. They can be embedded disjointly in a
link $L$ if there are embeddings $\phi_i:B^3\longrightarrow S^3$
such that $\phi_i(B^3)\cap L=\phi_i(T_i)$ for all $i$ and
$\phi_i(B^3)\cap\phi_j(B^3)=\emptyset$ for all $i,j$ with $i\neq
j$. A necessary condition similar to that of Krebes' will be given
for the existence of such a disjoint embedding of tangles in a
link (see Theorem 3.10).

In order to prove this generalization of Krebes' theorem, we will
study a class of topological objects in $S^3$ called {\it
$n$-punctured ball tangles}. This class of topological objects has
rich contents in the theory of {\it operads} \cite{MSS}, which is
beyond the scope of this thesis. Our main interest lies in a
special type of $n$-punctured ball tangles, which in the case of
$n=0$, corresponds exactly to Conway's notion of tangles in the
3-ball $B^3$. Using the Kauffman bracket at $A=e^{i \pi/4}$, we
will define an invariant for this special type of $n$-punctured
ball tangles. For an $n$-punctured ball tangle $T$, this invariant
$F^n(T)$ is an element in $PM_{2\times2^n}(\mathbb Z)$, that is
the set of $2\times 2^n$ matrices over $\mathbb Z$ modulo the
scalar multiplication of $\pm1$. When $n=0$, $F^n(T)$ is Krebes'
invariant.

Suppose now that we have $k$ tangles $T_i,\,i=1,\dots,k,$ embedded
disjointly in a link $L$. Let
$$F^0(T_i)=\left[\begin{matrix} p_i \\ q_i \end{matrix}\right],
\,i=1,\dots,k,$$ and let $\langle L\rangle$ be the Kauffman
bracket of $L$ at $A=e^{i \pi/4}$. Then Theorem 3.10 says that
$$\prod_{i=1}^k\,{\rm g.c.d.}\,(p_i,q_i)$$
divides $|\langle L\rangle |$. When $k=1$, this is exactly Krebes'
theorem.

The proof of Theorem 3.10 is based on the fact that the invariant
$F^n$ behaves well under the operadic composition of $n$-punctured
ball tangles.

In the second part of this thesis, we study the invariant $F^n$ in
some more details when $n=1$. In this case, $n$-punctured ball
tangles are called {\it spherical tangles}. For a given spherical
tangle $S$, ${\rm det}\,F^1(S)$ is a well-defined integer. Using a
theorem of S. Matveev, H. Murakami and Y. Nakanishi in \cite{Mat,
M-Y}, we will show that ${\rm det}\,F^1(S)$ is either 0 or 1
modulo 4 (Theorem 4.34). Thus, not every element in $PM_{2\times
2}(\mathbb Z)$ can be realized as $F^1(S)$ for some spherical
tangle $S$. This is in contrary with the case of $n=0$, where the
invariant $F^0$ is onto.

We organize the thesis as follows: In Chapter 2, we formally
define the notion of $n$-punctured ball tangles. We also recall
the Kauffman bracket at $A=e^{i \pi/4}$ and Krebes' theorem in
this Chapter. In Chapter 3, we define our invariant $F^n$ for a
special class of $n$-punctured ball tangles. A key result is about
the behave of the invariant $F^n$ under operadic composition of
$n$-punctured ball tangles (Theorems 3.6). Our generalization of
Krebes' theorem (Theorem 3.10) will follow easily from this
result. Furthermore, given an $n$-punctured ball tangle $T^n$ and
$k_1,\dots,k_n$-punctured ball tangles
$T^{k_1(1)},\dots,T^{k_n(n)}$, respectively, we consider the $k_1
+ \cdots + k_n$-punctured ball tangle
$T^n(T^{k_1(1)},\dots,T^{k_n(n)})$, where $n \in \mathbb N$ and
$k_1,\dots,k_n \in \mathbb N \cup \{0\}$. In this case, we show
how to compute the invariant
$F^{k_1+\cdots+k_n}(T^n(T^{k_1(1)},\dots,T^{k_n(n)}))$ if
$F^n(T^n), F^{k_1}(T^{k_1(1)}),\dots,F^{k_n}(T^{k_n(n)})$ are
given (Theorem 3.8). Also, we consider the horizontal connect sum
$T^{k_1(1)} +_h T^{k_2(2)}$ and the vertical connect sum
$T^{k_1(1)} +_v T^{k_2(2)}$ of $k_1$ and $k_2$-punctured ball
tangles $T^{k_1(1)}$ and $T^{k_2(2)}$, respectively, and provide
the formulas for the invariants $F^{k_1+k_2}(T^{k_1(1)} +_h
T^{k_2(2)})$ and $F^{k_1+k_2}(T^{k_1(1)} +_v T^{k_2(2)})$ when
$F^{k_1}(T^{k_1(1)})$ and $F^{k_2}(T^{k_2(2)})$ are given,
respectively (Theorem 3.9). In Chapter 4, we introduce $2$ similar
monoid operations on the class of $0$-punctured ball tangles which
are connect sums and $2$ similar monoid actions on it by the
spherical tangles which is also a monoid with respect to the
composition. For the spherical tangles, we consider the
composition as a monoid operation and $8$ connect sums as monoid
actions by the $0$-punctured ball tangles since spherical tangles
have $2$ holes which are inside and outside. In addition, we show
that the group $G(F)$ generated by the elementary operations on
$PM_{2\times2}(\mathbb Z)$ induced by those on the spherical
tangles is isomorphic to the Coxeter
group $C_M$ with the Coxeter matrix $M=\begin{pmatrix} 1 & 4 & 2 \\
4 & 1 & 2 \\ 2 & 2 & 1 \end{pmatrix}$ (Theorem 4.25). Also, we
study the surjectivity of the invariant $F^n$ in the case of
$n=0,1$. As mentioned before, we will show that $F^n$ is
surjective when $n=0$ but not surjective when $n=1$. In the final
Chapter, we pose some questions related with this work which we do
not know how to answer at this moment.

Notice that D. Ruberman has given a topological interpretation of
Krebes' theorem \cite{R}. We don't know if our generalization of
Krebes' theorem could have a similar topological interpretation.
In particular, it will be very nice if there is a topological
interpretation of the restriction on ${\rm det}\,F^1(S)$ for
spherical tangles $S$ (Theorem 4.34).

\chapter{General definitions}

\section{$n$-punctured ball tangles}

We define a topological object in the 3-dimensional sphere $S^3$
called an {\it $n$-punctured ball tangle} or, simply, an {\it
$n$-tangle}. To study this object, we consider a model for a class
of objects and an equivalence relation on it.

\begin{df} Let $n$ be a nonnegative integer, and let $H_0$ be a
3-dimensional closed ball, and let $H_1,\dots,H_n$ be pairwise
disjoint 3-dimensional closed balls contained in the interior
${\rm Int}(H_0)$ of $H_0$. For each $k \in \{0,1,\dots,n\}$, take
$2m_k$ distinct points $a_{k1},\dots,a_{k2m_k}$ of $\partial H_k$
for some positive integer $m_k$. Then a 1-dimensional proper
submanifold $T$ of $H_0 - \bigcup_{i=1}^n {\rm Int}(H_i)$ is
called an $n$-punctured ball tangle with respect to $(H_k)_{0 \leq
k \leq n}$, $(m_k)_{0 \leq k \leq n}$, and
$((a_{k1},\dots,a_{k2m_k}))_{0 \leq k \leq n}$ if $\partial T =
\bigcup_{k=0}^n \{a_{k1},\dots,a_{k2m_k}\}$. Hence, $\partial T
\cap \partial H_k = \{a_{k1},\dots,a_{k2m_k}\}$ for each $k \in
\{0,1,\dots,n\}$. Note that an $n$-tangle $T$ with respect to
$(H_k)_{0 \leq k \leq n}$, $(m_k)_{0 \leq k \leq n}$, and
$((a_{k1},\dots,a_{k2m_k}))_{0 \leq k \leq n}$ can be regarded as
a $5$-tuple $(n,(H_k)_{0 \leq k \leq n},(m_k)_{0 \leq k \leq
n},((a_{k1},\dots,a_{k2m_k}))_{0 \leq k \leq n},T)$.
\end{df}

\begin{pro} Let $n \in \mathbb N \cup \{0\}$, and let $\textbf{\textit{nPBT}}$
be the class of all $n$-punctured ball tangles with respect to
$(H_k)_{0 \leq k \leq n}$, $(m_k)_{0 \leq k \leq n}$, and
$((a_{k1},\dots,a_{k2m_k}))_{0 \leq k \leq n}$, and let $X = H_0 -
\bigcup_{i=1}^n {\rm Int}(H_i)$. Define $\cong$ on
$\textbf{\textit{nPBT}}$ by $T_1 \cong T_2$ if and only if there
is a homeomorphism $h:X \rightarrow X$ such that $h|_{\partial
X}=Id_X|_{\partial X}$, $h(T_1)=T_2$, and $h$ is isotopic to
$Id_X$ relative to the boundary $\partial X$ for all $T_1,T_2 \in
\textbf{\textit{nPBT}}$. Then $\cong$ is an equivalence relation
on $\textbf{\textit{nPBT}}$, where $Id_X$ is the identity map from
$X$ to $X$.
\end{pro}

\begin{proof} Note that $T_1 \cong T_2$ if and only if there are a
homeomorphism $h:X \rightarrow X$ with $h|_{\partial
X}=Id_X|_{\partial X}$ and $h(T_1)=T_2$ and a continuous function
$H:X \times I \rightarrow X$ such that $H(\_,t):X \rightarrow X$
is a homeomorphism with $H(\_,t)|_{\partial X}=Id_X|_{\partial X}$
for each $t \in I$ and $H(\_,0)=Id_X$ and $H(\_,1)=h$, where
$I=[0,1]$. Let us denote $H(x,t)$ by $H_t(x)$ for all $x \in X$
and $t \in I$, so $H(\_,t)=H_t$ for each $t \in I$.

For every $T \in \textbf{\textit{nPBT}}$, $T \cong T$ since $Id_X$
and the 1st projection $\pi_1:X \times I \rightarrow X$ satisfy
the condition. Suppose that $T_1 \cong T_2$ and $h:X \rightarrow
X$ is a homeomorphism with $h|_{\partial X}=Id_X|_{\partial X}$
such that $h(T_1)=T_2$ and $H:X \times I \rightarrow X$ is a
continuous function such that $H_t:X \rightarrow X$ is a
homeomorphism with $H_t|_{\partial X}=Id_X|_{\partial X}$ for each
$t \in I$ and $H_0=Id_X$ and $H_1=h$. Define $H':X \times I
\rightarrow X$ by $H'(x,t)=H^{-1}_t(x)$ for all $x \in X$ and $t
\in I$. Then $h^{-1}$ and $H'$ make $T_2 \cong T_1$. To show the
transitivity of $\cong$, suppose that $T_1 \cong T_2$ and $h:X
\rightarrow X$ is a homeomorphism with $h|_{\partial
X}=Id_X|_{\partial X}$ such that $h(T_1)=T_2$ and $H:X \times I
\rightarrow X$ is a continuous function such that $H_t:X
\rightarrow X$ is a homeomorphism with $H_t|_{\partial
X}=Id_X|_{\partial X}$ for each $t \in I$ and $H_0=Id_X$ and
$H_1=h$ and $T_2 \cong T_3$ and $h':X \rightarrow X$ is a
homeomorphism with $h'|_{\partial X}=Id_X|_{\partial X}$ such that
$h'(T_2)=T_3$ and $H':X \times I \rightarrow X$ is a continuous
function such that $H'_t:X \rightarrow X$ is a homeomorphism with
$H'_t|_{\partial X}=Id_X|_{\partial X}$ for each $t \in I$ and
$H'_0=Id_X$ and $H'_1=h'$. Define $H'':X \times I \rightarrow X$
by $H''(x,t)=(H'_t \circ H_t)(x)$ for all $x \in X$ and $t \in I$.
Then $h' \circ h$ and $H''$ make $T_1 \cong T_3$. Therefore,
$\cong$ is an equivalence relation on $\textbf{\textit{nPBT}}$.
\end{proof}

\begin{df} Let $T_1$ and $T_2$ be $n$-punctured ball tangles
in $\textbf{\textit{nPBT}}$. Then $T_1$ and $T_2$ are said to be
equivalent or of the same isotopy type if $T_1 \cong T_2$. Also,
for each $n$-punctured ball tangle $T$ in
$\textbf{\textit{nPBT}}$, the equivalence class of $T$ with
respect to $\cong$ is denoted by $[T]$. By the context, without
any confusion, we will also use $T$ for $[T]$.
\end{df}

There are many models for a class of $n$-punctured ball tangles.
It is convenient to use normalized ones. One model for a class of
$n$-punctured ball tangles is as follows:

(1) $H_0 = B((\frac{n+1}{2}, 0,0),\frac{n+1}2)$ and $H_i = B((i,
0,0),\frac{1}{3})$ for each $i \in \{1,\dots,n\}$.

(2) $a_{k1},\dots,a_{k2m_k}$ are $2m_k$ distinct points of
$\partial H_k$ in the $xy$-plane for each $k \in \{0,1,\dots,n\}$.

(3) $T$ is a 1-dimensional proper submanifold of $H_0 -
\bigcup_{i=1}^n {\rm Int}(H_i)$ such that $\partial T =
\bigcup_{k=0}^n \{a_{k1},\dots,a_{k2m_k}\}$, where $B((x,y,z),r)$
is the 3-ball in $\mathbb{R}^3$ with radius $r$ at $(x,y,z)$.

In order to study an $n$-punctured ball tangle $T$ through its
diagram $D$, we consider the $xy$-projection $P_{xy}:\mathbb{R}^3
\rightarrow \mathbb{R}^3$ defined by $P_{xy}(x,y,z)=(x,y,0)$ for
all $x,y,z, \in \mathbb R$. A point $p$ of the image $P_{xy}(T)$
is called a multiple point of $T$ if the cardinality of
$P_{xy}^{-1}(p) \cap T$ is greater than 1. In particular, $p$ is
called a double point of $T$ if the cardinality of $P_{xy}^{-1}(p)
\cap T$ is 2. If $p$ is a double point of $T$, then
$P_{xy}^{-1}(p) \cap T$ is called the crossing of $T$
corresponding to $p$ and the point in the crossing whose
$z$-coordinate is greater is called the overcrossing of $T$
corresponding to $p$ and the other is called the undercrossing.

An $n$-punctured ball tangle $T$ is said to be in regular position
if the only multiple points of $T$ are double points and each
double point of $T$ is a transversal intersection of the images of
two arcs of $T$ and $P_{xy}(T - \partial T) \cap (\partial
(P_{xy}(H_0)) \cup \bigcup_{i=1}^n P_{xy}(H_i)) = \emptyset$.

Note that for each $T \in \textbf{\textit{nPBT}}$, there is $T'
\in \textbf{\textit{nPBT}}$ such that $T'$ is in regular position
and $T' \cong T$. Furthermore, $T'$ has a finite number of
crossings.

Consider the image $P_{xy}(T)$ of an $n$-punctured ball tangle $T$
in regular position. For each double point of $T$, take a
sufficiently small closed ball centered at the double point such
that the intersection of $P_{xy}(T)$ and the closed ball is an
X-shape on the $xy$-plane. We may assume that the closed balls are
pairwise disjoint. Now, modify the interiors of the closed balls
keeping the image $P_{xy}(T)$ to assign crossings corresponding to
the crossings of $T$. As a result, we have a representative $D$ of
$T$ which is `almost planar' and $P_{xy}(D)=P_{xy}(T)$. $D$ is
called a diagram of $T$ and we usually use this representative.

To deal with diagrams of $n$-punctured ball tangles in the same
isotopy type, we need Reidemeister moves among them. For link
diagrams or ball tangle diagrams, we have 3 kinds of Reidemeister
moves. However, we need one and only one more kind of moves which
are called the Reidemeister moves of type IV.

The Reidemeister moves for diagrams of $n$-punctured ball tangles
are illustrated in the following figure.

\vskip 0.2in

\bigskip
\centerline{\epsfxsize=5.8 in \epsfbox{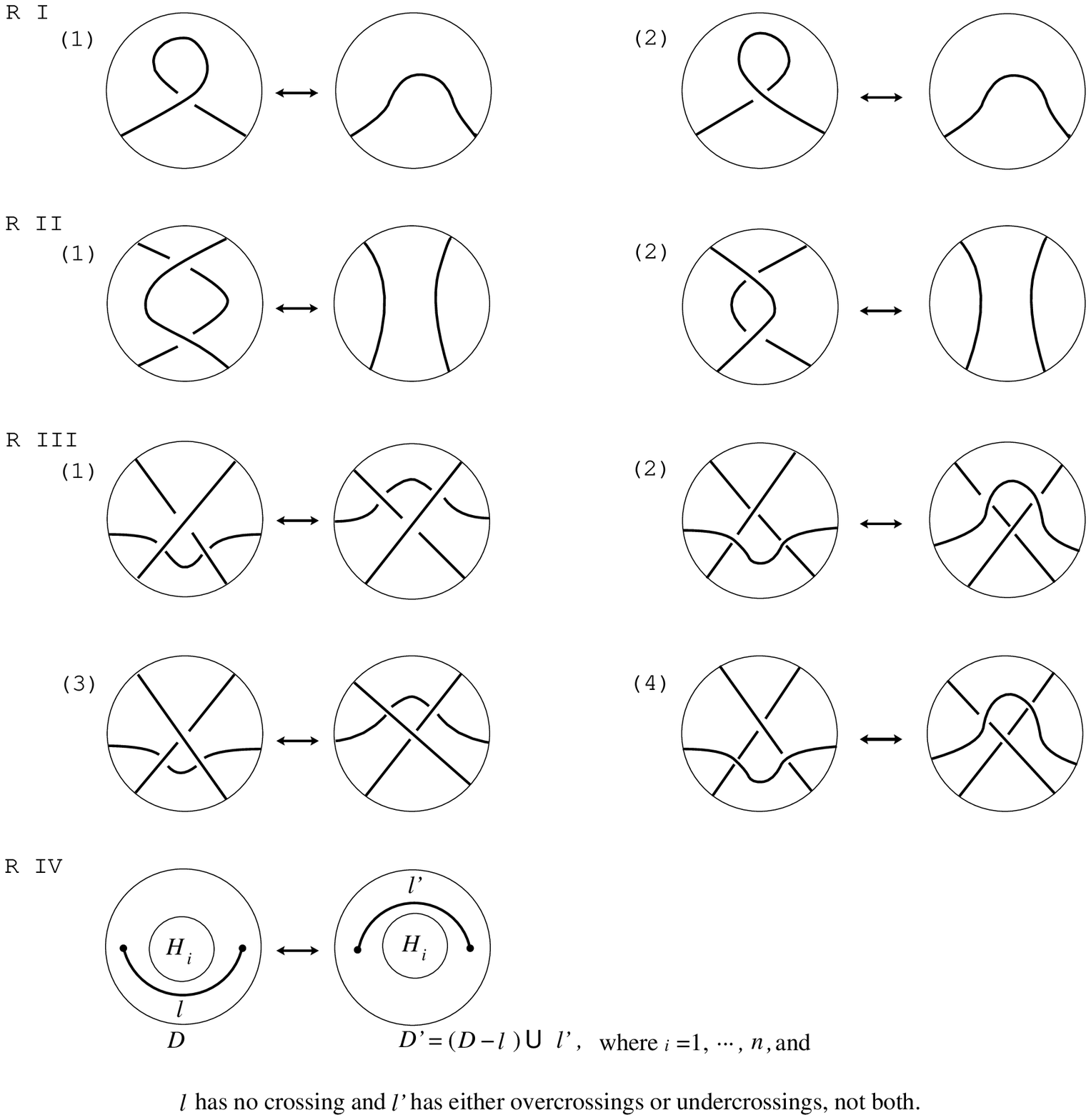}}
\medskip
\centerline{\small Figure 1. Tangle Reidemeister moves.}
\bigskip

Like link diagrams, tangle diagrams also have Reidemeister Theorem
involving the Reidemeister moves of type IV. Let us call
Reidemeister moves including type IV Tangle Reidemeister moves.

\begin{thm} Let $n$ be a nonnegative integer, and let $D_1$ and $D_2$
be diagrams of $n$-punctured ball tangles. Then $D_1 \cong D_2$ if
and only if $D_2$ can be obtained from $D_1$ by a finite sequence
of Tangle Reidemeister moves.
\end{thm}

Remark that, even though we may have different models for
$n$-punctured ball tangles, we may regard them as the same
$n$-punctured ball tangle if there are suitable model equivalences
among them.

\section{Kauffman bracket at $A=e^{i \pi/4}$ and monocyclic state
of link diagram}

Our invariant is based on the Kauffman bracket at $A=e^{i \pi/4}$.
In this section, we recall the Kauffman bracket which is a regular
isotopy invariant of link diagrams. That is, it will not be
changed under Reidemeister moves of type II and III.

Assume that $L$ is a link diagram with $n$ crossings and $c$ is a
crossing of $L$. Take a sufficiently small disk at the projection
of $c$ to get an X-shape on the projection plane of $L$. Now, we
have 4 regions in the disk. Rotate counterclockwise the projection
of the over-strand in the disk which is an arc of $L$ containing
the overcrossing for $c$ to pass over 2 regions. These 2 regions
and the other 2 regions are called the $A$-regions and the
$B$-regions of $c$, respectively. We consider 2 ways of splitting
the double point in the disk. $A$-type splitting is to open a
channel between the $A$-regions so that we have 1 $A$-region and 2
$B$-regions in the disk and $B$-type splitting is to open a
channel between the $B$-regions so that we have 2 $A$-regions and
1 $B$-region in the disk. A choice of how to destroy all of $n$
double points in the projection of $L$ by $A$-type or $B$-type
splitting is called a state of $L$.

Notice that we regard a state $\sigma$ of the link diagram $L$
with $n$ crossings as a function $\sigma:\{c_1,\dots,c_n\}
\rightarrow \{A,B\}$, where $\{c_1,\dots,c_n\}$ is the set of all
crossings of $L$ and $\{A,B\}$ is the set of $A$-type and $B$-type
splitting functions, respectively. Therefore, a link diagram $L$
with $n$ crossings has exactly $2^n$ states of it. Apply a state
$\sigma$ to $L$ in order to change $L$ to a diagram $L_\sigma$
without any crossing.

\begin{df} Let $L$ be a link diagram. Then the Kauffman bracket
$\langle L \rangle_A$, or simply, $\langle L \rangle$, is defined
by $$\langle L \rangle_A=\sum_{\sigma \in S}
A^{\alpha(\sigma)}(A^{-1})^{\beta(\sigma)}
(-A^2-A^{-2})^{d(\sigma)-1},$$ where $S$ is the set of all states
of $L$, $\alpha(\sigma)=|\sigma^{-1}(A)|$,
$\beta(\sigma)=|\sigma^{-1}(B)|$, and $d(\sigma)$ is the number of
circles in $L_\sigma$.
\end{df}

We have the following \textit{skein relation of the Kauffman
bracket}.

\begin{pro} Let $L$ be a link diagram, and let $c$ be a crossing
of $L$. Then if $L_A$ and $L_B$ are link diagrams obtained from
$L$ by $A$-type splitting and $B$-type splitting only at $c$,
respectively, then $\langle L \rangle = A\langle L_A \rangle +
A^{-1}\langle L_B \rangle$.
\end{pro}

\begin{proof} Suppose that $S$ is the set of all states
of $L$ and $S_A=\{\sigma \in S|\sigma(c)=A\}$ and $S_B=\{\tau \in
S|\tau(c)=B\}$. Then $\langle L \rangle = A\sum_{\sigma \in S_A}
A^{\alpha(\sigma)-1}(A^{-1})^{\beta(\sigma)}
(-A^2-A^{-2})^{d(\sigma)-1} + A^{-1}\sum_{\tau \in S_B}
A^{\alpha(\tau)}(A^{-1})^{\beta(\tau)-1} (-A^2-A^{-2})^{d(\tau)-1}
= A\langle L_A \rangle + A^{-1}\langle L_B \rangle$ because $S$ is
the disjoint union of $S_A$ and $S_B$. This proves the
proposition.
\end{proof}

The following lemma is useful in our discussion of the Kauffman
bracket.

\begin{lm} Let $L$ be a link diagram. Then states $\sigma$ and
$\sigma'$ of $L$ are of the same parity, i.e., $d(\sigma) \equiv
d(\sigma')$ {\rm mod} 2, if and only if $\sigma$ and $\sigma'$
differ at an even number of crossings, where $d(\sigma)$ and
$d(\sigma')$ are the numbers of circles in $L_\sigma$ and
$L_{\sigma'}$, respectively.
\end{lm}

\begin{proof} Let $\sigma$ be a state of a link diagram $L$ with
$n$ crossings $c_1,\dots,c_n$. Change the value of $\sigma$ at
only one crossing $c_i$ to get another state $\sigma_i$ and
observe what happens to $d(\sigma_i)$, where $1 \leq i \leq n$. We
claim that $\sigma$ and $\sigma_i$ have different parities, more
precisely, $d(\sigma)=d(\sigma_i) \pm1$. Hence, we will have
$d(\sigma) \equiv d(\sigma_i)+1$ {\rm mod} 2. Now, to consider
$\sigma(c_i)$ and $\sigma_i(c_i)$, take a sufficiently small
neighborhood $B_i$ at the projection of $c_i$ so that the
intersection of ${\rm Int}(B_i)$ and the set of all double points
of $L$ is the projection of $c_i$ and the intersection of
$\partial B_i$ and the projection of $L$ has exactly 4 points on
the projection plane of $L$ which are not double points of $L$.

{\it Case 1}. If these 4 points are on a circle in $L_\sigma$,
then
$$d(\sigma_i)=d(\sigma)+1.$$

{\it Case 2}. If two of 4 points are on a circle and the other
points are on another circle in $L_\sigma$, then
$$d(\sigma_i)=d(\sigma)-1.$$

Now, it is easy to show the lemma. Suppose that $\sigma$ and
$\sigma'$ are states of $L$ which differ at $k$ crossings of $L$
for some $k \in \{0,1,\dots,n\}$. Then $d(\sigma') \equiv
d(\sigma)+k$ {\rm mod} 2. If $d(\sigma) \equiv d(\sigma')$ {\rm
mod} 2, then $k$ is even. Conversely, if $d(\sigma) \equiv
d(\sigma')+1$ {\rm mod} 2, then $k+1$ is even, that is, $k$ is
odd. This proves the lemma.
\end{proof}

Following \cite{K}, a state $\sigma$ of a link diagram $L$ is
called a monocyclic state of $L$ if $d(\sigma)=1$. That is, we
have only one circle when we remove all crossings of $L$ by
$\sigma$.

From now on, we consider only the Kauffman brackets at $A=e^{i
\pi/4}$. Since $|A|=1$, the determinant $|\langle L \rangle|$ of
$L$ is an isotopy invariant. It is easy to show that Reidemeister
move of type I dose not change $|\langle L \rangle|$ by the skein
relation of Kauffman bracket and $|A|=1$.

Notice that $-A^2-A^{-2}=0$ if $A=e^{i \pi/4}$. Therefore,
$$\langle L \rangle=\sum_{\sigma \in M}
A^{\alpha(\sigma)-\beta(\sigma)},$$ where $M$ is the set of all
monocyclic states of $L$.

As a corollary of Lemma 2.7, monocyclic states $\sigma$ and
$\sigma'$ of $L$ differ at an even number of crossings.

\begin{lm} If $L$ is a link diagram, then there are $p \in \mathbb Z$ and
$u \in \mathbb C$ such that $u^8=1$ and $\langle L \rangle=pu$.
\end{lm}

\begin{proof} Suppose that $\sigma$ and $\sigma'$ are states of a
link diagram $L$ such that $\sigma$ and $\sigma'$ differ at only
one crossing. Then either
$\alpha(\sigma')-\beta(\sigma')=(\alpha(\sigma)+1)-(\beta(\sigma)-1)
=(\alpha(\sigma)-\beta(\sigma))+2$ or
$\alpha(\sigma')-\beta(\sigma')=(\alpha(\sigma)-1)-(\beta(\sigma)+1)
=(\alpha(\sigma)-\beta(\sigma))-2$. Hence, either
$A^{\alpha(\sigma')-\beta(\sigma')}=+iA^{\alpha(\sigma)-\beta(\sigma)}$
or
$A^{\alpha(\sigma')-\beta(\sigma')}=-iA^{\alpha(\sigma)-\beta(\sigma)}$.
That is, $A^{\alpha(\sigma')-\beta(\sigma')}=\pm
iA^{\alpha(\sigma)-\beta(\sigma)}$.

If $0 \leq k \leq c(L)$ and $\sigma''$ is a state of $L$ such that
$\sigma$ and $\sigma''$ differ at $k$ crossings, where $c(L)$ is
the number of crossings of $L$, then
$A^{\alpha(\sigma'')-\beta(\sigma'')}=\pm i^k
A^{\alpha(\sigma)-\beta(\sigma)}$ because there are exactly $2^k$
sequences with $k$ terms consisting of $+i$ and $-i$ and the
product of all terms of each of the sequences is either $+i^k$ or
$-i^k$. Hence, $A^{\alpha(\sigma'')-\beta(\sigma'')}=\pm
iA^{\alpha(\sigma)-\beta(\sigma)}$ if $k$ is odd and
$A^{\alpha(\sigma'')-\beta(\sigma'')}=\pm
A^{\alpha(\sigma)-\beta(\sigma)}$ if $k$ is even.

Now, let us take a monocyclic state $\sigma_0$ of $L$, and let
$u=A^{\alpha(\sigma_0)-\beta(\sigma_0)}$. Then $u^8=1$ and
$\langle L \rangle=pu$ for some $p \in \mathbb Z$ by the corollary
above. This proves the lemma.
\end{proof}

\section{Krebes' Theorem}

In this subsection, we introduce some notations and Krebes'
Theorem \cite{K}.

\begin{df} A ball tangle $B$, which is a 0-punctured ball tangle with
$m_0=2$, is said to be embedded in a link $L$ if there are a
diagram $D_B$ of $B$, a diagram $D_L$ of $L$, and a 3-dimensional
closed ball $H$ such that $H \cap D_L$ and $D_B$ are of the same
isotopy type.
\end{df}

Given a ball tangle diagram $B$, we consider 3 kinds of closures
as in Figure 2 a).

\bigskip
\centerline{\epsfxsize=5 in \epsfbox{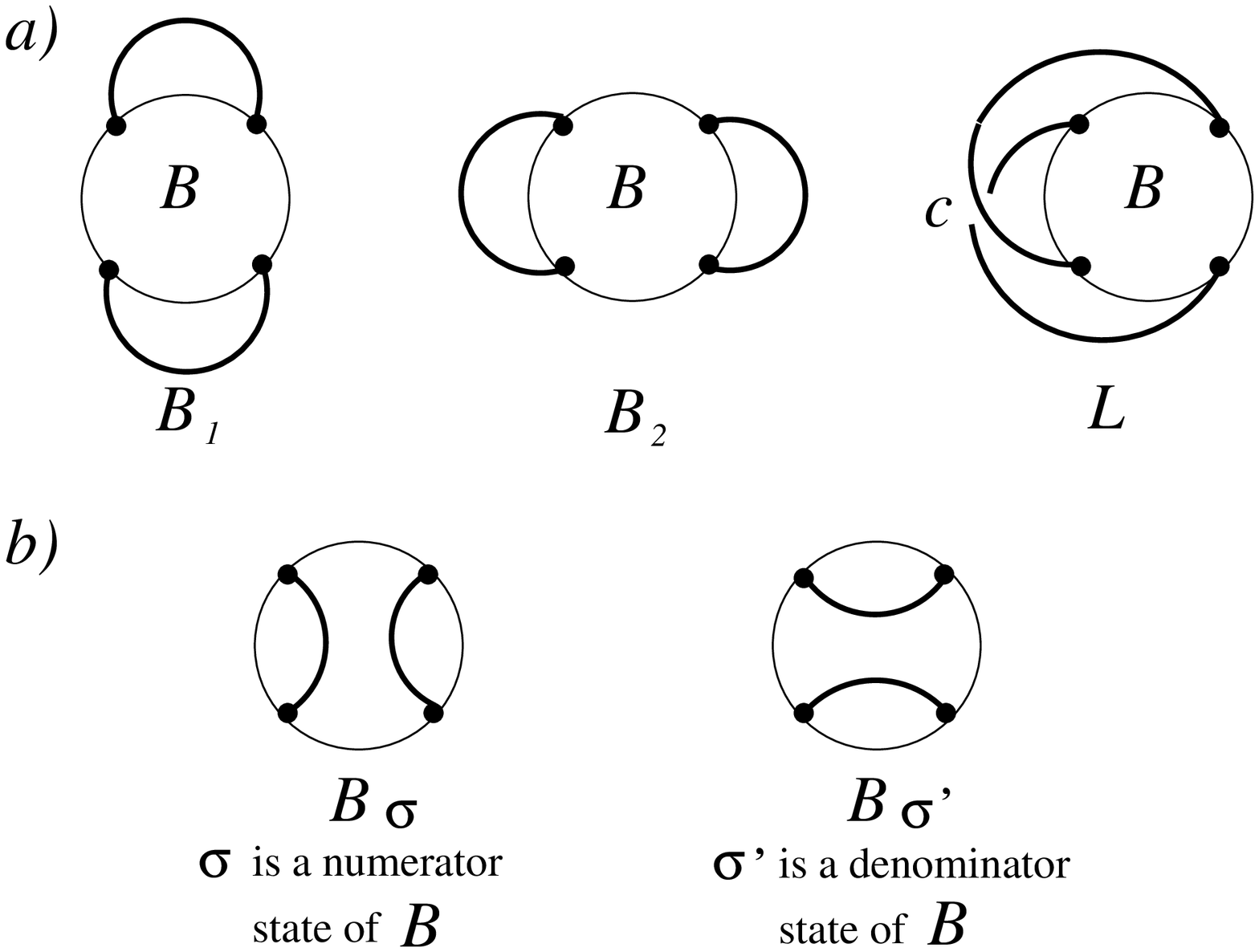}}
\medskip
\centerline{\small Figure 2. a) Closures, b) Diagrams by a
numerator state and a denominator state.}
\bigskip

The link diagrams $B_1$ and $B_2$ are called the numerator closure
and the denominator closure of $B$, respectively. A monocyclic
state of $B_1$ is called a numerator state of $B$ and that of
$B_2$ is a denominator state of $B$.

Notice that a numerator state $\sigma$ and a denominator state
$\sigma'$ of a ball tangle diagram $B$ differ at an odd number of
crossings. To see this, we think of a diagram of another closure
$L$ of $B$ which has only one more crossing $c$ at the outside of
ball containing $B$ (See $L$ in Figure 2). We have two monocyclic
states of $L$ from the numerator state $\sigma$ and the
denominator state $\sigma'$, respectively, which differ at $c$.
Hence, $\sigma$ and $\sigma'$ differ at an odd number of
crossings.

The following notations throughout the rest of the thesis:

\noindent$\bullet$ $\Phi = \{z \in \mathbb{C}\,|\,z^8 =1\} =
\{A^k\,|\,k \in\mathbb{Z}\}$, i.e. $\Phi$ is the set of 8-th roots
of unity; and $\mathbb{Z}\Phi=\{kz\,|\,k\in\mathbb{Z},z\in\Phi\}$.

\noindent$\bullet$ $M_{n\times m}(\mathbb{Z})$ is the set of all
$n\times m$ matrices over $\mathbb{Z}$, and $PM_{n\times
m}(\mathbb{Z})$ is the quotient of $M_{n\times m}(\mathbb{Z})$
under the scalar multiplication by $\pm1$.

\noindent$\bullet$ $\textbf{\textit{BT}}$ is the class of diagrams
of 0-punctured ball tangles with $m_0=2$ (i.e. ball tangles).

\noindent$\bullet$ $\textbf{\textit{ST}}$ is the class of diagrams
of 1-punctured ball tangles with $m_0=m_1=2$ (they will be called
spherical tangles).

\begin{pro} If $a,b,k,l \in \mathbb{Z}$, then $aA^k + bA^l \in \mathbb{Z}\Phi$
if and only if $ab=0$ or $k \equiv l$ {\rm mod} 4.
\end{pro}

\begin{proof} Suppose that $ab \neq 0$ and
$k-l \equiv t$ {\rm mod} 4 for some $t \in \{1,2,3\}$. then
$k-l=4m+t$ for some $m \in \mathbb{Z}$, hence, $A^k=(-1)^mA^tA^l$.
Therefore, $aA^k + bA^l$ is not in $\mathbb{Z}\Phi$. Conversely,
if $ab=0$ or $k \equiv l$ {\rm mod} 4, then $aA^k + bA^l \in
\mathbb{Z}\Phi$.
\end{proof}

We have
$$\langle L\rangle=A\langle B_1\rangle+A^{-1}\langle B_2\rangle\in\mathbb{Z}
\Phi$$ for the links $L$, $B_1$, and $B_2$ in Figure 2. If
$\langle B_1\rangle = pA^k$ and $\langle B_2\rangle=qA^l$, by
Proposition 2.10, we have $l\equiv k+2$ {\rm mod} 4. So there is a
unique $(\alpha,\beta) \in \mathbb{Z}^2$ such that
$$\left\{\begin{pmatrix} z\langle B_1 \rangle \\ iz\langle B_2 \rangle
\end{pmatrix}\,|\,z \in \Phi\right\} \cap M_{2\times 1}(\mathbb{Z})=\left\{
\begin{pmatrix} \alpha \\ \beta \end{pmatrix}, \begin{pmatrix} -\alpha \\ -\beta
\end{pmatrix}\right\}:=\left[\begin{matrix} \alpha \\ \beta
\end{matrix}\right] \in PM_{2\times 1}(\mathbb{Z}).$$

\begin{df} Define $f:\textbf{\textit{BT}} \rightarrow PM_{2\times1}
(\mathbb{Z})$ by
$$f(B)=\left\{\begin{pmatrix} z\langle B_1 \rangle \\ iz\langle B_2
\rangle \end{pmatrix}\,|\,z \in \Phi\right\} \cap
M_{2\times1}(\mathbb{Z}) \in PM_{2\times 1}(\mathbb{Z})$$ for each
$B \in \textbf{\textit{BT}}$. This is Krebes' tangle invariant.
\end{df}

Notice that Reidemeister move of type I dose not change $f(B)$. So
$f(B)$ is a ball tangle invariant. The following lemmas about the
ball tangle invariant $f$ are proved in \cite{K}.

\begin{lm} If $B^{(1)}$ and $B^{(2)}$ are diagrams of ball tangles with
$f(B^{(1)})=\left[\begin{matrix} p \\ q \end{matrix}\right]$ and
$f(B^{(2)})=\left[\begin{matrix} r \\ s \end{matrix}\right]$, then
$f(B^{(1)} +_h B^{(2)})=\left[\begin{matrix} ps+qr \\ qs
\end{matrix}\right]$,
where $B^{(1)} +_h B^{(2)}$ stands for the horizontal addition of
ball tangles (see Figure 10 a)).
\end{lm}

\begin{lm} If $B$ is a diagram of ball tangle with
$f(B)=\left[\begin{matrix} p \\ q \end{matrix}\right]$, then we
have
$$f(B^*)=\left[\begin{matrix} p \\ -q \end{matrix}\right]\qquad\text{and}\qquad
f(B^R)=\left[\begin{matrix} q \\ -p \end{matrix}\right],$$ where
$B^*$ is the mirror image of $B$ and $B^R$ is the $90^{\circ}$
counterclockwise rotation of $B$ on the projection plane.
\end{lm}

\begin{thm} {\rm (Krebes \cite{K})} If $L$ is a link and $B$ is a ball
tangle embedded in $L$ with $f(B)=\left[\begin{matrix} p \\ q
\end{matrix}\right]$, then ${\rm g.c.d.}\,(p,q)$ divides
$|\langle L \rangle|$.
\end{thm}

\chapter{The special case of $n$-punctured ball tangles}

Let $n$ be a positive integer. Then an $n$-punctured ball tangle
$T^n$ with $(H_k)_{0 \leq k \leq n}$ and $(m_k)_{0 \leq k \leq n}$
can be regarded as an $n$-variable function
$T^n:\textbf{\textit{A}}_1\times\cdots\times\textbf{\textit{A}}_n
\rightarrow \textbf{\textit{T}}$ defined as $T^n(X_1,\dots,X_n)$
is a tangle filled up in the $i$-th hole $H_i$ of $T^n$ by $X_i
\in \textbf{\textit{A}}_i$ for each $i \in \{1,\dots,n\}$, where
$\textbf{\textit{A}}_i$ is a class of $t_i$-punctured ball tangles
with $(m_k^i)_{0 \leq k \leq t_i}$ such that $m_i=m_0^i$ for each
$i \in \{1,\dots,n\}$ and $\textbf{\textit{T}}$ is a class of
tangles. However, this representation of $n$-punctured ball
tangles as $n$-variable functions is not perfect in the sense that
$n$-punctured ball tangles are equivalent only if they induce the
same function. On the other hand, $n$-punctured ball tangles which
induce the same function need not be equivalent. That is, we can
say that tangles are stronger than functions. Roughly speaking,
the class of $n$-punctured ball tangles as only $n$-variable
functions gives us an \textit{operad}, a mathematical device which
describes algebraic structure of many varieties and in various
categories. See \cite{MSS}.

\section{$n$-punctured ball tangles with $m_k=2$ for each $k \in
\{0,1,\dots,n\}$ and their invariants $F^n$}

From now on, we consider only $n$-punctured ball tangles with
$m_0=m_1=\cdots=m_n=2$.

To construct the invariant $F^n$ of $n$-punctured ball tangle
$T^n$, let us regard $T^n$ as a `hole-filling function', in sense
described as above $T^n:\textbf{\textit{BT}}^n \rightarrow
\textbf{\textit{BT}}$, where
$\textbf{\textit{BT}}^n=\textbf{\textit{BT}}_1\times\cdots\times\textbf{\textit{BT}}_n$
with
$\textbf{\textit{BT}}_1=\cdots=\textbf{\textit{BT}}_n=\textbf{\textit{BT}}$.

\vskip 0.2in

\bigskip
\centerline{\epsfxsize=2.7 in \epsfbox{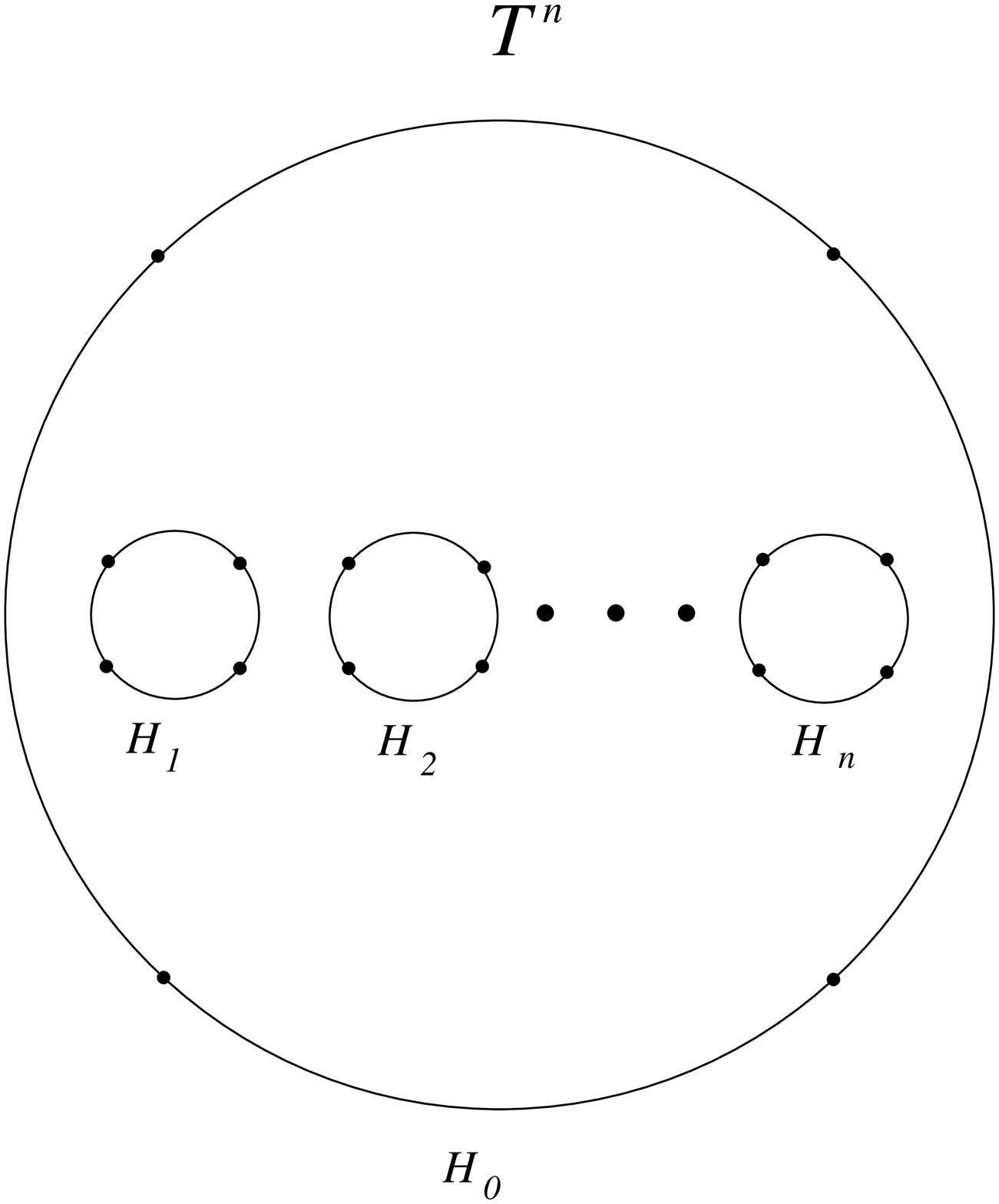}}
\medskip
\centerline{\small Figure 3. An $n$-punctured ball tangle with
$m_0=m_1=\cdots=m_n=2$.}
\bigskip

\bigskip
\centerline{\epsfxsize=4.7 in \epsfbox{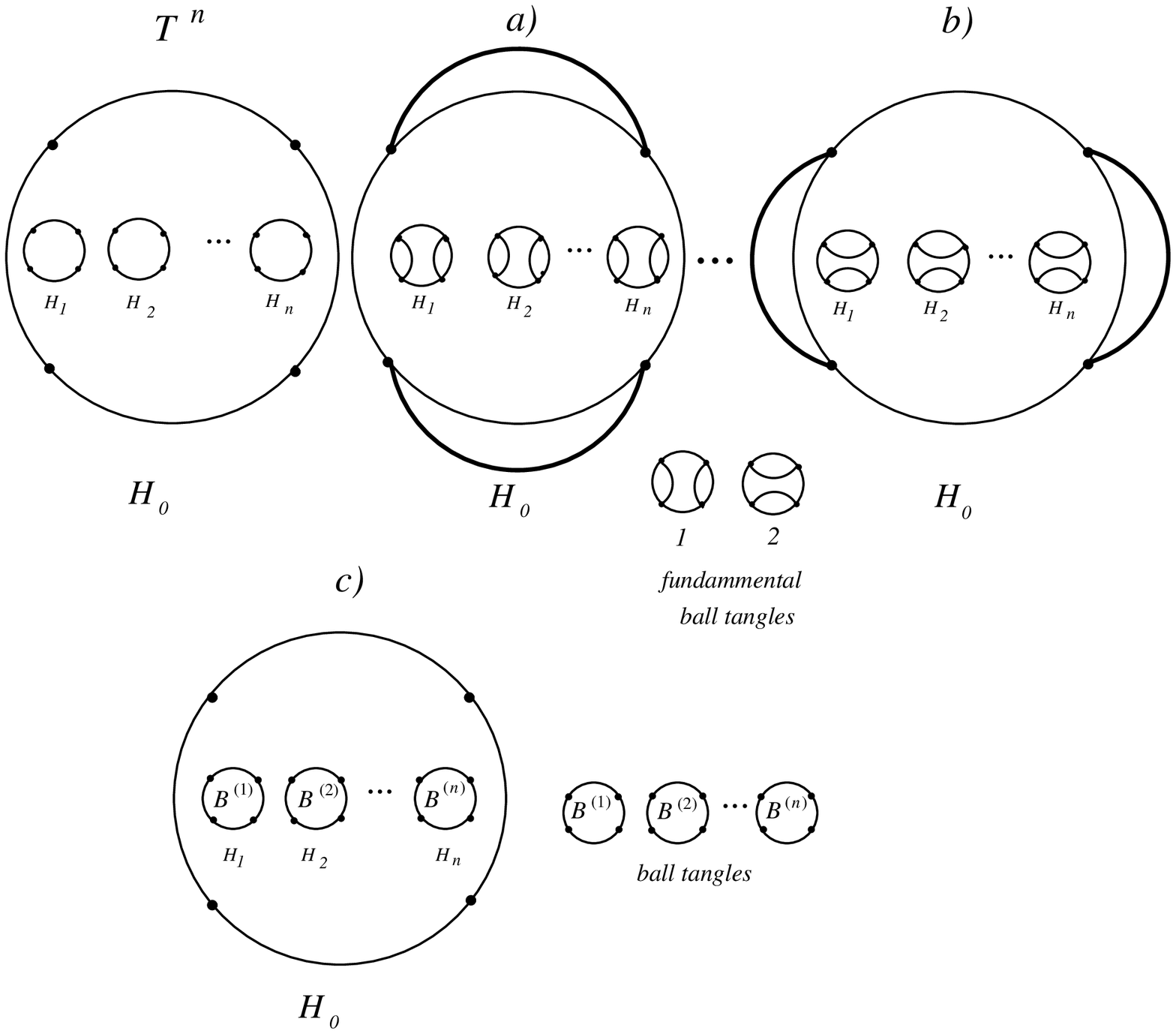}}
\medskip
\centerline{\small Figure 4. A hole-filling function $T^n$, a)
$T_{1\alpha_1^n}^n$, b) $T_{2\alpha_{2^n}^n}^n$, c)
$T^n(B^{(1)},\dots,B^{(n)})$.}
\bigskip

To construct our invariant of $n$-punctured ball tangles with
$m_0=m_1=\cdots=m_n=2$, we need to use some quite complicated
notations. Let us start with a gentle introduction to our
notations:

(1) For a diagram of 0-punctured ball tangle $T^0$ (a ball
tangle), we can produce 2 links $T_1^0$ and $T_2^0$, which are the
numerator closure and the denominator closure of $T^0$,
respectively.

(2) For a diagram of 1-punctured ball tangle $T^1$ (a spherical
tangle), we can produce $2^{1+1}$ links $T_{1(1)}^1$,
$T_{1(2)}^1$; $T_{2(1)}^1$, $T_{2(2)}^1$, where the subscript 1(1)
means to take the numerator closure of $T$ with its hole filled by
the fundamental tangle 1.

(3) For a diagram of 2-punctured ball tangle $T^2$, we can produce
$2^{2+1}$ links $T_{1(11)}^2$, $T_{1(12)}^2$, $T_{1(21)}^2$,
$T_{1(22)}^2$; $T_{2(11)}^2$, $T_{2(12)}^2$, $T_{2(21)}^2$,
$T_{2(22)}^2$.

If $n$ is a positive integer, $J_1=\cdots=J_n=\{1,2\}$, and
$J(n)=\prod_{k=1}^n J_k$, then $J(n)$ is linearly ordered by a
dictionary order, or lexicographic order, consisting of $2^n$
ordered $n$-tuples each of whose components is either 1 or 2. That
is, if $x,y \in J(n)$ and $x=(x_1,\dots,x_n)$,
$y=(y_1,\dots,y_n)$, then $x<y$ if and only if $x_1<y_1$ or there
is $k \in \{1,\dots,n-1\}$ such that
$x_1=y_1,\dots,x_k=y_k,x_{k+1}<y_{k+1}$.

(4) $J(n)=\{\alpha_i^n|1 \leq i \leq 2^n\}$ and
$\alpha_1^n<\alpha_2^n<\cdots<\alpha_{2^n}^n$, where $<$ is the
dictionary order on $J(n)$. Hence, $\alpha_1^n$ is the least
element $(1,1,\dots,1)$ and $\alpha_{2^n}^n$ is the greatest
element $(2,2,\dots,2)$ of $J(n)$. Let us denote
$\alpha_i^n=(\alpha_{i1}^n,\dots,\alpha_{in}^n)$ for each $i \in
\{1,\dots,2^n\}$.

(5) For a diagram of $n$-punctured ball tangle $T^n$, we can
produce $2^{n+1}$ links
$T_{1\alpha_1^n}^n,\dots,T_{1\alpha_{2^n}^n}^n$;
$T_{2\alpha_1^n}^n,\dots,T_{2\alpha_{2^n}^n}^n$.

(6) The sequence $(a_n)_{n \geq 0}=((t_k)_{1 \leq k \leq 2^n})_{n
\geq 0}$ is defined recursively as follows:

1) $a_0=(0)$;

2) If $a_{k-1}=(t_1,\dots,t_{2^{k-1}})$, then
$a_k=(t_1,\dots,t_{2^{k-1}},t_1+1,\dots,t_{2^{k-1}}+1)$ for each
$k \in N$. Note that $t_{2^n}=n$ for each $n \in \mathbb N \cup
\{0\}$.

Now, we define our invariant of $n$-punctured ball tangles with
$m_0=m_1=\cdots=m_n=2$ inductively.

\newpage

\begin{thm} For each $n \in \mathbb N$, define $F^n:\textbf{\textit{nPBT}} \rightarrow
PM_{2\times2^n}(\mathbb{Z})$ by
$$F^n(T^n)=\left\{\begin{pmatrix} (-i)^{t_1}z\langle T_{1\alpha_1^n}^n
\rangle & 
\cdots & (-i)^{t_{2^n}}z\langle T_{1\alpha_{2^n}^n}^n \rangle \\
(-i)^{t_1}iz\langle T_{2\alpha_1^n}^n \rangle &
\cdots & (-i)^{t_{2^n}}iz\langle T_{2\alpha_{2^n}^n}^n \rangle
\end{pmatrix}\,|\,z \in \Phi\right\}
\cap M_{2\times2^n}(\mathbb Z)$$ for each $T^n \in
\textbf{\textit{nPBT}}$. Then $F^n$ is an isotopy invariant of
$n$-punctured ball tangle diagrams. In particular, $F^0$ is
Krebes' ball tangle invariant $f$.
\end{thm}

\begin{proof} Let
$X(T^n)=\begin{pmatrix} (-i)^{t_1}\langle T_{1\alpha_1^n}^n
\rangle & (-i)^{t_2}\langle T_{1\alpha_2^n}^n \rangle & \cdots &
(-i)^{t_{2^n}}\langle T_{1\alpha_{2^n}^n}^n \rangle \\
(-i)^{t_1}i\langle T_{2\alpha_1^n}^n \rangle & (-i)^{t_2}i\langle
T_{2\alpha_2^n}^n \rangle & \cdots & (-i)^{t_{2^n}}i\langle
T_{2\alpha_{2^n}^n}^n \rangle
\end{pmatrix}$.
Then Tangle Reidemeister moves of type II, III, and IV do not
change $X(T^n)$ because Kauffman bracket is a regular invariant of
link diagrams. Also, it is easy to show that Tangle Reidemeister
move of type I does not change $\{zX(T^n)\,|\,z \in \Phi\}$ by the
skein relation of Kauffman bracket. Hence, it is enough to show
that $\{zX(T^n)\,|\,z \in \Phi\}\cap M_{2\times2^n}(\mathbb{Z})$
consists of two elements differ by a scalar multiplication of
$-1$. By Lemma 2.8, for each $k \in \{1,\dots,2^n\}$, there are
$p_k,q_k \in \mathbb{Z}$ and $u_k,v_k \in \Phi$ such that $\langle
T_{1\alpha_k^n}^n \rangle=p_ku_k$ and $\langle T_{2\alpha_k^n}^n
\rangle=q_kv_k$. Notice that, for each $k \in \{1,\dots,2^n\}$,
$\alpha_1^n$ and $\alpha_k^n$ differ at only $t_k$ coordinates.
Hence, $T_{1\alpha_1^n}^n$ and $T_{1\alpha_k^n}^n$ differ at only
$t_k$ holes. If $l,m \in \{1,\dots,2^n\}$ and $T_{1\alpha_l^n}^n$
and $T_{1\alpha_m^n}^n$ differ at only 1 hole, then $u_m=\pm iu_l$
by Proposition 2.6, Lemma 2.8, and Proposition 2.10. Thus,
$u_k=\pm i^{t_k}u_1$. Since $v_k=\pm iu_k$ for each $k \in
\{1,\dots,2^n\}$, $X(T^n) \in u_1M_{2\times2^n}(\mathbb{Z})$. This
shows that $F^n(T^n):=\{zX(T^n)\,|\,z \in \Phi\} \cap
M_{2\times2^n}(\mathbb{Z})=\{u_1^{-1}X(T^n),-u_1^{-1}X(T^n)\}$.
Therefore, $F^n(T_1)=F^n(T_2)$ if $T_1$ and $T_2$ are isotopic
$n$-punctured ball tangle diagrams.
\end{proof}

\bigskip
\centerline{\epsfxsize=5.3 in \epsfbox{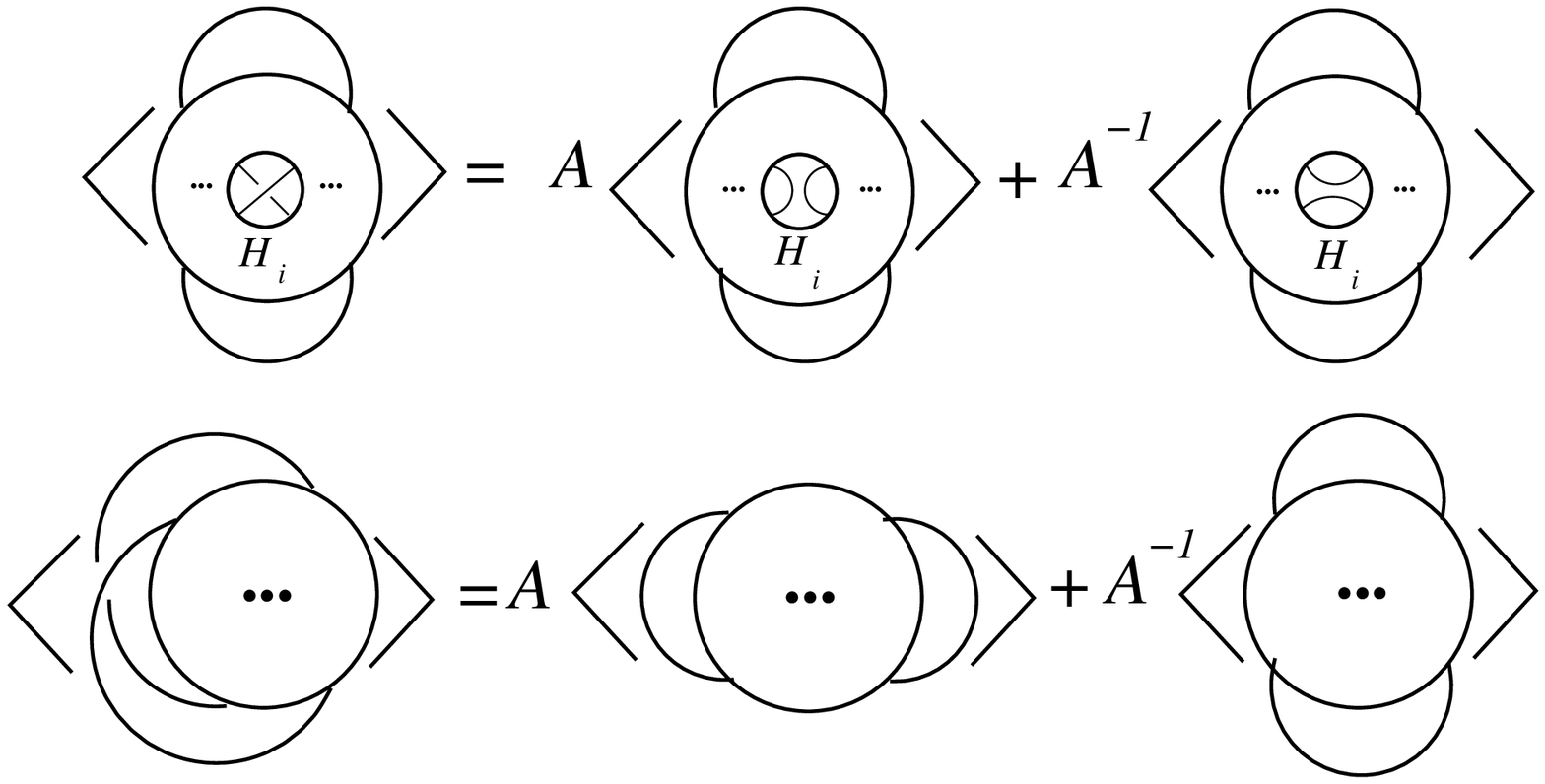}}
\medskip
\centerline{\small Figure 5. The skein relation of Kauffman
bracket at the $i$-th hole.}
\bigskip

\begin{df} For each nonnegative integer $n$, $F^n$ is called the
$n$-punctured ball tangle invariant, simply, the $n$-tangle
invariant.
\end{df}

Now, in order to think of an $n$-punctured ball tangle $T^n$ as a
`hole-filling function', we define a function which makes a
dictionary order on complex numbers.

Let $n$ be a positive integer, and let $(k_1,\dots,k_n)$ be an
$n$-tuple of positive integers, and let
$J(n,k_1,\dots,k_n)=\prod_{i=1}^n I_{k_i}$. Then
$J(n,k_1,\dots,k_n)$ is linearly ordered by a dictionary order,
where $I_k=\{1,\dots,k\}$ for each $k \in \mathbb{N}$.

(4$^*$) $J(n,k_1,\dots,k_n)= \{\alpha_i^{n,k_1,\dots,k_n}|1 \leq i
\leq k_1\cdots k_n\}$ and $\alpha_1^{n,k_1,\dots,k_n}<\cdots<
\alpha_{k_1\cdots k_n}^{n,k_1,\dots,k_n}$, where $<$ is the
dictionary order on $J(n,k_1,\dots,k_n)$. Hence,
$\alpha_1^{n,k_1,\dots,k_n}$ is the least element $(1,1,\dots,1)$
and $\alpha_{k_1\cdots k_n}^{n,k_1,\dots,k_n}$ is the greatest
element $(k_1,k_2,\dots,k_n)$ of $J(n,k_1,\dots,k_n)$. Let us
denote $\alpha_i^{n,k_1,\dots,k_n}=(\alpha_{i1}^{n,k_1,\dots,k_n},
\dots,\alpha_{in}^{n,k_1,\dots,k_n})$ for each $i \in
\{1,\dots,k_1\cdots k_n\}$.

\begin{df}
For each $n \in \mathbb{N}$ and $n$-tuple $(k_1,\dots,k_n)$ of
positive integers, define
$$\xi^{n,k_1,\dots,k_n}:\mathbb{C}^{k_1}\times\cdots\times
\mathbb{C}^{k_n} \rightarrow \mathbb{C}^{k_1\cdots k_n}$$ by
$$\xi^{n,k_1,\dots,k_n}((v_1^1,\dots,v_{k_1}^1),
\dots,(v_1^n,\dots,v_{k_n}^n))=\left(\prod_{j=1}^n
v_{\alpha_{1j}^{n,k_1,\dots,k_n}}^j, \dots, \prod_{j=1}^n
v_{\alpha_{k_1\cdots k_nj}^{n,k_1,\dots,k_n}}^j\right)$$ for all
$(v_1^1,\dots,v_{k_1}^1) \in
\mathbb{C}^{k_1},\dots,(v_1^n,\dots,v_{k_n}^n) \in
\mathbb{C}^{k_n}$. Then $\xi^{n,k_1,\dots,k_n}$ is well-defined
and called the dictionary order function on $\mathbb C$ with
respect to $k_1,\dots,k_n$. Also, the $i$-th projection of
$\xi^{n,k_1,\dots,k_n}$ is denoted by $\xi_i^{n,k_1,\dots,k_n}$
for each $i \in \{1,\dots,k_1\cdots k_n\}$. In particular, we
simply denote $\xi^{n,k_1,\dots,k_n}$ by $\xi^n$ when
$k_1=\cdots=k_n=2$.
\end{df}

Denote by $\mathbb{C}^{k\dag}$ the $k$-dimensional column vector
space over $\mathbb{C}$, so the map
$$(v_1,\dots,v_k)\mapsto(v_1,\dots,v_k)^\dag:\mathbb{C}^k\longrightarrow\mathbb{C}^{k\dag}$$
is to transpose row vectors to column vectors. Let
$P\mathbb{C}^{k\dag}=\mathbb{C}^{k\dag}/\pm1$. If
$(v_1,\dots,v_k)^\dag \in \mathbb{C}^{k\dag}$, then we denote by
$$[v_1,\dots,v_k]^\dag=\{(v_1,\dots,v_k)^\dag,(-v_1,\dots,-v_k)^\dag\}$$
the corresponding element in $P\mathbb{C}^{k\dag}$.

Remark that, we may extend the above notation to matrices modulo
$\pm1$. Under this extension, matrix multiplication is
well-defined. That is, if $A$ and $B$ are matrices and $AB$ is
defined, then $[A][B]=[A][-B]=[-A][B]=[-A][-B]=[-AB]=[AB]$.

\begin{lm}
For each $n \in\mathbb{N}$ and $n$-tuple $(k_1,\dots,k_n)$ of
positive integers, define
$$[\xi^{n,k_1,\dots,k_n}]:P\mathbb{C}^{k_1\dag}\times\cdots\times
P\mathbb{C}^{k_n\dag} \longrightarrow P\mathbb{C}^{k_1\cdots
k_n\dag}$$ by $$[\xi^{n,k_1,\dots,k_n}](\left[\begin{matrix} v_1^1
\\ \cdot \\ \cdot \\ \cdot \\ v_{k_1}^1 \end{matrix}\right], \dots,
\left[\begin{matrix} v_1^n \\ \cdot \\ \cdot \\ \cdot \\ v_{k_n}^n
\end{matrix}\right])=\left[\begin{matrix}
\prod_{j=1}^n v_{\alpha_{1j}^{n,k_1,\dots,k_n}}^j
\\ \cdot \\ \cdot \\ \cdot \\
\prod_{j=1}^n v_{\alpha_{k_1\cdots k_nj}^{n,k_1,\dots,k_n}}^j
\end{matrix}\right]$$ for all $(v_1^1,\dots,v_{k_1}^1) \in
\mathbb{C}^{k_1},\dots,(v_1^n,\dots,v_{k_n}^n) \in
\mathbb{C}^{k_n}$. Then $[\xi^{n,k_1,\dots,k_n}]$ is well-defined
and called the dictionary order function induced by
$\xi^{n,k_1,\dots,k_n}$.
\end{lm}

\begin{proof} Suppose that $(X_1,\dots,X_n)$ and
$(Y_1,\dots,Y_n)$ are in $\mathbb{C}^{k_1}\times\cdots\times
\mathbb{C}^{k_n}$ such that
$([X_1]^\dag,\dots,[X_n]^\dag)=([Y_1]^\dag,\dots,[Y_n]^\dag) \in
P\mathbb{C}^{k_1\dag}\times\cdots\times P\mathbb{C}^{k_n\dag}$.
Then $X_1=\pm Y_1,\dots,X_n=\pm Y_n$ and
$\xi^{n,k_1,\dots,k_n}(X_1,\dots,X_n)=
\pm\xi^{n,k_1,\dots,k_n}(Y_1,\dots,Y_n)$. Hence,
$$[\xi^{n,k_1,\dots,k_n}(X_1,\dots,X_n)]^\dag=
[\xi^{n,k_1,\dots,k_n}(Y_1,\dots,Y_n)]^\dag.$$

Therefore, $[\xi^{n,k_1,\dots,k_n}]([X_1]^\dag,\dots,[X_n]^\dag)=
[\xi^{n,k_1,\dots,k_n}]([Y_1]^\dag,\dots,[Y_n]^\dag)$. This shows
that $[\xi^{n,k_1,\dots,k_n}]$ is well-defined.
\end{proof}

As another notation, if $L$ is a link diagram and $T^n$ is a
diagram of $n$-punctured ball tangle for some $n \in\mathbb{ N}
\cup \{0\}$, then the sets of all crossings of $L$ and $T^n$ are
denoted by $c(L)$ and $c(T^n)$, respectively.

\begin{lm} If $n \in \mathbb{N}$ and $T^n$ is an $n$-punctured ball tangle
diagram and $B^{(1)},\dots,B^{(n)}$ are ball tangle diagrams, then
$$\begin{pmatrix} \langle T^n(B^{(1)},\dots,B^{(n)})_1 \rangle \\
\langle T^n(B^{(1)},\dots,B^{(n)})_2 \rangle
\end{pmatrix}=\begin{pmatrix} \sum_{i=1}^{2^n}
\langle T_{1\alpha_i^n}^n \rangle \langle B_{\alpha_{i1}^n}^{(1)}
\rangle \cdots \langle B_{\alpha_{in}^n}^{(n)} \rangle \\
\sum_{i=1}^{2^n} \langle T_{2\alpha_i^n}^n \rangle \langle
B_{\alpha_{i1}^n}^{(1)} \rangle \cdots \langle
B_{\alpha_{in}^n}^{(n)} \rangle
\end{pmatrix}.$$
\end{lm}

\begin{proof} We denote the set of all monocyclic states of
a link diagram $L$ by $M(L)$. Let $B=T^n(B^{(1)},\dots,B^{(n)})$.
Then $\sigma$ is a monocyclic state of $B_1$ if and only if there
is a unique $i \in \{1,\dots,2^n\}$ such that $\sigma|_{c(T^n)}
\in M(T_{1\alpha_i^n}^n)$, $\sigma|_{c(B^{(1)})} \in
M(B_{\alpha_{i1}^n}^{(1)}),\dots,\sigma|_{c(B^{(n)})} \in
M(B_{\alpha_{in}^n}^{(n)})$. Note that $c(B_1)=c(T^n) \coprod
c(B^{(1)}) \coprod \cdots \coprod c(B^{(n)})$. Let us denote a
state $\sigma$ of $B_1$ by $\sigma_0\sigma_1\cdots\sigma_n$, where
$\sigma_0=\sigma|_{c(T^n)},
\sigma_1=\sigma|_{c(B^{(1)})},\dots,\sigma_n=\sigma|_{c(B^{(n)})}$.
Then
$$M(B_1)=\coprod_{i=1}^{2^n}
M(T_{1\alpha_i^n}^n)M(B_{\alpha_{i1}^n}^{(1)})\cdots
M(B_{\alpha_{in}^n}^{(n)})$$ and
$$\begin{aligned}
&\sum_{\sigma \in M(B_1)}A^{\alpha(\sigma)-\beta(\sigma)}\\
&=\sum_{i=1}^{2^n} \sum_{\sigma_0\sigma_1\cdots\sigma_n \in
M(T_{1\alpha_i^n}^n)M(B_{\alpha_{i1}^n}^{(1)})\cdots
M(B_{\alpha_{in}^n}^{(n)})}
A^{\alpha(\sigma_0)+\alpha(\sigma_1)+\cdots+
\alpha(\sigma_n)-\beta(\sigma_0)-\beta(\sigma_1)-\cdots-\beta(\sigma_n)}\\
&=\sum_{i=1}^{2^n} \sum_{\sigma_0 \in M(T_{1\alpha_i^n}^n)}
A^{\alpha(\sigma_0)-\beta(\sigma_0)} \sum_{\sigma_1 \in
M(B_{\alpha_{i1}^n}^{(1)})} A^{\alpha(\sigma_1)-\beta(\sigma_1)}
\cdots \sum_{\sigma_n \in M(B_{\alpha_{in}^n}^{(n)})}
A^{\alpha(\sigma_n)-\beta(\sigma_n)}\\
&=\sum_{i=1}^{2^n} \langle T_{1\alpha_i^n}^n \rangle \langle
B_{\alpha_{i1}^n}^{(1)} \rangle \cdots \langle
B_{\alpha_{in}^n}^{(n)} \rangle
\end{aligned}
$$

Similarly, $\langle B_2 \rangle=\sum_{i=1}^{2^n} \langle
T_{2\alpha_i^n}^n \rangle \langle B_{\alpha_{i1}^n}^{(1)} \rangle
\cdots \langle B_{\alpha_{in}^n}^{(n)} \rangle$. This proves the
lemma.
\end{proof}

\begin{thm} For each $n \in \mathbb{N}$, $F^n$ is an $n$-punctured ball tangle
invariant such that
$F^0(T^n(B^{(1)},\dots,B^{(n)}))=F^n(T^n)[\xi^n](F^0(B^{(1)}),
\dots,F^0(B^{(n)}))$ for all $B^{(1)},\dots,B^{(n)} \in
\textbf{\textit{BT}}$.
\end{thm}

\begin{proof} Suppose that $T^n$ is an $n$-punctured ball tangle
such that $F^n(T^n)=[zX(T^n)]$ for some $z \in \Phi$ and
$B^{(1)},\dots,B^{(n)}$ are ball tangles such that
$$F^0(B^{(1)})=\left[\begin{matrix} z_1\langle B_1^{(1)} \rangle \\
iz_1\langle B_2^{(1)} \rangle \end{matrix}\right],\dots,
F^0(B^{(n)})=\left[\begin{matrix} z_n\langle B_1^{(n)} \rangle \\
iz_n\langle B_2^{(n)} \rangle \end{matrix}\right]$$ for some
$z_1,\dots,z_n \in \Phi$, where $\langle B_1^{(i)} \rangle$ and
$\langle B_2^{(i)} \rangle$ are the numerator closure and the
denominator closure of $B^{(i)}$, respectively, for each $i \in
\{1,\dots,n\}$. Then
$$
\begin{aligned}
&F^n(T^n)[\xi^n](F^0(B^{(1)}),\dots,F^0(B^{(n)}))\\
&=\left[\begin{matrix} \begin{pmatrix} (-i)^{t_1}z\langle
T_{1\alpha_1^n}^n \rangle & \cdots &
(-i)^{t_{2^n}}z\langle T_{1\alpha_{2^n}^n}^n \rangle \\
(-i)^{t_1}iz\langle T_{2\alpha_1^n}^n \rangle & \cdots &
(-i)^{t_{2^n}}iz\langle T_{2\alpha_{2^n}^n}^n \rangle
\end{pmatrix}
\begin{pmatrix} i^{t_1}z_1\cdots z_n\langle B_{\alpha_{11}^n}^{(1)} \rangle
\cdots \langle B_{\alpha_{1n}^n}^{(n)} \rangle \\ \cdot \\ \cdot \\ \cdot \\
i^{t_{2^n}}z_1\cdots z_n\langle B_{\alpha_{2^n1}^n}^{(1)} \rangle
\cdots \langle B_{\alpha_{2^nn}^n}^{(n)} \rangle
\end{pmatrix} \end{matrix}\right]
\end{aligned}$$
$$
\begin{aligned}
&=\left[\begin{matrix} zz_1\cdots z_n(\langle T_{1\alpha_1^n}^n
\rangle \langle B_{\alpha_{11}^n}^{(1)} \rangle \cdots \langle
B_{\alpha_{1n}^n}^{(n)} \rangle + \cdots\cdots + \langle
T_{1\alpha_{2^n}^n}^n \rangle \langle B_{\alpha_{2^n1}^n}^{(1)}
\rangle \cdots \langle B_{\alpha_{2^nn}^n}^{(n)} \rangle) \\
izz_1\cdots z_n(\langle T_{2\alpha_1^n}^n \rangle \langle
B_{\alpha_{11}^n}^{(1)} \rangle \cdots \langle
B_{\alpha_{1n}^n}^{(n)} \rangle + \cdots\cdots + \langle
T_{2\alpha_{2^n}^n}^n \rangle \langle B_{\alpha_{2^n1}^n}^{(1)}
\rangle \cdots \langle B_{\alpha_{2^nn}^n}^{(n)} \rangle)
\end{matrix}\right]\\
&=\left[\begin{matrix} zz_1\cdots z_n \sum_{i=1}^{2^n} \langle
T_{1\alpha_i^n}^n \rangle \langle B_{\alpha_{i1}^n}^{(1)}
\rangle \cdots \langle B_{\alpha_{in}^n}^{(n)} \rangle \\
izz_1\cdots z_n \sum_{i=1}^{2^n} \langle T_{2\alpha_i^n}^n \rangle
\langle B_{\alpha_{i1}^n}^{(1)} \rangle \cdots \langle
B_{\alpha_{in}^n}^{(n)} \rangle
\end{matrix}\right]=\left[\begin{matrix} zz_1\cdots z_n
\langle T^n(B^{(1)},\dots,B^{(n)})_1 \rangle \\
izz_1\cdots z_n\langle T^n(B^{(1)},\dots,B^{(n)})_2 \rangle
\end{matrix}\right]\\
&=F^0(T^n(B^{(1)},\dots,B^{(n)}))
\end{aligned}$$
by Lemma 3.5.
\end{proof}

Notice that an $n$-punctured ball tangle $T^n$ can be regarded as
an $n$ variable function about not only $0$-punctured ball tangles
but also various $n$-punctured ball tangles. Given an
$n$-punctured ball tangle diagram $T^n$ and
$k_1,\dots,k_n$-punctured ball tangle diagrams
$T^{k_1(1)},\dots,T^{k_n(n)}$, respectively, we consider the $k_1
+ \cdots + k_n$-punctured ball tangle diagram
$T^n(T^{k_1(1)},\dots,T^{k_n(n)})$, where $n \in \mathbb N$ and
$k_1,\dots,k_n \in \mathbb N \cup \{0\}$. We show how to calculate
the invariant
$F^{k_1+\cdots+k_n}(T^n(T^{k_1(1)},\dots,T^{k_n(n)}))$ of it if
$F^n(T^n), F^{k_1}(T^{k_1(1)}),\dots,F^{k_n}(T^{k_n(n)})$ are
given (Theorem 3.8). Let us start from the following notations:

Let $n \in \mathbb N$. Then

(1) $e_i^n=\left[\begin{matrix} v_1 \\ \cdot \\ \cdot \\ \cdot \\
v_{2^n} \end{matrix}\right]$ such that $v_i=1$ and $v_j=0$ if $j
\neq i$ for each $i \in \{1,\dots,2^n\}$. In particular,
$e_1^1=\left[\begin{matrix} 1 \\ 0
\end{matrix}\right]$ and $e_2^1=\left[\begin{matrix} 0 \\ 1
\end{matrix}\right]$. Hence, $e_i^n=[\xi^n]
(e_{\alpha_{i1}^n}^1,\dots,e_{\alpha_{in}^n}^1)$ for each $i \in
\{1,\dots,2^n\}$.

(2) $x=\left[\begin{matrix} 1 \\ 1
\end{matrix}\right]$.

(3) $E_j^n$ is the set of all $[\xi^n](y_1,\dots,y_n)$ such that
$j$ components of $(y_1,\dots,y_n)$ are $x$ and each of the others
is $e_1^1$ or $e_2^1$ for each $j \in \{0,1,\dots,n\}$. In
particular,
$$E_0^n=\{[\xi^n](e_{\alpha_{i1}^n}^1,\dots,e_{\alpha_{in}^n}^1)|i
\in \{1,\dots,2^n\}\}$$ and
$$E_n^n=\{[\xi^n](y_1,\dots,y_n)|y_1=\cdots=y_n=x\}.$$
Notice that $\{E_0^n,E_1^n,\dots,E_n^n\}$ is pairwise disjoint and
$|E_j^n|=\,_nC_j\,2^{n-j}$ for each $j \in \{0,1,\dots,n\}$, where
$\,_nC_j\,=\frac{n!}{(n-j)!j\,!}$. Hence,
$$|\coprod_{j=0}^n E_j^n|=\,_nC_0\,2^n + \,_nC_1\,2^{n-1} + \cdots
+ \,_nC_{n-1}\,2^1 + \,_nC_n\,2^0=(2+1)^n=3^n.$$ Note that
$$\coprod_{j=0}^n E_j^n=[\xi^n](\{e_1^1,e_2^1,x\}^n).$$

For example, when $n=3$, we have
$$E_0^3=\{e_1^3,e_2^3,e_3^3,e_4^3,e_5^3,e_6^3,e_7^3,e_8^3\},$$
$$E_1^3=\{\,[\,1\,0\,0\,0\,1\,0\,0\,0\,]^\dag,[\,0\,1\,0\,0\,0\,1\,0\,0\,]^\dag,
[\,0\,0\,1\,0\,0\,0\,1\,0\,]^\dag,[\,0\,0\,0\,1\,0\,0\,0\,1\,]^\dag,$$
$$[\,1\,0\,1\,0\,0\,0\,0\,0\,]^\dag,[\,0\,1\,0\,1\,0\,0\,0\,0\,]^\dag,
[\,0\,0\,0\,0\,1\,0\,1\,0\,]^\dag,[\,0\,0\,0\,0\,0\,1\,0\,1\,]^\dag,$$
$$[\,1\,1\,0\,0\,0\,0\,0\,0\,]^\dag,[\,0\,0\,1\,1\,0\,0\,0\,0\,]^\dag,
[\,0\,0\,0\,0\,1\,1\,0\,0\,]^\dag,[\,0\,0\,0\,0\,0\,0\,1\,1\,]^\dag\},$$
$$E_2^3=\{\,[\,1\,0\,1\,0\,1\,0\,1\,0\,]^\dag,[\,0\,1\,0\,1\,0\,1\,0\,1\,]^\dag,$$
$$[\,1\,1\,0\,0\,1\,1\,0\,0\,]^\dag,[\,0\,0\,1\,1\,0\,0\,1\,1\,]^\dag,$$
$$[\,1\,1\,1\,1\,0\,0\,0\,0\,]^\dag,[\,0\,0\,0\,0\,1\,1\,1\,1\,]^\dag\},$$
$$E_3^3=\{\,[\,1\,1\,1\,1\,1\,1\,1\,1\,]^\dag\}.$$

Now, we have the following lemma which supports our next theorems.

\begin{lm} If $n \in \mathbb N$ and $A,B \in PM_{2\times2^n}(\mathbb Z)$ and
$AX=BX$ for each $X \in [\xi^n](\{e_1^1,e_2^1,x\}^n)$, then $A=B$.
\end{lm}

\begin{proof} We prove the statement by induction on $n \in \mathbb N$.

Step 1. We show that the statement is true for $n=1$.

Let $A=\left[\begin{matrix} a_{11} & a_{12} \\ a_{21} & a_{22}
\end{matrix}\right]$, and let $B=\left[\begin{matrix}
b_{11} & b_{12} \\ b_{21} & b_{22} \end{matrix}\right]$. Since
$E_0^1=\{\left[\begin{matrix} 1 \\ 0
\end{matrix}\right],\left[\begin{matrix} 0 \\ 1
\end{matrix}\right]\}$ and $E_1^1=\{\left[\begin{matrix} 1 \\ 1
\end{matrix}\right]\}$ and $AX=BX$ for each $X \in [\xi^1](\{e_1^1,e_2^1,x\}^1)$,
$\left[\begin{matrix} a_{11} \\ a_{21}
\end{matrix}\right]=\left[\begin{matrix} b_{11} \\ b_{21}
\end{matrix}\right]$ and $\left[\begin{matrix} a_{12} \\ a_{22}
\end{matrix}\right]=\left[\begin{matrix} b_{12} \\ b_{22}
\end{matrix}\right]$ and $\left[\begin{matrix} a_{11}+a_{12} \\
a_{21}+a_{22} \end{matrix}\right]=\left[\begin{matrix} b_{11}+b_{12} \\
b_{21}+b_{22} \end{matrix}\right]$. Hence, $\begin{pmatrix}
a_{11} \\ a_{21} \end{pmatrix}=\epsilon\begin{pmatrix} b_{11} \\
b_{21} \end{pmatrix}$ and $\begin{pmatrix}
a_{12} \\ a_{22} \end{pmatrix}=\epsilon'\begin{pmatrix} b_{12} \\
b_{22} \end{pmatrix}$ and $\begin{pmatrix} a_{11}+a_{12} \\
a_{21}+a_{22} \end{pmatrix}= \epsilon_1\begin{pmatrix}
b_{11}+b_{12} \\ b_{21}+b_{22} \end{pmatrix}$ for some
$\epsilon,\epsilon',\epsilon_1 \in \{1,-1\}$. Suppose that
$\epsilon\epsilon'=-1$.

{\it Case 1}. $\begin{pmatrix} a_{11} \\ a_{21}
\end{pmatrix}=\begin{pmatrix} b_{11} \\ b_{21}
\end{pmatrix}$ and $\begin{pmatrix} a_{12} \\ a_{22}
\end{pmatrix}=\begin{pmatrix} -b_{12} \\ -b_{22} \end{pmatrix}$.
If $\epsilon_1=1$, then
$a_{11}+a_{12}=b_{11}+b_{12}=b_{11}-b_{12}$ and
$a_{21}+a_{22}=b_{21}+b_{22}=b_{21}-b_{22}$, so $b_{12}=b_{22}=0$.
If $\epsilon_1=-1$, then
$a_{11}+a_{12}=-b_{11}-b_{12}=b_{11}-b_{12}$ and
$a_{21}+a_{22}=-b_{21}-b_{22}=b_{21}-b_{22}$, so
$b_{11}=b_{21}=0$. Hence, $A=B$.

{\it Case 2}. $\begin{pmatrix} a_{11} \\ a_{21}
\end{pmatrix}=\begin{pmatrix} -b_{11} \\ -b_{21}
\end{pmatrix}$ and $\begin{pmatrix} a_{12} \\ a_{22}
\end{pmatrix}=\begin{pmatrix} b_{12} \\ b_{22} \end{pmatrix}$.
If $\epsilon_1=1$, then
$a_{11}+a_{12}=b_{11}+b_{12}=-b_{11}+b_{12}$ and
$a_{21}+a_{22}=b_{21}+b_{22}=-b_{21}+b_{22}$, so
$b_{11}=b_{21}=0$. If $\epsilon_1=-1$, then
$a_{11}+a_{12}=-b_{11}-b_{12}=-b_{11}+b_{12}$ and
$a_{21}+a_{22}=-b_{21}-b_{22}=-b_{21}+b_{22}$, so
$b_{12}=b_{22}=0$. Hence, $A=B$.

Step 2. Suppose that the statement is true for $n \in \mathbb N$.
We show that the statement is also true for $n+1$.

Suppose that $$A=\left[\begin{matrix} a_{11} & \cdots & a_{12^n} &
a_{12^n+1} & \cdots & a_{12^{n+1}} \\ a_{21} & \cdots & a_{22^n} &
a_{22^n+1} & \cdots & a_{22^{n+1}} \end{matrix}\right]$$ and
$$B=\left[\begin{matrix} b_{11} & \cdots & b_{12^n} & b_{12^n+1} &
\cdots & b_{12^{n+1}} \\ b_{21} & \cdots & b_{22^n} & b_{22^n+1} &
\cdots & b_{22^{n+1}} \end{matrix}\right]$$ and
$$A_1=\begin{pmatrix} a_{11} & \cdots & a_{12^n} \\
a_{21} & \cdots & a_{22^n} \end{pmatrix},\,
A_2=\begin{pmatrix} a_{12^n+1} & \cdots & a_{12^{n+1}} \\
a_{22^n+1} & \cdots & a_{22^{n+1}} \end{pmatrix},$$
$$A_3=\begin{pmatrix}
a_{11} & \cdots & a_{12^{n-1}} & a_{12^n+1} & \cdots &
a_{12^n+2^{n-1}} \\ a_{21} & \cdots & a_{22^{n-1}} & a_{22^n+1} &
\cdots & a_{22^n+2^{n-1}} \end{pmatrix},$$
$$A_4=\begin{pmatrix}
a_{12^{n-1}+1} & \cdots & a_{12^n} & a_{12^n+2^{n-1}+1} & \cdots &
a_{12^{n+1}} \\ a_{22^{n-1}+1} & \cdots & a_{22^n} &
a_{22^n+2^{n-1}+1} & \cdots & a_{22^{n+1}}
\end{pmatrix}$$ and $$B_1=\begin{pmatrix} b_{11} & \cdots & b_{12^n} \\
b_{21} & \cdots & b_{22^n} \end{pmatrix},\,
B_2=\begin{pmatrix} b_{12^n+1} & \cdots & b_{12^{n+1}} \\
b_{22^n+1} & \cdots & b_{22^{n+1}} \end{pmatrix},$$
$$B_3=\begin{pmatrix}
b_{11} & \cdots & b_{12^{n-1}} & b_{12^n+1} & \cdots &
b_{12^n+2^{n-1}} \\ b_{21} & \cdots & b_{22^{n-1}} & b_{22^n+1} &
\cdots & b_{22^n+2^{n-1}} \end{pmatrix},$$
$$B_4=\begin{pmatrix}
b_{12^{n-1}+1} & \cdots & b_{12^n} & b_{12^n+2^{n-1}+1} & \cdots &
b_{12^{n+1}} \\ b_{22^{n-1}+1} & \cdots & b_{22^n} &
b_{22^n+2^{n-1}+1} & \cdots & b_{22^{n+1}}
\end{pmatrix}.$$

Then $A=\left[\begin{matrix} A_1 & A_2
\end{matrix}\right]$ and $B=\left[\begin{matrix} B_1 & B_2
\end{matrix}\right]$. Notice that $$[\xi^{n+1}](\{e_1^1,e_2^1,x\}^{n+1})=$$
$$[\xi^{n+1}](\{e_1^1\} \times \{e_1^1,e_2^1,x\}^n) \coprod
[\xi^{n+1}](\{e_2^1\} \times \{e_1^1,e_2^1,x\}^n) \coprod
[\xi^{n+1}](\{x\} \times \{e_1^1,e_2^1,x\}^n)$$ and
$[\xi^{n+1}](\{e_1^1\} \times \{e_1^1,e_2^1,x\}^n)$,
$[\xi^{n+1}](\{e_2^1\} \times \{e_1^1,e_2^1,x\}^n)$,
$[\xi^{n+1}](\{x\} \times \{e_1^1,e_2^1,x\}^n)$ have exactly $3^n$
elements, respectively.

Since $AX=BX$ for each $X \in [\xi^{n+1}](\{e_1^1\} \times
\{e_1^1,e_2^1,x\}^n)$, $[A_1]X=[B_1]X$ for each $X \in
[\xi^n](\{e_1^1,e_2^1,x\}^n)$.

Similarly, since $AX=BX$ for each $X \in [\xi^{n+1}](\{e_2^1\}
\times \{e_1^1,e_2^1,x\}^n)$, $[A_2]X=[B_2]X$ for each $X \in
[\xi^n](\{e_1^1,e_2^1,x\}^n)$.

Also, since $AX=BX$ for each $X \in [\xi^{n+1}](\{e_1^1,e_2^1,x\}
\times \{e_1^1\} \times \{e_1^1,e_2^1,x\}^{n-1})$, $[A_3]X=[B_3]X$
for each $X \in [\xi^n](\{e_1^1,e_2^1,x\}^n)$.

Similarly, since $AX=BX$ for each $X \in
[\xi^{n+1}](\{e_1^1,e_2^1,x\} \times \{e_2^1\} \times
\{e_1^1,e_2^1,x\}^{n-1})$, $[A_4]X=[B_4]X$ for each $X \in
[\xi^n](\{e_1^1,e_2^1,x\}^n)$.

By induction hypothesis, we have $$[A_1]=[B_1],\,\,
[A_2]=[B_2],\,\, [A_3]=[B_3],\,\, [A_4]=[B_4].$$ Hence, $A_1 =
\epsilon B_1$ and $A_2 = \epsilon' B_2$ for some
$\epsilon,\epsilon' \in \{1,-1\}$. Now, we claim that, if
$\epsilon\epsilon'=-1$, then $B_1$ or $B_2$ is the $2\times2^n$
zero matrix.

Suppose that $\epsilon\epsilon'=-1$. Without loss of generality,
we may assume that $$A_1=B_1\,\, {\rm and}\,\, A_2=-B_2.$$

Suppose that $B_1$ is not the $2\times2^n$ zero matrix. Then there
is $i \in \{1,\dots,2^n\}$ such that $$\begin{pmatrix} b_{1i} \\
b_{2i} \end{pmatrix} \neq \begin{pmatrix} 0 \\ 0 \end{pmatrix}.$$

{\it Case 1}. If $1 \leq i \leq 2^{n-1}$, then
$\begin{pmatrix} b_{12^n+1} & \cdots & b_{12^n+2^{n-1}} \\
b_{22^n+1} & \cdots & b_{22^n+2^{n-1}} \end{pmatrix}$ is the
$2\times2^{n-1}$ zero matrix since $[A_3]=[B_3]$. We claim that
$\begin{pmatrix} b_{12^n+2^{n-1}+1} & \cdots & b_{12^{n+1}} \\
b_{22^n+2^{n-1}+1} & \cdots & b_{22^{n+1}} \end{pmatrix}$ is also
the $2\times2^{n-1}$ zero matrix.

Suppose that $\begin{pmatrix} b_{12^n+2^{n-1}+1} & \cdots & b_{12^{n+1}} \\
b_{22^n+2^{n-1}+1} & \cdots & b_{22^{n+1}} \end{pmatrix}$ is not
the $2\times2^{n-1}$ zero matrix. Then there is $j \in
\{2^n+2^{n-1}+1,\dots,2^{n+1}\}$ such that $\begin{pmatrix} b_{1j}
\\ b_{2j} \end{pmatrix} \neq \begin{pmatrix} 0 \\ 0 \end{pmatrix}$.

Since $[A_4]=[B_4]$, $\begin{pmatrix} b_{12^{n-1}+1} & \cdots &
b_{12^n} \\ b_{22^{n-1}+1} & \cdots & b_{22^n} \end{pmatrix}$ is
the $2\times2^{n-1}$ zero matrix. In this case, the fact that
$AX=BX$ for each $X \in [\xi^{n+1}](\{e_1^1,e_2^1,x\}^{n+1})$
implies
$$\left[\begin{matrix} a_{11} & \cdots & a_{12^{n-1}} &
a_{12^n+2^{n-1}+1} & \cdots & a_{12^{n+1}} \\ a_{21} & \cdots &
a_{22^{n-1}} & a_{22^n+2^{n-1}+1} & \cdots & a_{22^{n+1}}
\end{matrix}\right]X$$
$$=\left[\begin{matrix} b_{11} & \cdots & b_{12^{n-1}} &
b_{12^n+2^{n-1}+1} & \cdots & b_{12^{n+1}} \\ b_{21} & \cdots &
b_{22^{n-1}} & b_{22^n+2^{n-1}+1} & \cdots & b_{22^{n+1}}
\end{matrix}\right]X$$ for each $X \in
[\xi^n](\{e_1^1,e_2^1,x\}^n)$. Hence, by induction hypothesis, we
have $$\left[\begin{matrix} a_{11} & \cdots & a_{12^{n-1}} &
a_{12^n+2^{n-1}+1} & \cdots & a_{12^{n+1}} \\ a_{21} & \cdots &
a_{22^{n-1}} & a_{22^n+2^{n-1}+1} & \cdots & a_{22^{n+1}}
\end{matrix}\right]$$
$$=\left[\begin{matrix} b_{11} & \cdots & b_{12^{n-1}} &
b_{12^n+2^{n-1}+1} & \cdots & b_{12^{n+1}} \\ b_{21} & \cdots &
b_{22^{n-1}} & b_{22^n+2^{n-1}+1} & \cdots & b_{22^{n+1}}
\end{matrix}\right].$$

Since $A_1=B_1$ and $A_2=-B_2$,
$$\left[\begin{matrix} b_{11} & \cdots & b_{12^{n-1}} &
-b_{12^n+2^{n-1}+1} & \cdots & -b_{12^{n+1}} \\ b_{21} & \cdots &
b_{22^{n-1}} & -b_{22^n+2^{n-1}+1} & \cdots & -b_{22^{n+1}}
\end{matrix}\right]$$
$$=\left[\begin{matrix} b_{11} & \cdots & b_{12^{n-1}} &
b_{12^n+2^{n-1}+1} & \cdots & b_{12^{n+1}} \\ b_{21} & \cdots &
b_{22^{n-1}} & b_{22^n+2^{n-1}+1} & \cdots & b_{22^{n+1}}
\end{matrix}\right].$$

Since $\begin{pmatrix} b_{1i} \\ b_{2i} \end{pmatrix} \neq
\begin{pmatrix} 0 \\ 0 \end{pmatrix}$,
$\begin{pmatrix} b_{12^n+2^{n-1}+1} & \cdots & b_{12^{n+1}} \\
b_{22^n+2^{n-1}+1} & \cdots & b_{22^{n+1}} \end{pmatrix}$ is the
$2\times2^{n-1}$ zero matrix. This is a contradiction. Therefore,
$B_2$ is the $2\times2^n$ zero matrix.

Similarly, we show the other case.

{\it Case 2}. If $2^{n-1}+1 \leq i \leq 2^n$, then
$\begin{pmatrix} b_{12^n+2^{n-1}+1} & \cdots & b_{12^{n+1}} \\
b_{22^n+2^{n-1}+1} & \cdots & b_{22^{n+1}} \end{pmatrix}$ is the
$2\times2^{n-1}$ zero matrix since $[A_4]=[B_4]$. We claim that
$\begin{pmatrix} b_{12^n+1} & \cdots & b_{12^n+2^{n-1}} \\
b_{22^n+1} & \cdots & b_{22^n+2^{n-1}} \end{pmatrix}$ is also the
$2\times2^{n-1}$ zero matrix.

Suppose that $\begin{pmatrix} b_{12^n+1} & \cdots & b_{12^n+2^{n-1}} \\
b_{22^n+1} & \cdots & b_{22^n+2^{n-1}} \end{pmatrix}$ is not the
$2\times2^{n-1}$ zero matrix. Then there is $j \in
\{2^n+1,\dots,2^n+2^{n-1}\}$ such that $\begin{pmatrix} b_{1j}
\\ b_{2j} \end{pmatrix} \neq \begin{pmatrix} 0 \\ 0 \end{pmatrix}$.

Since $[A_3]=[B_3]$, $\begin{pmatrix} b_{11} & \cdots &
b_{12^{n-1}} \\ b_{21} & \cdots & b_{22^{n-1}} \end{pmatrix}$ is
the $2\times2^{n-1}$ zero matrix. In this case, the fact that
$AX=BX$ for each $X \in [\xi^{n+1}](\{e_1^1,e_2^1,x\}^{n+1})$
implies
$$\left[\begin{matrix} a_{12^{n-1}+1} & \cdots & a_{12^n} &
a_{12^n+1} & \cdots & a_{12^n+2^{n-1}} \\ a_{22^{n-1}+1} & \cdots
& a_{22^n} & a_{22^n+1} & \cdots & a_{22^n+2^{n-1}}
\end{matrix}\right]X$$
$$=\left[\begin{matrix} b_{12^{n-1}+1} & \cdots & b_{12^n} &
b_{12^n+1} & \cdots & b_{12^n+2^{n-1}} \\ b_{22^{n-1}+1} & \cdots
& b_{22^n} & b_{22^n+1} & \cdots & b_{22^n+2^{n-1}}
\end{matrix}\right]X$$ for each $X \in [\xi^n](\{e_1^1,e_2^1,x\}^n)$. Hence, by
induction hypothesis and $A_1=B_1$ and $A_2=-B_2$, we have
$$\left[\begin{matrix} b_{12^{n-1}+1} & \cdots & b_{12^n} &
-b_{12^n+1} & \cdots & -b_{12^n+2^{n-1}} \\ b_{22^{n-1}+1} &
\cdots & b_{22^n} & -b_{22^n+1} & \cdots & -b_{22^n+2^{n-1}}
\end{matrix}\right]$$
$$=\left[\begin{matrix} b_{12^{n-1}+1} & \cdots & b_{12^n} &
b_{12^n+1} & \cdots & b_{12^n+2^{n-1}} \\ b_{22^{n-1}+1} & \cdots
& b_{22^n} & b_{22^n+1} & \cdots & b_{22^n+2^{n-1}}
\end{matrix}\right].$$
Since $\begin{pmatrix} b_{1i} \\ b_{2i} \end{pmatrix} \neq
\begin{pmatrix} 0 \\ 0 \end{pmatrix}$,
$\begin{pmatrix} b_{12^n+1} & \cdots & b_{12^n+2^{n-1}} \\
b_{22^n+1} & \cdots & b_{22^n+2^{n-1}} \end{pmatrix}$ is the
$2\times2^{n-1}$ zero matrix. This is a contradiction. Therefore,
$B_2$ is the $2\times2^n$ zero matrix.

Hence, for each case, we have $A=\left[\begin{matrix} A_1 & A_2
\end{matrix}\right]=\left[\begin{matrix} B_1 & B_2
\end{matrix}\right]=B$.
This proves the lemma.
\end{proof}

\begin{thm} Let $n \in \mathbb N$, and let $k_1,\dots,k_n \in \mathbb N \cup
\{0\}$, and let $T^n,T^{k_1(1)},\dots,T^{k_n(n)}$ be
$n,k_1,\dots,k_n$-punctured ball tangle diagrams, respectively.
Then

if $F^{k_1}(T^{k_1(1)})=\left[\begin{matrix} b_{11}^1 & \cdots &
b_{12^{k_1}}^1 \\ b_{21}^1 & \cdots & b_{22^{k_1}}^1
\end{matrix}\right],\dots,
F^{k_n}(T^{k_n(n)})=\left[\begin{matrix} b_{11}^n & \cdots &
b_{12^{k_n}}^n \\ b_{21}^n & \cdots & b_{22^{k_n}}^n
\end{matrix}\right]$, then

$F^{k_1+\cdots+k_n}(T^n(T^{k_1(1)},\dots,T^{k_n(n)}))=
F^n(T^n)[\eta^n](F^{k_1}(T^{k_1(1)}),\dots,F^{k_n}(T^{k_n(n)}))$,

where $[\eta^n](F^{k_1}(T^{k_1(1)}),\dots,F^{k_n}(T^{k_n(n)}))$

$=\left[\begin{matrix} \prod_{j=1}^n
b_{\alpha_{1j}^n\alpha_{1j}^{n,2^{k_1},\dots,2^{k_n}}}^j &
\prod_{j=1}^n
b_{\alpha_{1j}^n\alpha_{2j}^{n,2^{k_1},\dots,2^{k_n}}}^j & \cdots
& \prod_{j=1}^n b_{\alpha_{1j}^n\alpha_{2^{k_1+\cdots+
k_n}j}^{n,2^{k_1},\dots,2^{k_n}}}^j
\\ \prod_{j=1}^n
b_{\alpha_{2j}^n\alpha_{1j}^{n,2^{k_1},\dots,2^{k_n}}}^j &
\prod_{j=1}^n
b_{\alpha_{2j}^n\alpha_{2j}^{n,2^{k_1},\dots,2^{k_n}}}^j & \cdots
& \prod_{j=1}^n b_{\alpha_{2j}^n\alpha_{2^{k_1+\cdots+
k_n}j}^{n,2^{k_1},\dots,2^{k_n}}}^j
\\ \cdot & \cdot & \cdot & \cdot \\ \cdot & \cdot & \cdot & \cdot
\\ \cdot & \cdot & \cdot & \cdot \\
\prod_{j=1}^n
b_{\alpha_{2^nj}^n\alpha_{1j}^{n,2^{k_1},\dots,2^{k_n}}}^j &
\prod_{j=1}^n
b_{\alpha_{2^nj}^n\alpha_{2j}^{n,2^{k_1},\dots,2^{k_n}}}^j &
\cdots & \prod_{j=1}^n b_{\alpha_{2^nj}^n\alpha_{2^{k_1+\cdots+
k_n}j}^{n,2^{k_1},\dots,2^{k_n}}}^j
\end{matrix}\right]$.
\end{thm}

\begin{proof} Without loss of generality, we may assume that
$k_1,\dots,k_n \in \mathbb N$.

Let $T=T^n(T^{k_1(1)},\dots,T^{k_n(n)})$, and let
$B^{(11)},\dots,B^{(1k_1)},\dots\dots,B^{(n1)},\dots,B^{(nk_n)}
\in \textbf{\textit{BT}}$ with
$$F^0(B^{(11)})=\left[\begin{matrix} v_1^{11} \\ v_2^{11} \end{matrix}\right],
\dots, F^0(B^{(1k_1)})=\left[\begin{matrix} v_1^{1k_1} \\
v_2^{1k_1} \end{matrix}\right],$$
$$\dots\dots,$$
$$F^0(B^{(n1)})=\left[\begin{matrix} v_1^{n1} \\ v_2^{n1} \end{matrix}\right],
\dots, F^0(B^{(nk_n)})=\left[\begin{matrix} v_1^{nk_n} \\
v_2^{nk_n} \end{matrix}\right].$$

Then $T(B^{(11)},\dots,B^{(1k_1)},\dots\dots,
B^{(n1)},\dots,B^{(nk_n)})$

$=T^n(T^{k_1(1)}(B^{(11)},\dots,B^{(1k_1)}),\dots,
T^{k_n(n)}(B^{(n1)},\dots,B^{(nk_n)}))$ and

$F^0(T(B^{(11)},\dots,B^{(1k_1)},\dots\dots,
B^{(n1)},\dots,B^{(nk_n)}))$

$=F^0(T^n(T^{k_1(1)}(B^{(11)},\dots,B^{(1k_1)}),\dots,
T^{k_n(n)}(B^{(n1)},\dots,B^{(nk_1)})))$

$=F^n(T^n)[\xi^n](F^{k_1}(T^{k_1(1)})[\xi^{k_1}](F^0(B^{(11)}),\dots,F^0(B^{(1k_1)})),
\dots, F^{k_n}(T^{k_n(n)})[\xi^{k_n}]$

$(F^0(B^{(n1)}),\dots,F^0(B^{(nk_n)})))=$

$F^n(T^n)[\xi^n](\left[\begin{matrix} b_{11}^1 & \cdots &
b_{12^{k_1}}^1 \\ b_{21}^1 & \cdots & b_{22^{k_1}}^1
\end{matrix}\right]
\left[\begin{matrix} \prod_{j=1}^{k_1} v_{\alpha_{1j}^{k_1}}^{1j}
\\ \cdot \\ \cdot \\ \cdot \\
\prod_{j=1}^{k_1} v_{\alpha_{2^{k_1}j}^{k_1}}^{1j}
\end{matrix}\right],
\dots, \left[\begin{matrix} b_{11}^n & \cdots & b_{12^{k_n}}^n
\\ b_{21}^n & \cdots & b_{22^{k_n}}^n
\end{matrix}\right]
\left[\begin{matrix} \prod_{j=1}^{k_n} v_{\alpha_{1j}^{k_n}}^{nj}
\\ \cdot \\ \cdot \\ \cdot \\
\prod_{j=1}^{k_n} v_{\alpha_{2^{k_n}j}^{k_n}}^{nj}
\end{matrix}\right])$

$=F^n(T^n)[\xi^n](\left[\begin{matrix} b_{11}^1 \prod_{j=1}^{k_1}
v_{\alpha_{1j}^{k_1}}^{1j} + \cdots +
b_{12^{k_1}}^1 \prod_{j=1}^{k_1} v_{\alpha_{2^{k_1}j}^{k_1}}^{1j} \\
b_{21}^1 \prod_{j=1}^{k_1} v_{\alpha_{1j}^{k_1}}^{1j} + \cdots +
b_{22^{k_1}}^1 \prod_{j=1}^{k_1} v_{\alpha_{2^{k_1}j}^{k_1}}^{1j}
\end{matrix}\right],\dots,$

$\left[\begin{matrix} b_{11}^n \prod_{j=1}^{k_n}
v_{\alpha_{1j}^{k_n}}^{nj} + \cdots +
b_{12^{k_n}}^n \prod_{j=1}^{k_n} v_{\alpha_{2^{k_n}j}^{k_n}}^{nj} \\
b_{21}^n \prod_{j=1}^{k_n} v_{\alpha_{1j}^{k_n}}^{nj} + \cdots +
b_{22^{k_n}}^n \prod_{j=1}^{k_n} v_{\alpha_{2^{k_n}j}^{k_n}}^{nj}
\end{matrix}\right])$

$=F^n(T^n)\left[\begin{matrix}
\xi^{n,2^{k_1},\dots,2^{k_n}}((b_{\alpha_{11}^n1}^1,\dots,
b_{\alpha_{11}^n2^{k_1}}^1),\dots,(b_{\alpha_{1n}^n1}^n,
\dots,b_{\alpha_{1n}^n2^{k_n}}^n)) \\
\xi^{n,2^{k_1},\dots,2^{k_n}}((b_{\alpha_{21}^n1}^1,\dots,
b_{\alpha_{21}^n2^{k_1}}^1),\dots,(b_{\alpha_{2n}^n1}^n,
\dots,b_{\alpha_{2n}^n2^{k_n}}^n)) \\ \cdot \\ \cdot \\ \cdot \\
\xi^{n,2^{k_1},\dots,2^{k_n}}((b_{\alpha_{2^n1}^n1}^1,\dots,
b_{\alpha_{2^n1}^n2^{k_1}}^1),\dots,(b_{\alpha_{2^nn}^n1}^n,
\dots,b_{\alpha_{2^nn}^n2^{k_n}}^n))
\end{matrix}\right]$

$\left[\begin{matrix} \prod_{j=1}^{k_1} v_{\alpha_{1j}^{k_1}}^{1j}
\cdots \prod_{j=1}^{k_{n-2}} v_{\alpha_{1j}^{k_{n-2}}}^{n-2j}
\prod_{j=1}^{k_{n-1}} v_{\alpha_{1j}^{k_{n-1}}}^{n-1j}
\prod_{j=1}^{k_n} v_{\alpha_{1j}^{k_n}}^{nj} \\ \prod_{j=1}^{k_1}
v_{\alpha_{1j}^{k_1}}^{1j} \cdots \prod_{j=1}^{k_{n-2}}
v_{\alpha_{1j}^{k_{n-2}}}^{n-2j} \prod_{j=1}^{k_{n-1}}
v_{\alpha_{1j}^{k_{n-1}}}^{n-1j} \prod_{j=1}^{k_n}
v_{\alpha_{2j}^{k_n}}^{nj} \\ \cdot \\ \cdot \\ \cdot \\
\prod_{j=1}^{k_1} v_{\alpha_{1j}^{k_1}}^{1j} \cdots
\prod_{j=1}^{k_{n-2}} v_{\alpha_{1j}^{k_{n-2}}}^{n-2j}
\prod_{j=1}^{k_{n-1}} v_{\alpha_{1j}^{k_{n-1}}}^{n-1j}
\prod_{j=1}^{k_n} v_{\alpha_{2^{k_n}j}^{k_n}}^{nj} \\
\prod_{j=1}^{k_1} v_{\alpha_{1j}^{k_1}}^{1j} \cdots
\prod_{j=1}^{k_{n-2}} v_{\alpha_{1j}^{k_{n-2}}}^{n-2j}
\prod_{j=1}^{k_{n-1}} v_{\alpha_{2j}^{k_{n-1}}}^{n-1j}
\prod_{j=1}^{k_n} v_{\alpha_{1j}^{k_n}}^{nj} \\
\prod_{j=1}^{k_1} v_{\alpha_{1j}^{k_1}}^{1j} \cdots
\prod_{j=1}^{k_{n-2}} v_{\alpha_{1j}^{k_{n-2}}}^{n-2j}
\prod_{j=1}^{k_{n-1}} v_{\alpha_{2j}^{k_{n-1}}}^{n-1j}
\prod_{j=1}^{k_n} v_{\alpha_{2j}^{k_n}}^{nj}
\\ \cdot \\ \cdot \\ \cdot \\ \cdot \\ \cdot \\ \cdot \\
\prod_{j=1}^{k_1} v_{\alpha_{2^{k_1}j}^{k_1}}^{1j} \cdots
\prod_{j=1}^{k_{n-2}} v_{\alpha_{2^{k_{n-2}}j}^{k_{n-2}}}^{n-2j}
\prod_{j=1}^{k_{n-1}} v_{\alpha_{2^{k_{n-1}}j}^{k_{n-1}}}^{n-1j}
\prod_{j=1}^{k_n} v_{\alpha_{2^{k_n}j}^{k_n}}^{nj}
\end{matrix}\right]=$

$F^{k_1+\cdots+k_n}(T)[\xi^{k_1+\cdots+k_n}](F^0(B^{(11)}),
\dots,F^0(B^{(1k_1)}),\dots\dots,
F^0(B^{(n1)}),\dots,F^0(B^{(nk_n)}))$.

Notice that there are ball tangle diagrams
$B^{(1)},B^{(2)},B^{(3)}$ such that

$F^0(B^{(1)})=\left[\begin{matrix} 1 \\ 0
\end{matrix}\right],F^0(B^{(2)})=\left[\begin{matrix} 0 \\ 1
\end{matrix}\right],F^0(B^{(3)})=\left[\begin{matrix} 1 \\ 1
\end{matrix}\right]$, respectively (See Figure 7).

Therefore, by Lemma 3.7,

$F^{k_1+\cdots+k_n}(T^n(T^{k_1(1)},\dots,T^{k_n(n)}))$

$=F^n(T^n)\left[\begin{matrix}
\xi^{n,2^{k_1},\dots,2^{k_n}}((b_{\alpha_{11}^n1}^1,\dots,
b_{\alpha_{11}^n2^{k_1}}^1),\dots,(b_{\alpha_{1n}^n1}^n,
\dots,b_{\alpha_{1n}^n2^{k_n}}^n)) \\
\xi^{n,2^{k_1},\dots,2^{k_n}}((b_{\alpha_{21}^n1}^1,\dots,
b_{\alpha_{21}^n2^{k_1}}^1),\dots,(b_{\alpha_{2n}^n1}^n,
\dots,b_{\alpha_{2n}^n2^{k_n}}^n)) \\ \cdot \\ \cdot \\ \cdot \\
\xi^{n,2^{k_1},\dots,2^{k_n}}((b_{\alpha_{2^n1}^n1}^1,\dots,
b_{\alpha_{2^n1}^n2^{k_1}}^1),\dots,(b_{\alpha_{2^nn}^n1}^n,
\dots,b_{\alpha_{2^nn}^n2^{k_n}}^n))
\end{matrix}\right]=$

$F^n(T^n)\left[\begin{matrix} \prod_{j=1}^n
b_{\alpha_{1j}^n\alpha_{1j}^{n,2^{k_1},\dots,2^{k_n}}}^j &
\prod_{j=1}^n
b_{\alpha_{1j}^n\alpha_{2j}^{n,2^{k_1},\dots,2^{k_n}}}^j & \cdots
& \prod_{j=1}^n b_{\alpha_{1j}^n\alpha_{2^{k_1+\cdots+
k_n}j}^{n,2^{k_1},\dots,2^{k_n}}}^j
\\ \prod_{j=1}^n
b_{\alpha_{2j}^n\alpha_{1j}^{n,2^{k_1},\dots,2^{k_n}}}^j &
\prod_{j=1}^n
b_{\alpha_{2j}^n\alpha_{2j}^{n,2^{k_1},\dots,2^{k_n}}}^j & \cdots
& \prod_{j=1}^n b_{\alpha_{2j}^n\alpha_{2^{k_1+\cdots+
k_n}j}^{n,2^{k_1},\dots,2^{k_n}}}^j
\\ \cdot & \cdot & \cdot & \cdot \\ \cdot & \cdot & \cdot & \cdot
\\ \cdot & \cdot & \cdot & \cdot \\
\prod_{j=1}^n
b_{\alpha_{2^nj}^n\alpha_{1j}^{n,2^{k_1},\dots,2^{k_n}}}^j &
\prod_{j=1}^n
b_{\alpha_{2^nj}^n\alpha_{2j}^{n,2^{k_1},\dots,2^{k_n}}}^j &
\cdots & \prod_{j=1}^n b_{\alpha_{2^nj}^n\alpha_{2^{k_1+\cdots+
k_n}j}^{n,2^{k_1},\dots,2^{k_n}}}^j
\end{matrix}\right]$

$=F^n(T^n)[\eta^n](F^{k_1}(T^{k_1(1)}),\dots,F^{k_n}(T^{k_n(n)}))$.

This proves the theorem.
\end{proof}

Let us give the following example.

Suppose that $T^2,T^{2(1)},T^{1(2)}$ are $2,2,1$-punctured ball
tangle diagrams such that

$F^2(T^2)=\left[\begin{matrix} a_{11} & a_{12} & a_{13} & a_{14}
\\ a_{21} & a_{22} & a_{23} & a_{24}
\end{matrix}\right]$,
$F^2(T^{2(1)})=\left[\begin{matrix} b_{11}^1 & b_{12}^1 & b_{13}^1
& b_{14}^1 \\ b_{21}^1 & b_{22}^1 & b_{23}^1 & b_{24}^1
\end{matrix}\right]$, and

$F^1(T^{1(2)})=\left[\begin{matrix} b_{11}^2 & b_{12}^2 \\
b_{21}^2 & b_{22}^2 \end{matrix}\right]$, respectively. Then
$F^3(T^2(T^{2(1)},T^{1(2)}))=$

$\left[\begin{matrix} a_{11} & a_{12} & a_{13} & a_{14}
\\ a_{21} & a_{22} & a_{23} & a_{24}
\end{matrix}\right]
\left[\begin{matrix} b_{11}^1b_{11}^2 & b_{11}^1b_{12}^2 &
b_{12}^1b_{11}^2 & b_{12}^1b_{12}^2 & b_{13}^1b_{11}^2 &
b_{13}^1b_{12}^2 & b_{14}^1b_{11}^2 & b_{14}^1b_{12}^2 \\
b_{11}^1b_{21}^2 & b_{11}^1b_{22}^2 & b_{12}^1b_{21}^2 &
b_{12}^1b_{22}^2 & b_{13}^1b_{21}^2 & b_{13}^1b_{22}^2 &
b_{14}^1b_{21}^2 & b_{14}^1b_{22}^2 \\
b_{21}^1b_{11}^2 & b_{21}^1b_{12}^2 & b_{22}^1b_{11}^2 &
b_{22}^1b_{12}^2 & b_{23}^1b_{11}^2 & b_{23}^1b_{12}^2 &
b_{24}^1b_{11}^2 & b_{24}^1b_{12}^2\\
b_{21}^1b_{21}^2 & b_{21}^1b_{22}^2 & b_{22}^1b_{21}^2 &
b_{22}^1b_{22}^2 & b_{23}^1b_{21}^2 & b_{23}^1b_{22}^2 &
b_{24}^1b_{21}^2 & b_{24}^1b_{22}^2
\end{matrix}\right]$.

Now, let us consider `(outer) connect sums' of various
$n$-punctured ball tangle diagrams and their invariants. They will
be also very useful when we compute invariants of complicate
tangles. Given $k_1$ and $k_2$-punctured ball tangle diagrams
$T^{k_1(1)}$ and $T^{k_2(2)}$, we denote the \textit{horizontal}
and the \textit{vertical} connect sums of them by $T^{k_1(1)} +_h
T^{k_2(2)}$ and $T^{k_1(1)} +_v T^{k_2(2)}$, respectively (See
Figure 10 for an example).

\begin{thm} Let $k_1, k_2 \in \mathbb N \cup
\{0\}$, and let $T^{k_1(1)}, T^{k_2(2)}$ be $k_1, k_2$-punctured
ball tangle diagrams, respectively. Then

if $F^{k_1}(T^{k_1(1)})=\left[\begin{matrix} a_{11} & a_{12} &
\cdots & a_{12^{k_1}} \\ a_{21} & a_{22} & \cdots & a_{22^{k_1}}
\end{matrix}\right]$ and $F^{k_2}(T^{k_2(2)})=\left[\begin{matrix}
b_{11} & b_{12} & \cdots & b_{12^{k_2}} \\
b_{21} & b_{22} & \cdots & b_{22^{k_2}}
\end{matrix}\right]$, then

(1) $F^{k_1+k_2}(T^{k_1(1)} +_h
T^{k_2(2)}))=\left[\begin{matrix}\left(\begin{pmatrix}
a_{1i}b_{2j}+a_{2i}b_{1j} \\ a_{2i}b_{2j}
\end{pmatrix}_{j=1,\dots,2^{k_2}}\right)_{i=1,\dots,2^{k_1}}\end{matrix}\right]$,

(2) $F^{k_1+k_2}(T^{k_1(1)} +_v
T^{k_2(2)}))=\left[\begin{matrix}\left(\begin{pmatrix}
a_{1i}b_{1j} \\ a_{2i}b_{1j}+a_{1i}b_{2j}
\end{pmatrix}_{j=1,\dots,2^{k_2}}\right)_{i=1,\dots,2^{k_1}}\end{matrix}\right]$.
\end{thm}

\begin{proof} We denote $F^0(B^{(1)} +_h B^{(2)})$ by
$F^0(B^{(1)}) +_h F^0(B^{(2)})$ and $F^0(B^{(1)} +_v B^{(2)})$ by
$F^0(B^{(1)}) +_v F^0(B^{(2)})$ if $B^{(1)},B^{(2)} \in
\textbf{\textit{BT}}$. We will define these notations $+_h$ and
$+_v$ in Chapter 4 again (See Lemma 2.12 and Definition 4.5).

(1) Let $T=T^{k_1(1)} +_h T^{k_2(2)}$, and let
$B^{(11)},\dots,B^{(1k_1)},B^{(21)},\dots,B^{(2k_2)} \in
\textbf{\textit{BT}}$ with
$$F^0(B^{(11)})=\left[\begin{matrix} v_1^{11} \\ v_2^{11} \end{matrix}\right],
\dots, F^0(B^{(1k_1)})=\left[\begin{matrix} v_1^{1k_1} \\
v_2^{1k_1} \end{matrix}\right],$$
$$F^0(B^{(21)})=\left[\begin{matrix} v_1^{21} \\ v_2^{21} \end{matrix}\right],
\dots, F^0(B^{(2k_2)})=\left[\begin{matrix} v_1^{2k_2} \\
v_2^{2k_2} \end{matrix}\right].$$

Then $T(B^{(11)},\dots,B^{(1k_1)},B^{(21)},\dots,B^{(2k_2)})$

$=T^{k_1(1)}(B^{(11)},\dots,B^{(1k_1)}) +_h
T^{k_2(2)}(B^{(21)},\dots,B^{(2k_2)})$ and

$F^0(T(B^{(11)},\dots,B^{(1k_1)},B^{(21)},\dots,B^{(2k_2)}))$

$=F^0(T^{k_1(1)}(B^{(11)},\dots,B^{(1k_1)}) +_h
T^{k_2(2)}(B^{(21)},\dots,B^{(2k_2)}))$

$=F^0(T^{k_1(1)}(B^{(11)},\dots,B^{(1k_1)})) +_h
F^0(T^{k_2(2)}(B^{(21)},\dots,B^{(2k_2)}))$

$=F^{k_1}(T^{k_1(1)})[\xi^{k_1}](F^0(B^{(11)}),\dots,F^0(B^{(1k_1)}))$

$+_h
F^{k_2}(T^{k_2(2)})[\xi^{k_2}](F^0(B^{(21)}),\dots,F^0(B^{(2k_2)}))$

$=\left[\begin{matrix} a_{11} & \cdots & a_{12^{k_1}} \\
a_{21} & \cdots & a_{22^{k_1}}
\end{matrix}\right]
\left[\begin{matrix} \prod_{j=1}^{k_1} v_{\alpha_{1j}^{k_1}}^{1j}
\\ \cdot \\ \cdot \\ \cdot \\
\prod_{j=1}^{k_1} v_{\alpha_{2^{k_1}j}^{k_1}}^{1j}
\end{matrix}\right] +_h
\left[\begin{matrix} b_{11} & \cdots & b_{12^{k_2}} \\
b_{21} & \cdots & b_{22^{k_2}}
\end{matrix}\right]
\left[\begin{matrix} \prod_{j=1}^{k_2} v_{\alpha_{1j}^{k_2}}^{2j}
\\ \cdot \\ \cdot \\ \cdot \\
\prod_{j=1}^{k_2} v_{\alpha_{2^{k_2}j}^{k_2}}^{2j}
\end{matrix}\right]=$

$$\left[\begin{matrix} a_{11} \prod_{j=1}^{k_1}
v_{\alpha_{1j}^{k_1}}^{1j} + \cdots +
a_{12^{k_1}} \prod_{j=1}^{k_1} v_{\alpha_{2^{k_1}j}^{k_1}}^{1j} \\
a_{21} \prod_{j=1}^{k_1} v_{\alpha_{1j}^{k_1}}^{1j} + \cdots +
a_{22^{k_1}} \prod_{j=1}^{k_1} v_{\alpha_{2^{k_1}j}^{k_1}}^{1j}
\end{matrix}\right] +_h
\left[\begin{matrix} b_{11} \prod_{j=1}^{k_2}
v_{\alpha_{1j}^{k_2}}^{2j} + \cdots +
b_{12^{k_2}} \prod_{j=1}^{k_2} v_{\alpha_{2^{k_2}j}^{k_2}}^{2j} \\
b_{21} \prod_{j=1}^{k_2} v_{\alpha_{1j}^{k_2}}^{2j} + \cdots +
b_{22^{k_2}} \prod_{j=1}^{k_2} v_{\alpha_{2^{k_2}j}^{k_2}}^{2j}
\end{matrix}\right]$$

$=\left[\begin{matrix} a_{11}b_{21}+a_{21}b_{11} & a_{21}b_{21} \\
\cdot & \cdot \\ \cdot & \cdot \\ \cdot & \cdot \\
a_{11}b_{22^{k_2}}+a_{21}b_{12^{k_2}} & a_{21}b_{22^{k_2}} \\
\cdot & \cdot \\ \cdot & \cdot \\ \cdot & \cdot \\
\cdot & \cdot \\ \cdot & \cdot \\ \cdot & \cdot \\
a_{12^{k_1}}b_{21}+a_{22^{k_1}}b_{11} & a_{22^{k_1}}b_{21} \\
\cdot & \cdot \\ \cdot & \cdot \\ \cdot & \cdot \\
a_{12^{k_1}}b_{22^{k_2}}+a_{22^{k_1}}b_{12^{k_2}} &
a_{22^{k_1}}b_{22^{k_2}} \end{matrix}\right]^\dag
\left[\begin{matrix} \prod_{j=1}^{k_1} v_{\alpha_{1j}^{k_1}}^{1j}
\prod_{j=1}^{k_2} v_{\alpha_{1j}^{k_2}}^{2j} \\
\cdot \\ \cdot \\ \cdot \\
\prod_{j=1}^{k_1} v_{\alpha_{1j}^{k_1}}^{1j}
\prod_{j=1}^{k_2} v_{\alpha_{2^{k_2}j}^{k_2}}^{2j} \\
\cdot \\ \cdot \\ \cdot \\ \cdot \\ \cdot \\ \cdot \\
\prod_{j=1}^{k_1} v_{\alpha_{1j}^{k_1}}^{1j}
\prod_{j=1}^{k_2} v_{\alpha_{1j}^{k_2}}^{2j} \\
\cdot \\ \cdot \\ \cdot \\
\prod_{j=1}^{k_1} v_{\alpha_{2^{k_1}j}^{k_1}}^{1j}
\prod_{j=1}^{k_2} v_{\alpha_{2^{k_2}j}^{k_2}}^{2j}
\end{matrix}\right]$

$=F^{k_1+k_2}(T)[\xi^{k_1+k_2}](F^0(B^{(11)}),
\dots,F^0(B^{(1k_1)}),F^0(B^{(21)}),\dots,F^0(B^{(2k_2)}))$.

Notice that there are ball tangle diagrams
$B^{(1)},B^{(2)},B^{(3)}$ such that

$F^0(B^{(1)})=\left[\begin{matrix} 1 \\ 0
\end{matrix}\right],F^0(B^{(2)})=\left[\begin{matrix} 0 \\ 1
\end{matrix}\right],F^0(B^{(3)})=\left[\begin{matrix} 1 \\ 1
\end{matrix}\right]$, respectively (See Figure 7).

Therefore, by Lemma 3.7,

$F^{k_1+k_2}(T^{k_1(1)} +_h
T^{k_2(2)})=\left[\begin{matrix}\left(\begin{pmatrix}
a_{1i}b_{2j}+a_{2i}b_{1j} \\ a_{2i}b_{2j}
\end{pmatrix}_{j=1,\dots,2^{k_2}}\right)_{i=1,\dots,2^{k_1}}\end{matrix}\right]$.

(2) Similarly, we can show that

$F^{k_1+k_2}(T^{k_1(1)} +_v
T^{k_2(2)}))=\left[\begin{matrix}\left(\begin{pmatrix}
a_{1i}b_{1j} \\ a_{2i}b_{1j}+a_{1i}b_{2j}
\end{pmatrix}_{j=1,\dots,2^{k_2}}\right)_{i=1,\dots,2^{k_1}}\end{matrix}\right]$.

This proves the theorem.
\end{proof}

Notice that each of $+_h$ and $+_v$ gives us a monoid action if
one of $k_1$ and $k_2$ is fixed and the other is 0. Obviously,
they are binary operations which give monoid structures to
$\textbf{\textit{BT}}$ if $k_1=k_2=0$.

We will consider these connect sums again for spherical tangles
and ball tangles including `inner connect sums' to yield monoid
structures on $\textbf{\textit{ST}}$ in Chapter 4.

\vskip 0.2in

\bigskip
\centerline{\epsfxsize=5.5 in \epsfbox{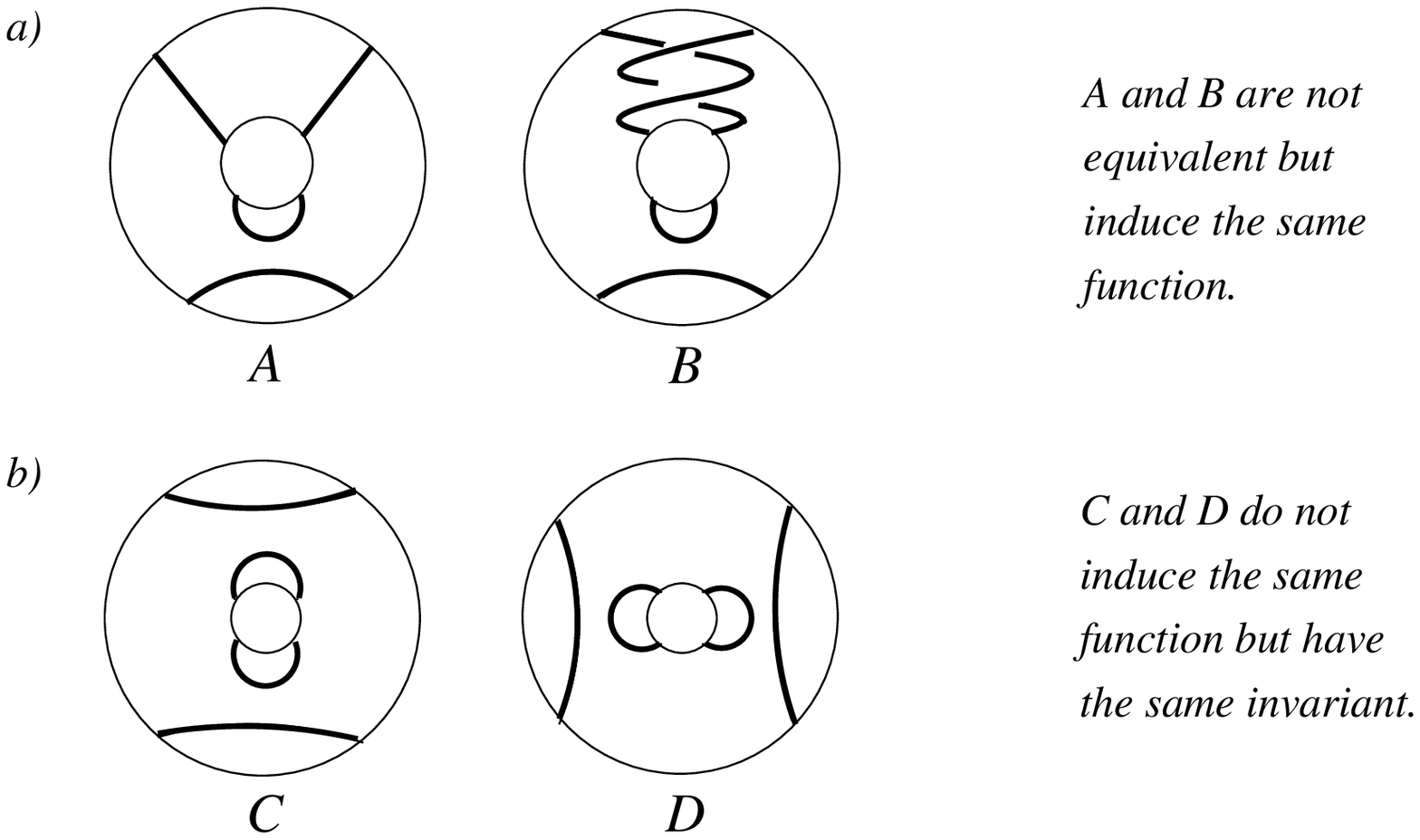}}
\medskip
\centerline{\small Figure 6. Tangles and functions.}
\bigskip

When we denote the statement that $n$-punctured ball tangles
$T^{n(1)}$ and $T^{n(2)}$ induces the same function from
$\textbf{\textit{BT}}^n$ to $\textbf{\textit{BT}}$ by $T^{n(1)}
\simeq T^{n(2)}$ and $F^n(T^{n(1)})=F^n(T^{n(2)})$ by $T^{n(1)}
\sim T^{n(2)}$, $\simeq$ and $\sim$ are clearly equivalence
relations on $\textbf{\textit{nPBT}}$ and we have
$$T^{n(1)} \cong T^{n(2)} \Longrightarrow T^{n(1)} \simeq T^{n(2)}
\Longrightarrow T^{n(1)} \sim T^{n(2)}.$$ The first implication
comes from the definition of $\cong$ and the second implication is
proved by Theorem 3.6 and Lemma 3.7 immediately.

Notice that neither the converse of the first implication nor that
of the second implication is true (See Figure 6). In particular,
the spherical tangles $C$ and $D$ in Figure 6 have the zero matrix
invariant. For a nonzero matrix invariant, we can take the
spherical tangle $A$ in Figure 6 and a spherical tangle $B'$
obtained from a single twist of two upper strands of $A$. When we
apply a ball tangle whose numerical closure is a noninvertible
knot for the holes, $A$ and $B'$ give us different ball tangles.
However, $A$ and $B'$ have the same invariant. By theses reasons,
we may consider the equivalence relation $\simeq$ instead of
$\cong$ for our $n$-punctured ball tangle invariant.

This aspect is quite similar to that in Algebraic Topology in the
sense as follows:

If $X$ and $Y$ are pathconnected topological spaces, then
$$X \cong Y \Longrightarrow X \simeq Y \Longrightarrow X \sim Y,$$
where $X \cong Y$, $X \simeq Y$, and $X \sim Y$ mean the
statements that $X$ and $Y$ are topologically equivalent, $X$ and
$Y$ are homotopically equivalent, and $\pi_1(X)$ and $\pi_1(Y)$
are isomorphic, respectively.

\section{A Generalization of Krebes' Theorem}

Suppose that a ball tangle diagram $B$ is embedded in a link
diagram $L$. Since the complement of a ball in $S^3$ is still a
ball, we may think of $B'=\overline{L-B}$ as another ball tangle
diagram embedded in $L$ and $L=(B +_h B')_1$, that is, $L$ is the
numerator closure of horizontal addition of $B$ and $B'$. If
$F^0(B)=\left[\begin{matrix} p_1 \\ q_1 \end{matrix}\right]$ and
$F^0(B')=\left[\begin{matrix} p \\ q \end{matrix}\right]$, then
$F^0(B +_h B')=\left[\begin{matrix} p_1q + q_1p \\ q_1q
\end{matrix}\right]$ and $|\langle L \rangle|=|\langle (B +_h B')_1
\rangle|=|p_1q + q_1p|$. Hence, ${\rm g.c.d.}\,(p_1,q_1)$ divides
$|\langle L \rangle|$. This is Krebes' Theorem (See Theorem 2.14).
We have the following generalization of Krebes' theorem.

\begin{thm} Let $L$ be a link, and let $B^{(1)},\dots,B^{(n)}$
be ball tangles with the invariants $\left[\begin{matrix} p_1 \\
q_1 \end{matrix}\right], \dots, \left[\begin{matrix} p_n \\ q_n
\end{matrix}\right]$, respectively. If $B^{(1)},\dots,B^{(n)}$
are embedded in $L$ disjointly, then $\prod_{i=1}^k\,{\rm
g.c.d.}\,(p_i,q_i)$ divides $|\langle L \rangle|$.
\end{thm}

\begin{proof} Denote by $d_i={\rm g.c.d.}\,(p_i,q_i)$ for each $i \in
\{1,\dots,n\}$. Let $B^{(n)'}=\overline{L-B^{(n)}}$, and let
$F^0(B^{(n)'})=\left[\begin{matrix} p \\ q
\end{matrix}\right]$. Then $L=(B^{(n)} +_h B^{(n)'})_1$ and
$F^0(B^{(n)} +_h B^{(n)'})=\left[\begin{matrix} p_nq + q_np \\
q_nq \end{matrix}\right]$, hence, $|\langle L \rangle|=|\langle
(B^{(n)} +_h B^{(n)'})_1 \rangle|=|p_nq + q_np|$. Notice that we
can regard $B^{(n)'}$ as an $(n-1)$-punctured ball tangle with its
holes filled up by $B^{(1)},\dots,B^{(n-1)}$. Hence,
$B^{(n)'}=T^{n-1}(B^{(1)},\dots,B^{(n-1)})$ for some
$(n-1)$-punctured ball tangle $T^{n-1}$.

Let $F^{n-1}(T^{n-1})=\left[\begin{matrix} a_{11} & a_{12} &
\cdots & a_{12^{n-1}} \\
a_{21} & a_{22} & \cdots & a_{22^{n-1}}
\end{matrix}\right]$. Then, by Lemma 3.4  and Theorem 3.6, we have
$$
\begin{aligned}
&F^0(B^{(n)'})=F^0(T^{n-1}(B^{(1)},\dots,B^{(n-1)}))\\
&=F^{n-1}(T^{n-1})[\xi^{n-1}](F^0(B^{(1)}),
\dots,F^0(B^{(n-1)}))\\
&=\left[\begin{matrix} \begin{pmatrix} a_{11} & a_{12} &
\cdots & a_{12^{n-1}} \\
a_{21} & a_{22} & \cdots & a_{22^{n-1}}
\end{pmatrix}
\begin{pmatrix} p_1p_2\cdots p_{n-1} \\
p_1p_2\cdots q_{n-1}
\\ \cdot \\ \cdot \\ \cdot \\
q_1q_2\cdots q_{n-1}
\end{pmatrix} \end{matrix}\right]
\end{aligned}$$
$$=\left[\begin{matrix} a_{11}p_1p_2\cdots p_{n-1}+
\cdots\cdots + a_{12^{n-1}}q_1q_2\cdots q_{n-1} \\
a_{21}p_1p_2\cdots p_{n-1}+ \cdots\cdots +
a_{22^{n-1}}q_1q_2\cdots q_{n-1} \end{matrix}\right].
$$

Let $d'={\rm g.c.d.}\,(p,q)$. Then
$$
\begin{aligned}
d'={\rm g.c.d.}\,(&a_{11}p_1p_2\cdots p_{n-1}+ \cdots\cdots +
a_{12^{n-1}}q_1q_2\cdots
q_{n-1},\\
&a_{21}p_1p_2\cdots p_{n-1}+ \cdots\cdots +
a_{22^{n-1}}q_1q_2\cdots q_{n-1})
\end{aligned}$$
and there are $k,l \in \mathbb{Z}$ such that
$d'=k(a_{11}p_1p_2\cdots p_{n-1}+ \cdots\cdots +
a_{12^{n-1}}q_1q_2\cdots q_{n-1})+l(a_{21}p_1p_2\cdots p_{n-1}+
\cdots\cdots + a_{22^{n-1}}q_1q_2\cdots q_{n-1})$. Since
$d_1\cdots d_{n-1}$ divides each term of the above linear
combination, $d_1\cdots d_{n-1}$ divides $d'$. Hence, $d_1\cdots
d_{n-1}d_n$ divides $d'd_n$ and $d'd_n$ divides $|\langle L
\rangle|$. Therefore, $\prod_{i=1}^k\,{\rm g.c.d.}\,(p_i,q_i)$
divides $|\langle L \rangle|$.
\end{proof}

\chapter{Surjectivity of invariants}

We use the following notation throughout this section:

(1) The subscripts 1,2 of ball tangles will no longer used to
denote different kinds of closures. They will be used simply to
distinguish different ball tangles.

(2) The ball tangle invariant $F^0$ will be denoted by $f$ with
values in $PM_2=PM_{2\times1}(\mathbb{Z})$ and the spherical
tangle invariant $F^1$ will be denoted by $F$ with values in
$PM_{2\times2}=PM_{2\times2}(\mathbb{Z})$.

\begin{df} Let $a_1, a_2, a_3, a_4$ be 4 points in $S^2$, and let
$x_i=\{a_i\}\times\{0\}$ and $y_i=\{a_i\}\times\{1\}$ for each $i
\in \{1,2,3,4\}$. Then a 1-dimensional proper submanifold $S$ of
$S^2 \times I$, $I=[0,1]$, is called a spherical tangle about
$a_1, a_2, a_3, a_4$ (or simply, spherical tangle) if $\partial S
\cap (S^2 \times \{0\})=\{x_1, x_2, x_3, x_4\}$ and $\partial S
\cap (S^2 \times \{1\})=\{y_1, y_2, y_3, y_4\}$.
\end{df}

Note that $B((0,0,0),2)-{\rm Int}(B((0,0,0),1))$ is homeomorphic
to $S^2 \times I$.

\newpage

\begin{df} Define $\cong$ on the class of all spherical tangles
about $a_1, a_2, a_3, a_4$ by $S_1 \cong S_2$ if and only if there
is a homeomorphism $h:S^2 \times I \rightarrow S^2 \times I$ such
that $h(S_1)=S_2$ and $h \simeq Id_{S_2 \times I}$ {\rm rel}
$\partial$ for spherical tangles $S_1$ and $S_2$. Then $\cong$ is
an equivalence relation on it. $S_1$ and $S_2$ are said to be
isotopic, or of the same isotopy type, if $S_1 \cong S_2$ and, for
each spherical tangle $S$, the equivalence class $[S]$ is called
the isotopy type of $S$.
\end{df}

Remark that we usually use $S$ for $[S]$ and consider only
diagrams for spherical tangles and ball tangles. Now, let us
define the product of spherical tangle diagrams as follows:
$$[S_2] \circ [S_1]=[S_2(S_1)],$$ or simply, $S_2 \circ S_1=S_2(S_1)$
for all spherical tangle diagrams $S_1$ and $S_2$, where, roughly
speaking, $S_2(S_1)$ means to put $S_1$ inside of $S_2$, using the
identification
$$\frac{(S^2\times[0,1])_1 \coprod (S^2\times[0,1])_2}
{(S^2\times \{1\})_1=(S^2\times \{0\})_2}=S^2\times [0,1].$$ It is
clear that $\circ$ is associative and
$\textbf{\textit{I}}=\coprod_{i=1}^4\,\{a_i\}\times I$ is the
identity spherical tangle for $\circ$. Thus, the class
$\textbf{\textit{ST}}$ of spherical tangle diagrams with $\circ$
forms a monoid.

\section{Surjectivity of the ball tangle invariant $f$}

Recall Lemma 2.12, and Lemma 2.13:

(1) If $B_1, B_2 \in \textbf{\textit{BT}}$ and
$f(B_1)=\left[\begin{matrix} p \\ q \end{matrix}\right],
f(B_2)=\left[\begin{matrix} r \\ s \end{matrix}\right]$, then
$f(B_1 +_h B_2)=\left[\begin{matrix} ps+qr \\ qs
\end{matrix}\right]$.
So if we denote $\left[\begin{matrix} ps+qr \\ qs
\end{matrix}\right]$ by $\left[\begin{matrix} p \\ q \end{matrix}\right]
+_h \left[\begin{matrix} r \\ s \end{matrix}\right]$, then we have
$f(B_1 +_h B_2)=f(B_1) +_h f(B_2)$.

(2) If $B \in \textbf{\textit{BT}}$ and $f(B)=\left[\begin{matrix}
p \\ q \end{matrix}\right]$, then $f(B^*)=\left[\begin{matrix} p
\\ -q \end{matrix}\right]$ and $f(B^R)=\left[\begin{matrix} q \\
-p \end{matrix}\right]$, where $B^*$ is the mirror image of $B$
and $B^R$ is the $90^{\circ}$ rotation of $B$ counterclockwise on
the projection plane. So if we denote $\left[\begin{matrix} p \\
-q \end{matrix}\right]$ by $\left[\begin{matrix} p \\ q
\end{matrix}\right]^*$ and $\left[\begin{matrix} q \\ -p
\end{matrix}\right]$ by $\left[\begin{matrix} p \\ q
\end{matrix}\right]^R$, then we have $f(B^*)=f(B)^*$ and
$f(B^R)=f(B)^R$.

To avoid complication, we use the same notations for $+_h$, $^*$,
and $^R$ applied to $\textbf{\textit{BT}}$ and $PM_2$. We shall be
able to understand the meaning of different operations by their
contexts.

(3) If $B \in \textbf{\textit{BT}}$, then $B^{**}=B$ but $B^{RR}$
need not be the same as $B$.

(4) If $A \in PM_2$, then $A^{**}=A$ and $A^{RR}=A$.

Notice that, if an element $A$ in $PM_2$ can be obtained by
applying $+_h$, $^*$, and $^R$ to finitely many invariants of ball
tangles, then $A$ itself is the invariant of a ball tangle.

Let us calculate $f$ for ball tangles in Figure 7.

1. The ball tangles \textbf{\textit{b}} and \textbf{\textit{c}}
have invariants $\left[\begin{matrix} 1
\\ 0 \end{matrix}\right]$ and $\left[\begin{matrix} 0 \\ 1
\end{matrix}\right]$, respectively.

2. The ball tangle \textbf{\textit{a}} has invariant
$\left[\begin{matrix} 0 \\ 0 \end{matrix}\right]$ because
$\textbf{\textit{a}}=\textbf{\textit{b}} +_h \textbf{\textit{b}}$
and $\left[\begin{matrix} 0 \\ 0 \end{matrix}\right]=
\left[\begin{matrix} 1 \\ 0 \end{matrix}\right] +_h
\left[\begin{matrix} 1 \\ 0 \end{matrix}\right]$.

\bigskip
\centerline{\epsfxsize=4.8 in \epsfbox{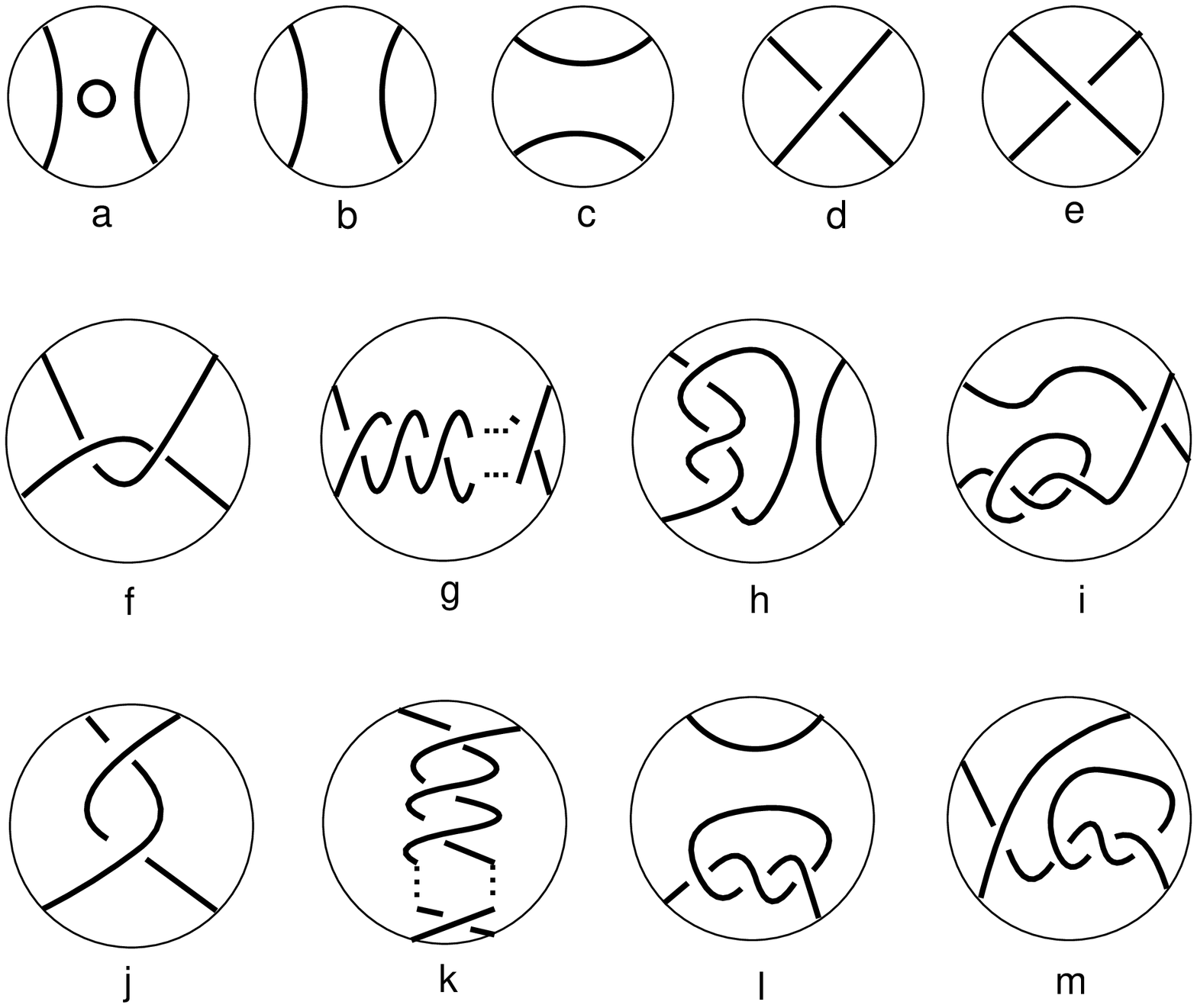}}
\medskip
\centerline{\small Figure 7. Ball tangle diagrams.}
\bigskip

3. The ball tangles \textbf{\textit{d}} and \textbf{\textit{e}}
have invariants $\left[\begin{matrix} 1 \\ 1
\end{matrix}\right]$ and $\left[\begin{matrix} 1
\\ -1 \end{matrix}\right]$, respectively. They are the mirror images each other.

4. The ball tangle \textbf{\textit{f}} has invariant
$\left[\begin{matrix} 1 \\ 1 \end{matrix}\right] +_h \left[\begin{matrix} 1 \\
1 \end{matrix}\right]=\left[\begin{matrix} 2 \\
1 \end{matrix}\right]$, and the ball tangle \textbf{\textit{g}}
has invariant $\left[\begin{matrix} p \\ 1 \end{matrix}\right]$,
where $p$ is the number of horizontal twists in
\textbf{\textit{g}}.

5. Ball tangles \textbf{\textit{j}} and \textbf{\textit{k}} have
invariants $\left[\begin{matrix} 1
\\ 2 \end{matrix}\right]$ and $\left[\begin{matrix} 1 \\ q
\end{matrix}\right]$, respectively, where $q$ is the number of
vertical twists in \textbf{\textit{k}}.

6. The ball tangle \textbf{\textit{h}} has invariant
$\left[\begin{matrix} 3 \\
0 \end{matrix}\right]$ because $\left[\begin{matrix} 3 \\
0 \end{matrix}\right]=\left[\begin{matrix} 1 \\
3 \end{matrix}\right] +_h \left[\begin{matrix} 1 \\
0 \end{matrix}\right]$. The ball tangle \textbf{\textit{l}} is
$\textbf{\textit{h}}^R$ and has invariant
$\left[\begin{matrix} 0 \\
-3 \end{matrix}\right]=\left[\begin{matrix} 0 \\
3 \end{matrix}\right]$.

7. The ball tangle \textbf{\textit{m}} has invariant
$\left[\begin{matrix} 3 \\ 3 \end{matrix}\right]$ because
$\textbf{\textit{m}}=\textbf{\textit{d}} +_h \textbf{\textit{l}}$
and $\left[\begin{matrix} 3 \\ 3 \end{matrix}\right]=
\left[\begin{matrix} 1 \\ 1 \end{matrix}\right] +_h
\left[\begin{matrix} 0 \\ 3 \end{matrix}\right]$.
$\textbf{\textit{i}}=\textbf{\textit{l}} +_h \textbf{\textit{d}}$.
Hence, the ball tangles \textbf{\textit{i}} and
\textbf{\textit{m}} have the same invariant but they are
apparently not isotopic.

To prove the surjectivity of $f$, we use Euclidean Algorithm.

\begin{pro}{\rm (Euclidean Algorithm)}
If $a,b \in \mathbb{N}$ and $a<b$, then there are uniquely $k \in
\mathbb{N}$ and $r_0,r_1,\dots,r_k,r_{k+1} \in \mathbb{N} \cup
\{0\}$ and $q_1,\dots,q_{k+1} \in \mathbb{N}$ such that $r_0=a$
and $r_{k+1}=0$ and $r_0>r_1>\cdots>r_k>r_{k+1}$ and $b=q_1a+r_1$
and $r_{i-2}=q_ir_{i-1}+r_i$ for each $i \in \{2,\dots,k+1\}$.
\end{pro}

\begin{thm} The ball tangle invariant
$f:\textbf{\textit{BT}} \rightarrow PM_2$ is onto.
\end{thm}

\begin{proof} It is enough to show that there is
$B \in \textbf{\textit{BT}}$ such that $f(B)=\left[\begin{matrix}
b \\ a \end{matrix}\right]$ if $a,b \in\mathbb{N}$ and $a<b$.

Suppose that $a,b \in\mathbb{N}$ and $a<b$. Then, by Euclidean
Algorithm, there are uniquely $k \in \mathbb{N}$ and
$r_0,r_1,\dots,r_k,r_{k+1} \in\mathbb{N} \cup \{0\}$ and
$q_1,\dots,q_{k+1} \in\mathbb{N}$ such that $r_0=a$ and

$$b=q_1a+r_1, \left[\begin{matrix} b \\ a \end{matrix}\right] =
\left[\begin{matrix} q_1 \\ 1
\end{matrix}\right] +_h
\left[\begin{matrix} r_1 \\ a \end{matrix}\right], 0<r_1<a,$$
$$a=q_2r_1+r_2, \left[\begin{matrix} a \\ r_1 \end{matrix}\right] =
\left[\begin{matrix} q_2 \\ 1
\end{matrix}\right] +_h
\left[\begin{matrix} r_2 \\ r_1 \end{matrix}\right], 0<r_2<r_1,$$
$$r_1=q_3r_2+r_3, \left[\begin{matrix} r_1 \\ r_2 \end{matrix}\right] =
\left[\begin{matrix} q_2 \\ 1
\end{matrix}\right] +_h
\left[\begin{matrix} r_3 \\ r_2 \end{matrix}\right], 0<r_3<r_2,$$
$$\cdots\cdots\cdots$$
$$r_{k-2}=q_kr_{k-1}+r_k, \left[\begin{matrix} r_{k-2} \\ r_{k-1}
\end{matrix}\right] =
\left[\begin{matrix} q_k \\ 1
\end{matrix}\right] +_h
\left[\begin{matrix} r_k \\ r_{k-1} \end{matrix}\right],
0<r_{k-1}<r_{k-2},$$
$$r_{k-1}=q_{k+1}r_k+r_{k+1}, \left[\begin{matrix} r_{k-1} \\ r_k
\end{matrix}\right] =
\left[\begin{matrix} q_{k+1} \\ 1
\end{matrix}\right] +_h
\left[\begin{matrix} 0 \\ r_k \end{matrix}\right], r_{k+1}=0.$$

Since $\left[\begin{matrix} q_1 \\ 1
\end{matrix}\right],\dots,\left[\begin{matrix} q_{k+1} \\
1 \end{matrix}\right]$, and $\left[\begin{matrix} 0 \\
r_k \end{matrix}\right]$ are realizable by ball tangles and
$\left[\begin{matrix} r_i \\ r_{i-1} \end{matrix}\right]
=\left[\begin{matrix} r_{i-1} \\ r_i \end{matrix}\right]^{R*}$ for
each $i \in \{1,\dots,k\}$, $\left[\begin{matrix} b \\ a
\end{matrix}\right]$ corresponds a ball tangle.

Therefore, there is $B \in \textbf{\textit{BT}}$ such that
$\left[\begin{matrix} b \\ a \end{matrix}\right]=f(B)$. This
proves the theorem.
\end{proof}

We can define vertical connect sum of ball tangle diagrams by
horizontal connect sum and rotations.

\begin{df} Define $+_v$ on $\textbf{\textit{BT}}$ by
$B_1 +_v B_2=(B_1^R +_h B_2^R)^{RRR}$ for all $B_1, B_2 \in
\textbf{\textit{BT}}$. Then $+_v$ is called the vertical connect
sum on $\textbf{\textit{BT}}$.
\end{df}

(2*) If $B_1, B_2 \in \textbf{\textit{BT}}$ and
$f(B_1)=\left[\begin{matrix} p \\ q \end{matrix}\right],
f(B_2)=\left[\begin{matrix} r \\ s \end{matrix}\right]$, then
$f(B_1 +_v B_2)=\left[\begin{matrix} pr \\ qr+ps
\end{matrix}\right]$.
So if we denote  $\left[\begin{matrix} pr \\ qr+ps
\end{matrix}\right]$ by $\left[\begin{matrix} p \\ q \end{matrix}\right]
+_v \left[\begin{matrix} r \\ s \end{matrix}\right]$, we have
$f(B_1 +_v B_2)=f(B_1) +_v f(B_2)$.

Note that $(\textbf{\textit{BT}}, +_h)$ and
$(\textbf{\textit{BT}}, +_v)$ are noncommutative monoids with
identities \textbf{\textit{c}} and \textbf{\textit{b}} in Figure
7, respectively. On the other hand, $(PM_2, +_h)$ and $(PM_2,
+_v)$
are commutative monoids with identities $\left[\begin{matrix} 0 \\
1 \end{matrix}\right]$ and $\left[\begin{matrix} 1 \\ 0
\end{matrix}\right]$, respectively. The ball tangle invariant $f$
is a monoid epimorphism from $(\textbf{\textit{BT}}, +_h)$ and
$(\textbf{\textit{BT}}, +_v)$ to $(PM_2, +_h)$ and $(PM_2, +_v)$,
respectively.

\section{The invariant $F$ of spherical tangles, connect sums,
and determinants}

The following lemma tells us a unique commutative square.

\begin{lm} For each $S \in \textbf{\textit{ST}}$, there is a
unique function $S_*:PM_2 \rightarrow PM_2$ such that $f \circ S =
S_* \circ f$. Furthermore, $S_*$ is the function from $PM_2$ to
$PM_2$ defined by $S_*(A)=F(S)A$ for each $A \in PM_2$.
\end{lm}

\begin{proof} Let $S \in \textbf{\textit{ST}}$. Then $f(S(B))=F(S)f(B)$
for each $B \in \textbf{\textit{BT}}$ (Theorem 3.6). Hence, $f
\circ S = S_* \circ f$. The uniqueness of $S_*$ follows from the
surjectivity of $f$. To show the uniqueness of $S_*$, suppose that
$F_1$ and $F_2$ are functions from $PM_2$ to $PM_2$ such that $f
\circ S = F_1 \circ f$ and $f \circ S = F_2 \circ f$,
respectively, and $A \in PM_2$. Then there is $B \in
\textbf{\textit{BT}}$ such that $A=f(B)$ by the surjectivity of
$f$ (Theorem 4.4) and $F_1(A)=F_1(f(B))=(F_1 \circ f)(B)=(F_2
\circ f)(B)=F_2(f(B))=F_2(A)$. Hence, $F_1 = F_2$, that is, $S_*$
is unique. This proves the lemma.
\end{proof}

\begin{lm} If $S_1, S_2 \in \textbf{\textit{ST}}$, then
$(S_2 \circ S_1)_* = S_{2*} \circ S_{1*}$.
\end{lm}

\begin{proof} Suppose that $S_1, S_2 \in \textbf{\textit{ST}}$. Then
$f \circ S_1 = S_{1*} \circ f$ and $f \circ S_2 = S_{2*} \circ f$.
Hence, $f \circ (S_2 \circ S_1) = (f \circ S_2) \circ S_1 =
(S_{2*} \circ f) \circ S_1 = S_{2*} \circ (f \circ S_1) = S_{2*}
\circ (S_{1*} \circ f) = (S_{2*} \circ S_{1*}) \circ f$.
Therefore, by the uniqueness of $(S_2 \circ S_1)_*$, $(S_2 \circ
S_1)_* = S_{2*} \circ S_{1*}$.
\end{proof}

Let us identify $S_*$ with $F(S)$ for each $S \in
\textbf{\textit{ST}}$. Since $S_2 \circ S_1$ is the composition of
$S_1$ and $S_2$, we have the following lemma immediately.

\begin{lm} If $S_1, S_2 \in \textbf{\textit{ST}}$, then
$F(S_2 \circ S_1)=F(S_2)F(S_1)$.
\end{lm}

Notice that Lemma 4.8 does not depend on the surjectivity of $f$.
We can prove Lemma 4.8 by Theorem 3.6 and the following lemma, a
corollary of Lemma 3.7, without using the surjectivity of $f$.

\begin{lm} If $A, B \in PM_{2\times2}$ and

$A\left[\begin{matrix} 1 \\ 0 \end{matrix}\right]=
B\left[\begin{matrix} 1 \\ 0 \end{matrix}\right],
A\left[\begin{matrix} 0 \\ 1 \end{matrix}\right]=
B\left[\begin{matrix} 0 \\ 1 \end{matrix}\right],
A\left[\begin{matrix} 1 \\ 1 \end{matrix}\right]=
B\left[\begin{matrix} 1 \\ 1 \end{matrix}\right]$, then $A=B$.
\end{lm}

\bigskip
\centerline{\epsfxsize=5 in \epsfbox{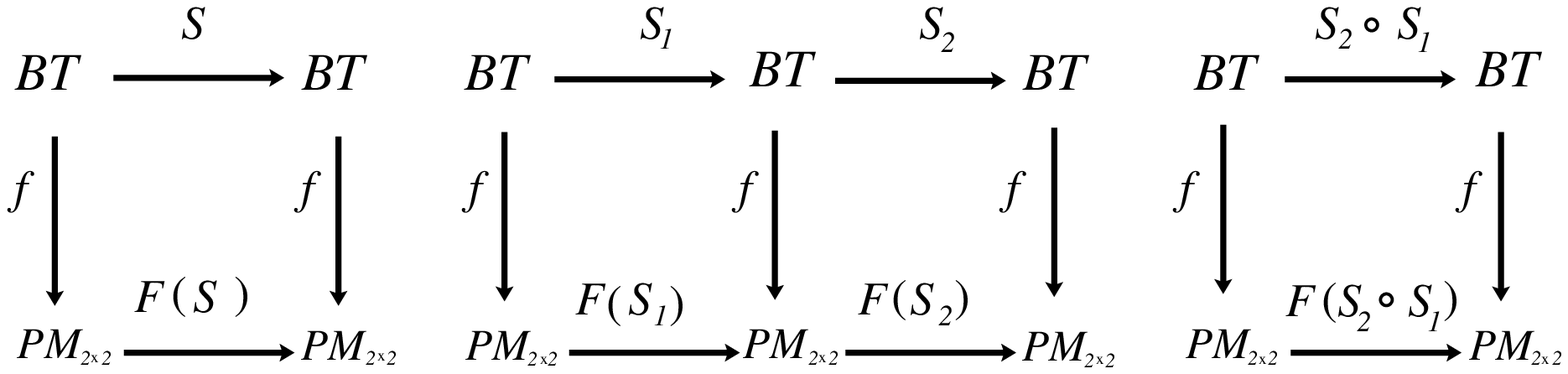}}
\medskip
\centerline{\small Figure 8. Commutative diagrams of invariants.}
\bigskip

By Theorem 3.6, we have that $F(S_2 \circ S_1)$ and $F(S_2)F(S_1)$
are matrices in $PM_{2\times2}$ such that $F(S_2 \circ
S_1)\left[\begin{matrix}
1 \\ 0 \end{matrix}\right]=F(S_2)F(S_1)\left[\begin{matrix} 1 \\
0 \end{matrix}\right]$, $F(S_2 \circ S_1)\left[\begin{matrix} 0 \\
1 \end{matrix}\right]=F(S_2)F(S_1)\left[\begin{matrix} 0 \\ 1
\end{matrix}\right]$, and $ F(S_2 \circ S_1)\left[\begin{matrix} 1 \\ 1
\end{matrix}\right]=F(S_2)F(S_1)\left[\begin{matrix} 1 \\ 1
\end{matrix}\right]$. Hence,
$F(S_2 \circ S_1)=F(S_2)F(S_1)$ by Lemma 4.9.

Let us introduce the elementary operations on
$\textbf{\textit{ST}}$.

\begin{df} Let $S$ be a spherical tangle diagram. Then

(1) $S^*$ is the mirror image of $S$,

(2) $S^-$ is the spherical tangle diagram obtained by
interchanging the inside hole with the outside hole of $S$,

(3) $S^{r_1}$ is the spherical tangle diagram obtained by only
rotating inside hole of $S$ $90^{\circ}$ counterclockwise on the
projection plane,

(4) $S^{r_2}$ is the spherical tangle diagram obtained by only
rotating outside hole of $S$ $90^{\circ}$ counterclockwise on the
projection plane,

(5) $S^R$ is the spherical tangle diagram obtained by the
$90^{\circ}$ rotation of $S$ itself counterclockwise on the
projection plane.
\end{df}

\bigskip
\centerline{\epsfxsize=4.8 in \epsfbox{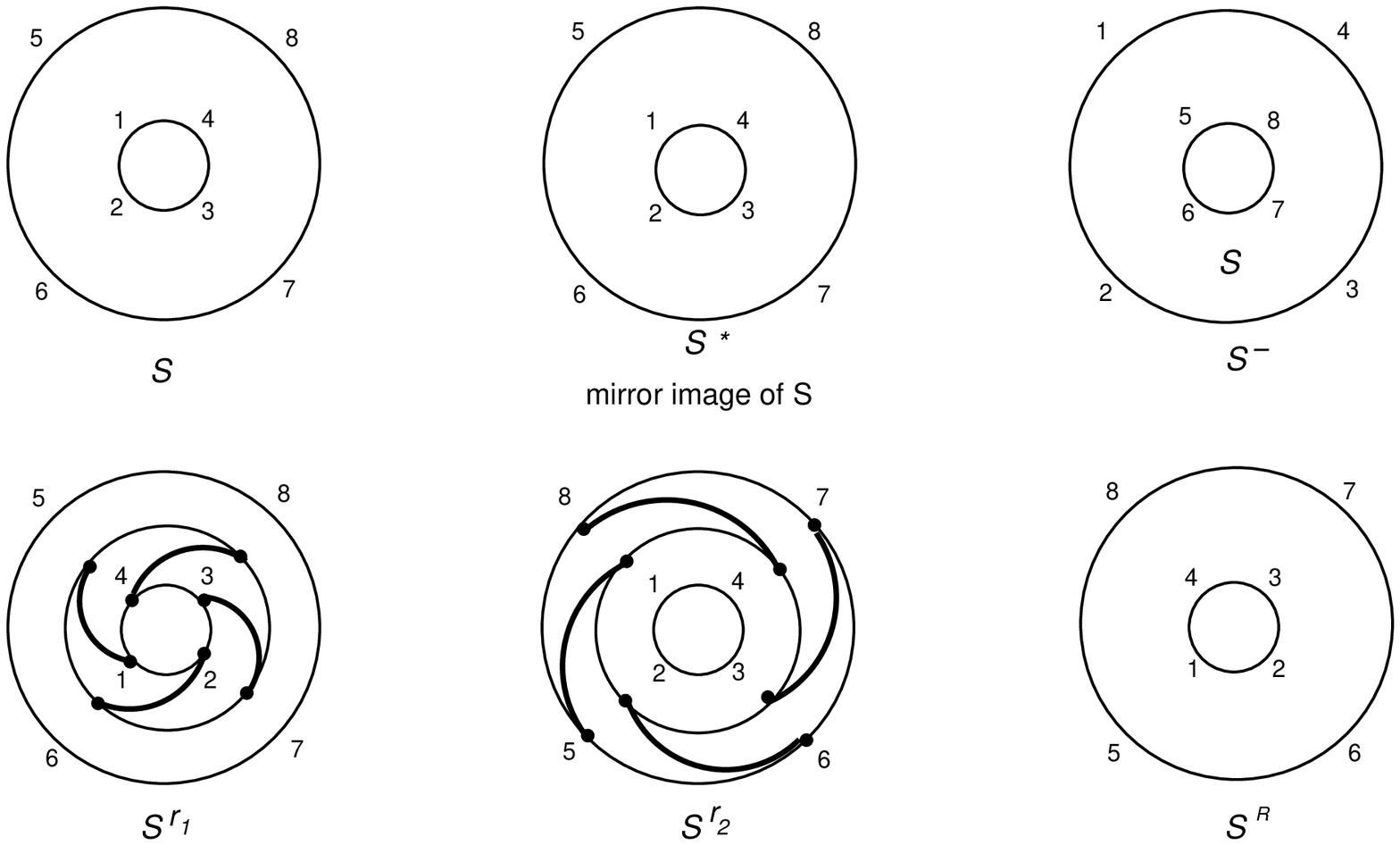}}
\medskip
\centerline{\small Figure 9. Elementary operations on
$\textbf{\textit{ST}}$.}
\bigskip

Note that $S^{r_2}=S^{-r_1-}$, $S^{r_1}=S^{-r_2-}$, and
$S^R=S^{r_1r_2}=S^{r_2r_1}$ for each $S \in \textbf{\textit{ST}}$.

\begin{lm} If $S \in \textbf{\textit{ST}}$ with the invariant
$F(S)=\left[\begin{matrix} \alpha & \gamma \\ \beta & \delta
\end{matrix}\right]$, then

(1) $F(S^*)=\left[\begin{matrix} \alpha & -\gamma \\
-\beta & \delta \end{matrix}\right]$,
(2) $F(S^-)=\left[\begin{matrix} \delta & \gamma \\
\beta & \alpha \end{matrix}\right]$,
(3) $F(S^{r_1})=\left[\begin{matrix} -\gamma & \alpha \\
-\delta & \beta \end{matrix}\right]$,

(4) $F(S^{r_2})=\left[\begin{matrix} -\beta & -\delta \\
\alpha & \gamma \end{matrix}\right]$,
(5) $F(S^R)=\left[\begin{matrix} \delta & -\beta \\
-\gamma & \alpha \end{matrix}\right]$.
\end{lm}

\begin{proof} Let $S \in \textbf{\textit{ST}}$ with
$F(S)=\left[\begin{matrix} \alpha & \gamma \\ \beta & \delta
\end{matrix}\right]$. Then there is $u \in \Phi$ such that
$\langle S_{11} \rangle=\alpha u$, $\langle S_{12} \rangle=\gamma
iu$, $\langle S_{21} \rangle=\beta (-i)u$, $\langle S_{22}
\rangle=\delta u$. Here the link $S_{ij}$, $i,j\in\{1,2\}$, is
obtained by taking the numerator closure ($i=1$) or the
denominator closure ($i=2$) of $S$ with its hole filled by the
fundamental tangle $j$. Therefore, $$\left[\begin{matrix} \alpha &
\gamma
\\ \beta & \delta
\end{matrix}\right]=\left[\begin{matrix} u^{-1}\alpha u &
u^{-1}(-i)\gamma iu
\\ u^{-1}i\beta (-i)u & u^{-1}\delta u
\end{matrix}\right].$$

Now we have

(1) $\langle S_{11}^* \rangle=\alpha u^{-1}, \langle S_{12}^*
\rangle=\gamma (iu)^{-1}, \langle S_{21}^* \rangle=\beta
(-iu)^{-1}, \langle S_{22}^* \rangle=\delta u^{-1}$,

(2) $\langle S_{11}^- \rangle=\langle S_{22} \rangle, \langle
S_{12}^- \rangle=\langle S_{12} \rangle, \langle S_{21}^-
\rangle=\langle S_{21} \rangle, \langle S_{22}^- \rangle=\langle
S_{11} \rangle$,

(3) $\langle S_{11}^{r_1} \rangle=\langle S_{12} \rangle, \langle
S_{12}^{r_1} \rangle=\langle S_{11} \rangle, \langle S_{21}^{r_1}
\rangle=\langle S_{22} \rangle, \langle S_{22}^{r_1}
\rangle=\langle S_{21} \rangle$.

Hence, $F(S^*)=\left[\begin{matrix} \alpha & -\gamma \\
-\beta & \delta \end{matrix}\right]$,
$F(S^-)=\left[\begin{matrix} \delta & \gamma \\
\beta & \alpha \end{matrix}\right]$,
$F(S^{r_1})=\left[\begin{matrix} \gamma & -\alpha \\
\delta & -\beta \end{matrix}\right]=
\left[\begin{matrix} -\gamma & \alpha \\
-\delta & \beta \end{matrix}\right]$.

Since $S^{r_2}=S^{-r_1-}$ and $S^R=S^{r_1r_2}$, (4) and (5) are
easily proved by (2) and (3).
\end{proof}

Like the case of ball tangle operations and invariants, it is
convenient to use the following notations.

Notation: Let $\left[\begin{matrix} \alpha & \gamma
\\ \beta & \delta
\end{matrix}\right] \in PM_{2\times2}$. Then

(1) $\left[\begin{matrix} \alpha & \gamma
\\ \beta & \delta
\end{matrix}\right]^*=\left[\begin{matrix} \alpha & -\gamma
\\ -\beta & \delta
\end{matrix}\right]$,
(2) $\left[\begin{matrix} \alpha & \gamma
\\ \beta & \delta
\end{matrix}\right]^-=\left[\begin{matrix} \delta & \gamma
\\ \beta & \alpha
\end{matrix}\right]$,
(3) $\left[\begin{matrix} \alpha & \gamma
\\ \beta & \delta
\end{matrix}\right]^{r_1}=\left[\begin{matrix} -\gamma & \alpha \\
-\delta & \beta \end{matrix}\right]$,

(4) $\left[\begin{matrix} \alpha & \gamma
\\ \beta & \delta
\end{matrix}\right]^{r_2}=\left[\begin{matrix} -\beta & -\delta \\
\alpha & \gamma \end{matrix}\right]$, (5) $\left[\begin{matrix}
\alpha & \gamma
\\ \beta & \delta
\end{matrix}\right]^R=\left[\begin{matrix} \delta & -\beta \\
-\gamma & \alpha \end{matrix}\right]$.

With these notations, we can write: $F(S^*)=F(S)^*, F(S^-)=F(S)^-,
F(S^{r_1})=F(S)^{r_1}, F(S^{r_2})=F(S)^{r_2}, F(S^R)=F(S)^R$ if $S
\in \textbf{\textit{ST}}$.

The determinant function {\rm det} is well-defined on
$PM_{2\times2}$ since ${\rm det}\,(-A)=(-1)^2\,{\rm det}\,A$ for
each $A \in PM_{2\times2}$. Notice that the 5 elementary
operations on $\textbf{\textit{ST}}$ do not change the determinant
of invariants of spherical tangles.

\begin{lm} If $S_1, S_2 \in \textbf{\textit{ST}}$, then

(1) $(S_1 \circ S_2)^*=S_1^* \circ S_2^*$, (2) $(S_1 \circ
S_2)^-=S_2^- \circ S_1^-$, (3) $(S_1 \circ S_2)^{r_1}=S_1 \circ
S_2^{r_1}$,

(4) $(S_1 \circ S_2)^{r_2}=S_1^{r_2} \circ S_2$, (5) $(S_1 \circ
S_2)^R=S_1^R \circ S_2^R$.
\end{lm}

Notice that a spherical tangle has exactly 2 holes which are
inside and outside.

\begin{df} Let $B \in \textbf{\textit{BT}}$, and let
$S \in \textbf{\textit{ST}}$. Then

(1) the 1st and the 2nd outer horizontal connect sums of $B$ and
$S$ are the spherical tangle diagrams denoted by $B +_h S$ and $S
+_h B$, respectively,

(2) the 1st and the 2nd outer vertical connect sums of $B$ and $S$
are the spherical tangle diagrams denoted by $B +_v S$ and $S +_v
B$, respectively (See Figure 10).
\end{df}

We also define the connect sums at the inside hole by $+_h, +_v,
^{-}$ as follows.

\begin{df} Let $B \in \textbf{\textit{BT}}$, and let
$S \in \textbf{\textit{ST}}$. Then

(1) the 1st and the 2nd inner horizontal connect sums of $B$ and
$S$ are the spherical tangle diagrams $B \overline{+}_h S$ and $S
\overline{+}_h B$ defined by $B \overline{+}_h S=(S^- +_h
B^{h-})^-$ and $S \overline{+}_h B=(B^{h-} +_h S^-)^-$,
respectively,

(2) the 1st and the 2nd inner vertical connect sums of $B$ and $S$
are the spherical tangle diagrams $B \overline{+}_v S$ and $S
\overline{+}_v B$ defined by $B \overline{+}_v S=(S^- +_v
B^{v-})^-$ and $S \overline{+}_h B=(B^{v-} +_v S^-)^-$,
respectively,

where $B^{h-}$ and $B^{v-}$ are the $180^{\circ}$ rotation of $B$
with respect to the vertical axis and the horizontal axis of the
projection plane, respectively (See Figure 11).
\end{df}

Let us give definitions of monoid actions. This is just a
generalization of group actions.

\begin{df} Let $M$ be a monoid with the identity $e$, and let
$X$ be a nonempty set. Then

(1) a function $*_l:M \times X \rightarrow X$ is called a left
monoid action of $M$ on $X$ if $*_l(m_1m_2,x)=*_l(m_1,*_l(m_2,x))$
and $*_l(e,x)=x$ for all $m_1,m_2 \in M$ and $x \in X$,

(2) a function $*_r:X \times M \rightarrow X$ is called a right
monoid action of $M$ on $X$ if $*_r(x,m_1m_2)=*_r(*_r(x,m_1),m_2)$
and $*_r(x,e)=x$ for all $m_1,m_2 \in M$ and $x \in X$.
\end{df}

In this sense, the connect sums of diagrams of ball tangles and
spherical tangles are monoid actions on $\textbf{\textit{ST}}$.
Hence, we have 8 monoid actions on $\textbf{\textit{ST}}$ by
$\textbf{\textit{BT}}$ which are similar. Also, the composition on
$\textbf{\textit{ST}}$ induces a left monoid action and a right
monoid action on $\textbf{\textit{BT}}$. In particular, a monoid
action is onto like a group action because of identity.

\begin{lm} Let $B \in \textbf{\textit{BT}}$, and let
$S \in \textbf{\textit{ST}}$. Then

(1) $(B +_h S)^R=S^R +_v B^R$, (2) $(S +_h B)^R=B^R +_v S^R$,

(3) $(B +_v S)^R=B^R +_h S^R$, (4) $(S +_v B)^R=S^R +_h B^R$.
\end{lm}

Note that $B^{RRRR}=B$ for each $B \in \textbf{\textit{BT}}$ and
$S^{RRRR}=S$ for each $S \in \textbf{\textit{ST}}$.

\bigskip
\centerline{\epsfxsize=4.5 in \epsfbox{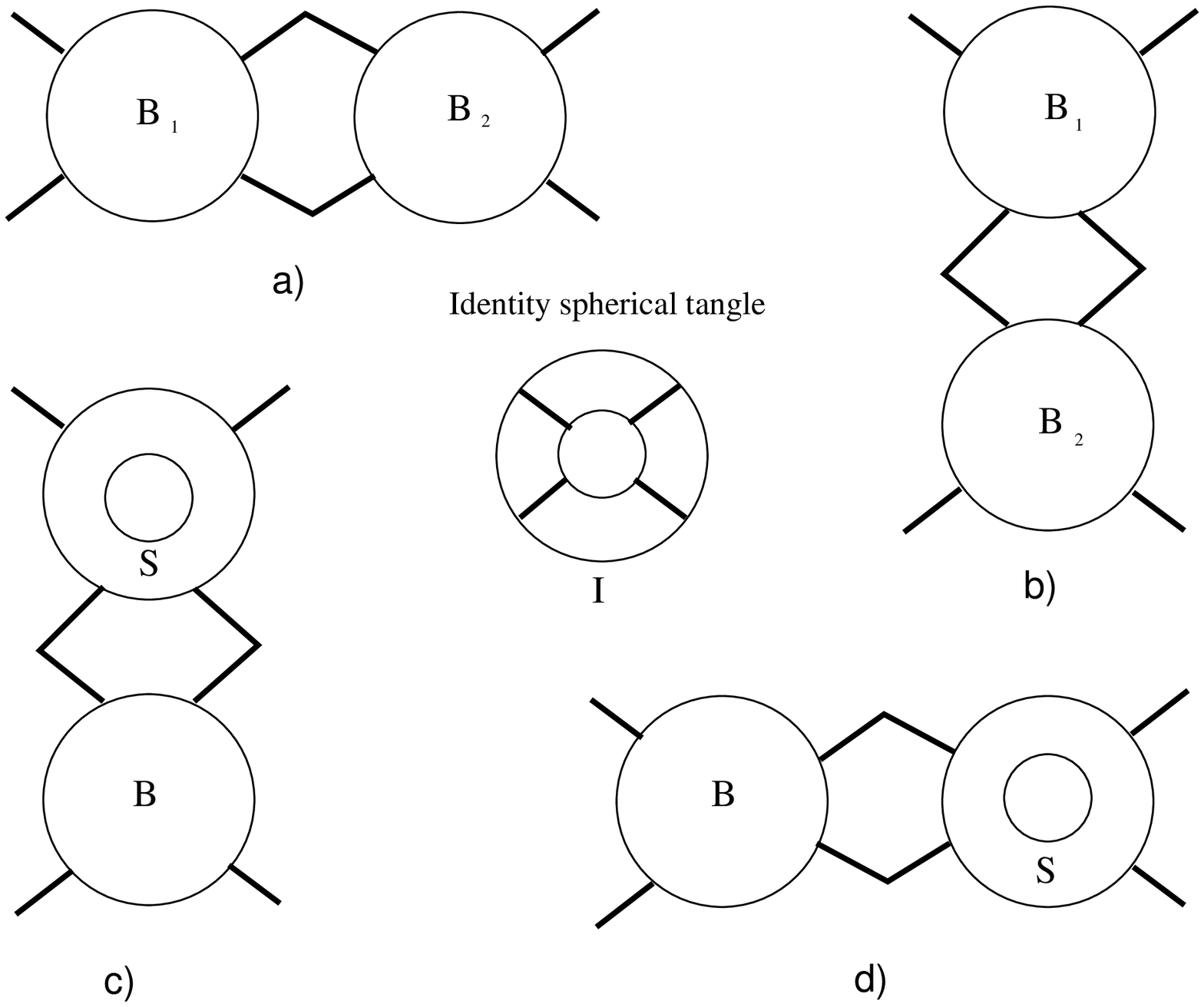}}
\medskip
\centerline{\small Figure 10. Connect sums of ball tangles and
outer connect sums.}
\bigskip

\vskip 0.2in

\bigskip
\centerline{\epsfxsize=4.5 in \epsfbox{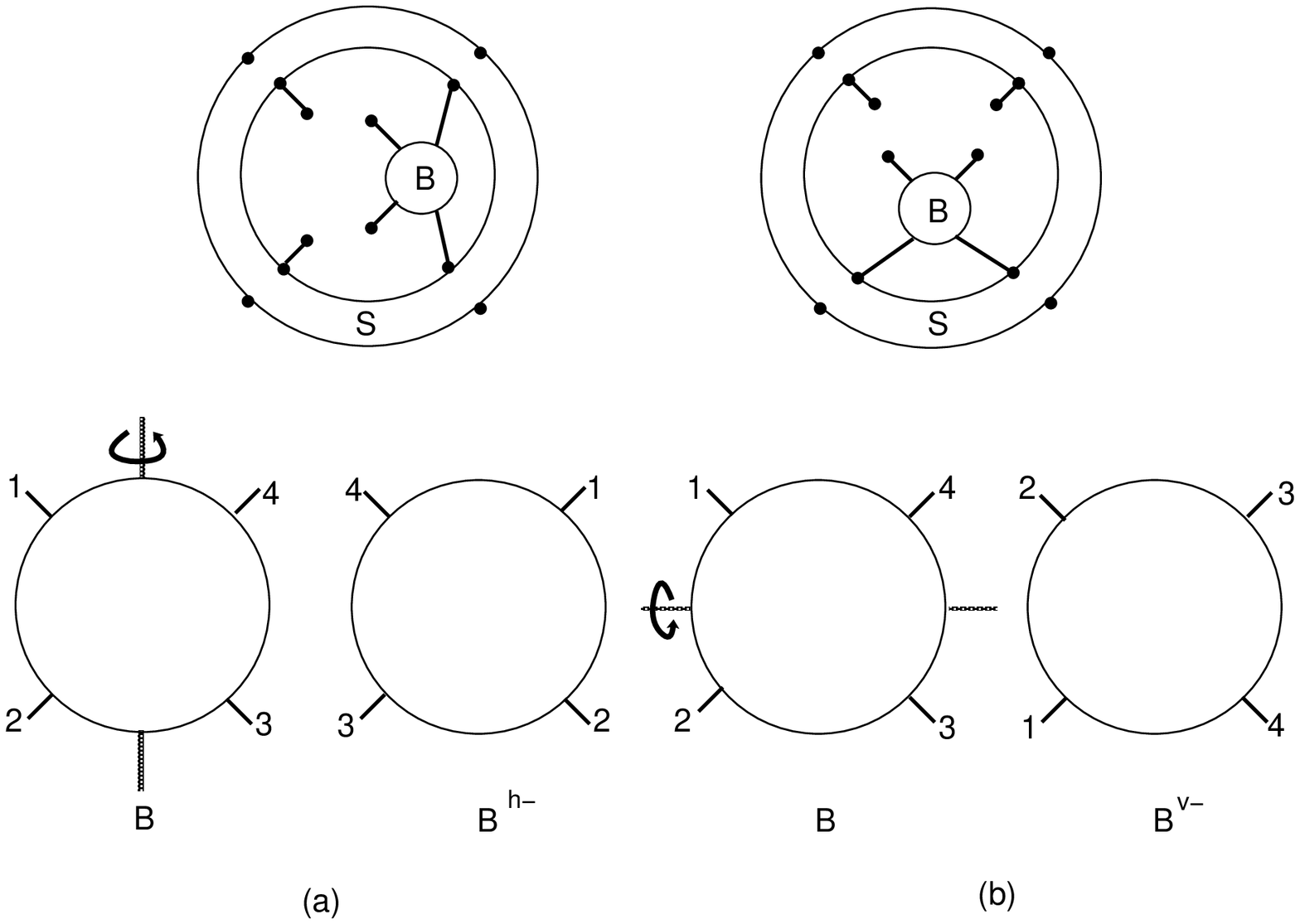}}
\medskip
\centerline{\small Figure 11. Inner connect sums and rotations of
ball tangles about axes.}
\bigskip

The following lemma tells us that the other outer horizontal sum
and two outer vertical sums can be expressed in terms of the 1st
horizontal sum and $^R$ which is the rotation for ball tangle
diagrams or spherical tangle diagrams.

\begin{lm} Let $B \in \textbf{\textit{BT}}$, and let
$S \in \textbf{\textit{ST}}$. Then

(1) $S +_h B=(B^{RR} +_h S^{RR})^{RR}$,

(2) $B +_v S=(B^R +_h S^R)^{RRR}$,

(3) $S +_v B=(B^{RRR} +_h S^{RRR})^R$.
\end{lm}

\begin{proof} Let $B \in \textbf{\textit{BT}}$, and let
$S \in \textbf{\textit{ST}}$. Then

(1) $S +_h B=(B' +_h S')^{RR}$ for some $B' \in
\textbf{\textit{BT}}, S' \in \textbf{\textit{ST}}$. Hence, $(B'
+_h S')^{RR}=(S'^R +_v B'^R)^R=S'^{RR} +_h B'^{RR}$. Take
$B'=B^{RR}$ and $S'=S^{RR}$. Then $S +_h B=(B^{RR} +_h
S^{RR})^{RR}$,

(2) $B +_v S=(B' +_h S')^{RRR}$ for some $B' \in
\textbf{\textit{BT}}, S' \in \textbf{\textit{ST}}$. Hence, $(B'
+_h S')^{RRR}=(S'^R +_v B'^R)^{RR}= (S'^{RR} +_h
B'^{RR})^R=B'^{RRR} +_v S'^{RRR}$. Take $B'=B^R$ and $S'=S^R$.
Then $B +_v S=(B^R +_h S^R)^{RRR}$,

(3) $S +_v B=(B'+_h S')^R$ for some $B' \in \textbf{\textit{BT}},
S' \in \textbf{\textit{ST}}$. Hence, $(B' +_h S')^R=S'^R +_v
B'^R$. Take $B'=B^{RRR}$ and $S'=S^{RRR}$. Then $S +_v B=(B^{RRR}
+_h S^{RRR})^R$.
\end{proof}

Now, we consider the invariants of spherical tangles obtained from
the various connect sums with ball tangles and their determinants.

Note that $f(B)=f(B^{h-})=f(B^{v-})$ for each $B \in
\textbf{\textit{BT}}$ and $S^{--}=S$ for each $S \in
\textbf{\textit{ST}}$. Also, $S^{**}=S$ for each $S \in
\textbf{\textit{ST}}$.

\begin{lm} If $B \in \textbf{\textit{BT}}$ with $f(B)=\left[\begin{matrix}
p \\ q \end{matrix}\right]$ and $S \in \textbf{\textit{ST}}$ with
$F(S)=\left[\begin{matrix} \alpha & \gamma \\ \beta & \delta
\end{matrix}\right]$, then

(1) $F(B +_h S)=F(S +_h B)=\left[\begin{matrix} p\beta + q\alpha &
p\delta + q\gamma \\ q\beta & q\delta
\end{matrix}\right]$, ${\rm det}\,F(B +_h S)=q^2\,{\rm det}\,F(S)$,

(2) $F(B +_v S)=F(S +_v B)=\left[\begin{matrix} p\alpha & p\gamma
\\ q\alpha + p\beta & q\gamma + p\delta
\end{matrix}\right]$, ${\rm det}\,F(B +_v S)=p^2\,{\rm det}\,F(S)$,

(3) $F(B \overline{+}_h S)=F(S \overline{+}_h
B)=\left[\begin{matrix} q\alpha & p\alpha + q\gamma \\ q\beta &
p\beta + q\delta
\end{matrix}\right]$, ${\rm det}\,F(B
\overline{+}_h S)=q^2\, {\rm det}\,F(S)$,

(4) $F(B \overline{+}_v S)=F(S \overline{+}_v
B)=\left[\begin{matrix} q\gamma + p\alpha & p\gamma
\\ q\delta + p\beta & p\delta
\end{matrix}\right]$, ${\rm det}\,F(B \overline{+}_v S)=p^2\,
{\rm det}\,F(S)$.
\end{lm}

\begin{proof} (1) Let $\left[\begin{matrix} x \\ y \end{matrix}\right] \in
PM_2$. Then there is $X \in \textbf{\textit{BT}}$ such that
$f(X)=\left[\begin{matrix} x \\ y \end{matrix}\right] \in PM_2$.
We have
$$
\begin{aligned}
&F(B +_h S)f(X)=f(B +_h S(X)) =\left[\begin{matrix} p \\ q
\end{matrix}\right] +_h \left[\begin{matrix}
\alpha & \gamma \\
\beta & \delta \end{matrix}\right] \left[\begin{matrix} x \\ y
\end{matrix}\right]\\
&=\left[\begin{matrix} p \\
q \end{matrix}\right] +_h \left[\begin{matrix} \alpha x + \gamma y \\
\beta x + \delta y \end{matrix}\right]
=\left[\begin{matrix} p\beta x + p\delta y + q\alpha x + q\gamma y \\
q\beta x + q\delta y \end{matrix}\right]\\
&=
\left[\begin{matrix} p\beta + q\alpha & p\delta + q\gamma \\
q\beta & q\delta
\end{matrix}\right] \left[\begin{matrix} x \\ y \end{matrix}\right]
=\left[\begin{matrix} p\beta + q\alpha & p\delta + q\gamma \\
q\beta & q\delta
\end{matrix}\right]f(X).
\end{aligned}$$
By Lemma 4.9, $F(B +_h S)=
\left[\begin{matrix} p\beta + q\alpha & p\delta + q\gamma \\
q\beta & q\delta \end{matrix}\right]$. Hence, ${\rm det}\,F(B +_h
S)=q^2\,{\rm det}\,F(S)$.

Also, $F(S +_h B)f(X)=f(S(X) +_h B)=f(B +_h S(X))=F(B +_h S)f(X)$.
Therefore, $F(S +_h B)=F(B +_h S)$. This proves (1).

(3) Since $B \overline{+}_h S=(S^- +_h B^{h-})^-$,
$F(S^-)=\left[\begin{matrix} \delta & \gamma \\ \beta & \alpha
\end{matrix}\right]$, and $f(B^{h-})=f(B)$, we have
$$
\begin{aligned}
&F(B \overline{+}_h S)=F((S^- +_h B^{h-})^-)=F(S^- +_h
B^{h-})^-\\
&=F(B^{h-} +_h S^-)^-=F((B^{h-} +_h S^-)^-)=F(S \overline{+}_h B).
\end{aligned}$$

Since $F(S^- +_h B^{h-})=\left[\begin{matrix} p\beta + q\delta &
p\alpha + q\gamma \\ q\beta & q\alpha
\end{matrix}\right]$, we have
$$F(B \overline{+}_h S)=\left[\begin{matrix} q\alpha &
p\alpha + q\gamma \\ q\beta & p\beta + q\delta
\end{matrix}\right]$$
and ${\rm det}\,F(B \overline{+}_h S)=q^2\,{\rm det}\,F(S)$. This
proves (3).

Similarly, (2) and (4) can be proved.
\end{proof}

\begin{df} A spherical tangle diagram $S$ is said to be
$\textbf{\textit{I}}$-reducible if there are $n \in\mathbb{N}$,
$A_1,\dots,A_{n+1} \in \textbf{\textit{BT}} \cup
\textbf{\textit{ST}}$, $*_1,\dots,*_n \in
\{+_h,+_h^{op},+_v,+_v^{op},\overline{+}_h,\overline{+}_h^{op},
\overline{+}_v,\overline{+}_v^{op}\}$ such that
$S=(\cdots(A_1*_1A_2)*_2\cdots A_n)*_nA_{n+1}$ and the only one of
$A_1,\dots,A_{n+1}$ is equivalent to $\textbf{\textit{I}}$ after
some possible operations of $(\cdot)^{r_1}$ and/or
$(\cdot)^{r_2}$, where $A +_h^{op} B = B +_h A$, $A +_v^{op} B = B
+_v A$, $A \overline{+}_h^{op} B = B \overline{+}_h A$, $A
\overline{+}_v^{op} B = B \overline{+}_v A$. In general, a
spherical tangle $S$ is $J$-reducible if, in above definition, we
replace $\textbf{\textit{I}}$ by a spherical tangle $J$.
\end{df}

In other words, a spherical tangle diagram $S$ is
$\textbf{\textit{I}}$-reducible if $S$ can be decomposed by
finitely many ball tangle diagrams and only one identity spherical
tangle diagram with respect to the inner and the outer connect
sums and their opposite operations. Note that, in Definition 4.19,
the expression of $S$ can be written as
$S=A_1*_1\cdots*_nA_{n+1}$. In this case, the order of operations
in the expression is important.

\begin{thm} If a spherical tangle $S$ is
$\textbf{\textit{I}}$-reducible, then ${\rm det}\,F(S)=n^2$ for
some $n \in \mathbb{Z}$. Furthermore, if $S$ is $J$-reducible,
then ${\rm det}\,F(S)=n^2\,{\rm det}\,F(J)$ for some $n \in
\mathbb{Z}$.
\end{thm}

\begin{proof} It follows from Lemma 4.18 immediately.
By Lemma 4.18, we know that a ball tangle connected to a spherical
tangle in the sense of Definition 4.13 and Definition 4.14
contributes a square of integer to the determinant of the
invariant of connect sum.
\end{proof}

Like rational ball tangles, we can define rational spherical
tangles, which will be $\textbf{\textit{I}}$-reducible.

\begin{df} A spherical tangle $S$ is said to be rational if $S$
can be obtained from the identity spherical tangle
$\textbf{\textit{I}}$ by moving boundary points of
$\textbf{\textit{I}}$ only on the spheres. In other words, A
spherical tangle $S$ is rational if the 1-submanifold $S$ in $S^2
\times [0,1]$ is isotopic to the identity spherical tangle
$\textbf{\textit{I}}$ provided that we are allowed to move the
boundary points of $S$ only on $S^2 \times \{0,1\}$.
\end{df}

\begin{lm} If a spherical tangle $S$ is rational, then
$S$ is $\textbf{\textit{I}}$-reducible. Therefore, ${\rm
det}\,F(S)$ is a square of integer.
\end{lm}

\begin{proof} Suppose that a spherical tangle $S$ is isotopic
to $\textbf{\textit{I}}$ by allowing points on $\partial S$ to
move. Then $S$ can be obtained from $\textbf{\textit{I}}$ by the
outer connect sums and the inner connect sums of finitely many
horizontal twists and vertical twists. Hence, $S$ can be written
as $S=\textbf{\textit{I}}*_0A_1*_1\dots*_kA_{k+1}$, where
$A_1,\dots,A_{k+1}$ are horizontal twists or vertical twists and
$*_0,*_1,\dots,*_{k+1}$ are connect sums of ball tangle diagrams
or outer connect sums or the inner connect sums or their opposite
operations. Therefore, $S$ is $\textbf{\textit{I}}$-reducible.
Hence, ${\rm det}\,F(S)$ is a square of integer.
\end{proof}

Some examples of spherical tangles are given in Figure 12.

\bigskip
\centerline{\epsfxsize=4.5 in \epsfbox{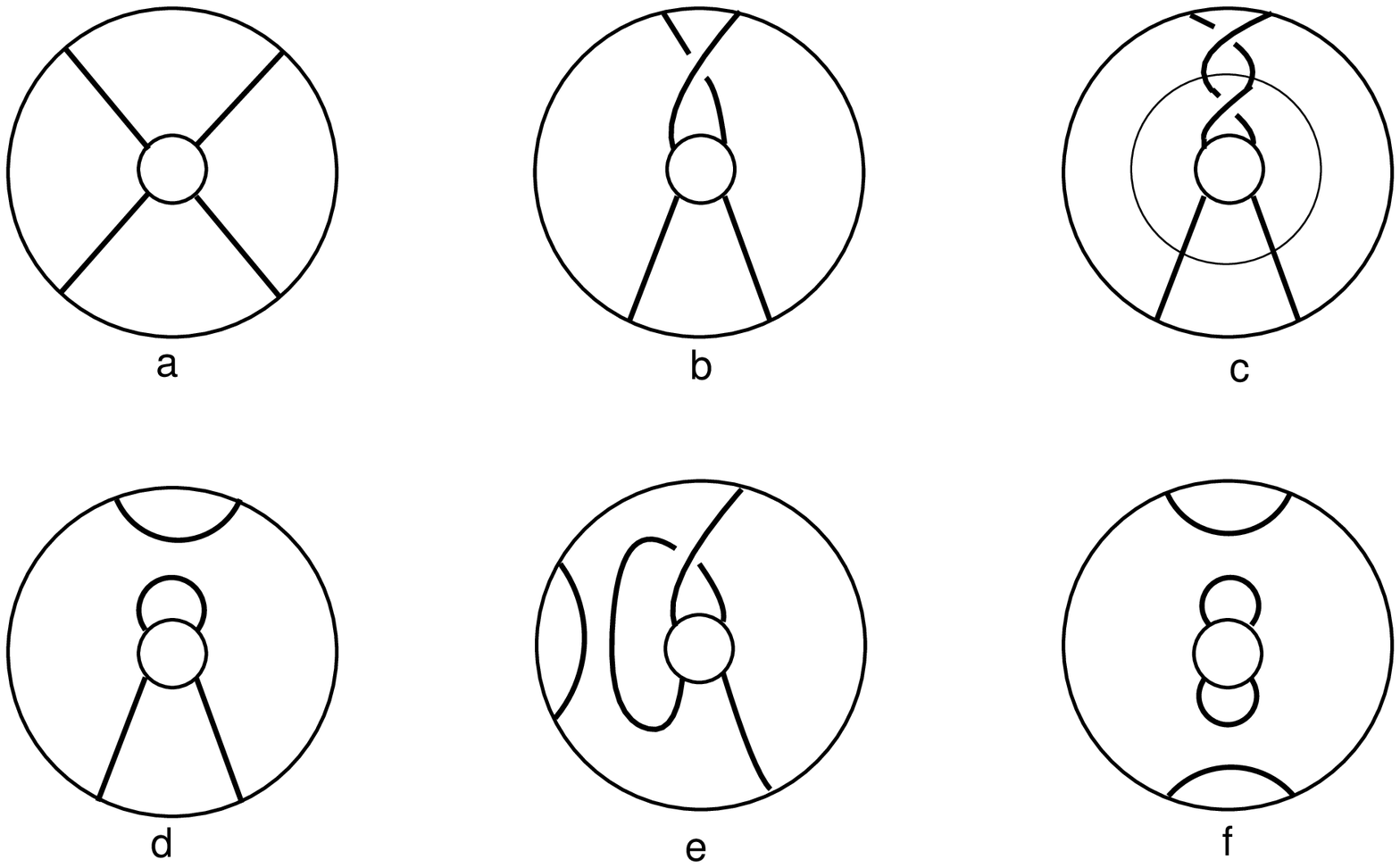}}
\medskip
\centerline{\small Figure 12. Spherical tangle diagrams.}
\bigskip

1. The spherical tangle \textbf{\textit{a}} is
$\textbf{\textit{I}}$ and has invariant $\left[\begin{matrix} 1 & 0 \\
0 & 1 \end{matrix}\right]$.

2. The spherical tangle \textbf{\textit{b}} has invariant
$\left[\begin{matrix} 1 & 0 \\ 1 & 1 \end{matrix}\right]$.

3. The spherical tangle \textbf{\textit{c}} is
$\textbf{\textit{b}} \circ \textbf{\textit{b}}$ and has invariant
$\left[\begin{matrix} 1 & 0 \\ 1 & 1 \end{matrix}\right]
\left[\begin{matrix} 1 & 0 \\ 1 & 1 \end{matrix}\right]=
\left[\begin{matrix} 1 & 0 \\ 2 & 1 \end{matrix}\right]$.

4. The spherical tangle \textbf{\textit{d}} has invariant
$\left[\begin{matrix} 0 & 0 \\ 1 & 0 \end{matrix}\right]$.

5. The spherical tangle \textbf{\textit{e}} has invariant
$\left[\begin{matrix} 1 \\ 0 \end{matrix}\right] +_{h_1}
\left[\begin{matrix} 1 & 0 \\ 1 & 1 \end{matrix}\right] =
\left[\begin{matrix} 1 & 1 \\ 0 & 0
\end{matrix}\right]$. (See lemma 4.18)

6. The spherical tangle \textbf{\textit{f}} has invariant
$\left[\begin{matrix} 0 & 0 \\ 0 & 0 \end{matrix}\right]$.

7. The spherical tangle $\textbf{\textit{b}} \circ
\textbf{\textit{e}}$ has invariant $\left[\begin{matrix} 1 & 0 \\
1 & 1 \end{matrix}\right] \left[\begin{matrix} 1 & 1
\\ 0 & 0 \end{matrix}\right] = \left[\begin{matrix}
1 & 1 \\ 1 & 1 \end{matrix}\right]$ and $(\textbf{\textit{b}}
\circ \textbf{\textit{e}})^{r_1}$ has invariant
$\left[\begin{matrix} 1 & -1 \\ 1 & -1 \end{matrix}\right]$.

As another example, let us consider the spherical tangle diagram
$S$ in Figure 13. By the definition of invariant, we have
$$F(S)=\{\begin{pmatrix} z5A^3 & -iz(-8A^5) \\ iz8A & z(-11A^3)
\end{pmatrix}|z \in \Phi\} \cap M_{2\times2}(Z)=\left[\begin{matrix}
5 & -8 \\ 8 & -11 \end{matrix}\right].$$ Hence, ${\rm
det}\,F(S)=-55-(-64)=9$. That is, ${\rm det}\,F(S)$ is a square of
integer. However, $S$ looks like a
non-$\textbf{\textit{I}}$-reducible spherical tangle although we
are not able to prove this fact.

\bigskip
\centerline{\epsfxsize=4.5 in \epsfbox{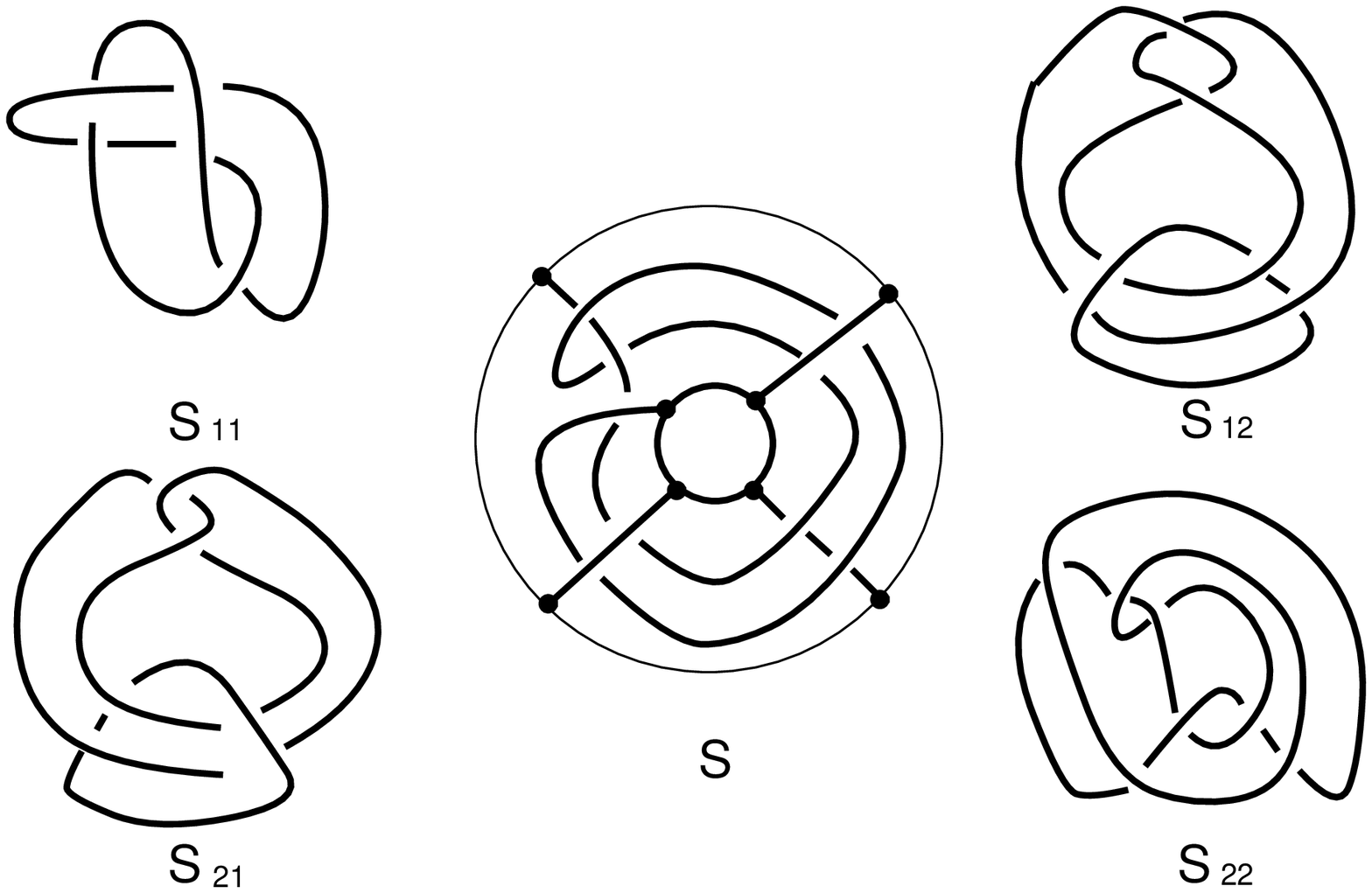}}
\medskip
\centerline{\small Figure 13. A spherical tangle $S$ which cannot
be decomposed.}
\bigskip

\section{The Elementary operations on $PM_{2\times2}$ and Coxeter groups}

In this section, we introduce the group structure generated by the
elementary operations on $PM_{2\times2}$ induced by the elementary
operations on $\textbf{\textit{ST}}$.

\begin{df} An $n \times n$ matrix $M$ is called a Coxeter matrix
if $M_{ii}=1$ and $M_{ij}=M_{ji}>1$ for all $i,j \in
\{1,\dots,n\}$ with $i \neq j$, where $M_{ij}$ is the
$(i,j)$-entry of $M$.
\end{df}

\begin{df} Let $M$ be an $n \times n$ Coxeter matrix. Then a group
presented by $$\langle\, x_1,\dots,x_n\, |\,
(x_ix_j)^{M_{ij}}=1\,\, {\rm for}\,\, {\rm all}\,\, i,j \in
\{1,\dots,n\}\, \rangle,$$ denoted by $C_M$, is called the Coxeter
group with the Coxeter matrix $M$.
\end{df}

Let us think of the $5$ elementary operations $*, -, r_1, r_2, R$
on $PM_{2\times2}$ induced by the elementary operations on
$\textbf{\textit{ST}}$ as functions from $PM_{2\times2}$ to
$PM_{2\times2}$, respectively. For convenience, we use the
opposite composition of functions for the binary operation. For
instance, $-r_1$ means the composition $r_1 \circ -$. Recall that
$S^{r_2}=S^{-r_1-}$, $S^{r_1}=S^{-r_2-}$, and
$S^R=S^{r_1r_2}=S^{r_2r_1}$ for each $S \in \textbf{\textit{ST}}$
and observe the followings:

Suppose that $\left[\begin{matrix} \alpha & \gamma
\\ \beta & \delta \end{matrix}\right] \in PM_{2\times2}$. Then
$$\left[\begin{matrix} \alpha & \gamma
\\ \beta & \delta \end{matrix}\right] \underrightarrow{-}
\left[\begin{matrix} \delta & \gamma
\\ \beta & \alpha \end{matrix}\right] \underrightarrow{r_1}
\left[\begin{matrix} -\gamma & \delta
\\ -\alpha & \beta \end{matrix}\right] \underrightarrow{-}
\left[\begin{matrix} \beta & \delta
\\ -\alpha & -\gamma \end{matrix}\right] \underrightarrow{r_1}
\left[\begin{matrix} -\delta & \beta
\\ \gamma & -\alpha \end{matrix}\right],$$
$$\left[\begin{matrix} \alpha & \gamma
\\ \beta & \delta \end{matrix}\right] \underrightarrow{r_1}
\left[\begin{matrix} -\gamma & \alpha
\\ -\delta & \beta \end{matrix}\right] \underrightarrow{-}
\left[\begin{matrix} \beta & \alpha
\\ -\delta & -\gamma \end{matrix}\right] \underrightarrow{r_1}
\left[\begin{matrix} -\alpha & \beta
\\ \gamma & -\delta \end{matrix}\right] \underrightarrow{-}
\left[\begin{matrix} -\delta & \beta
\\ \gamma & -\alpha \end{matrix}\right]$$ and
$$\left[\begin{matrix} \alpha & \gamma
\\ \beta & \delta \end{matrix}\right] \underrightarrow{-}
\left[\begin{matrix} \delta & \gamma
\\ \beta & \alpha \end{matrix}\right] \underrightarrow{*}
\left[\begin{matrix} \delta & -\gamma
\\ -\beta & \alpha \end{matrix}\right],$$
$$\left[\begin{matrix} \alpha & \gamma
\\ \beta & \delta \end{matrix}\right] \underrightarrow{*}
\left[\begin{matrix} \alpha & -\gamma
\\ -\beta & \delta \end{matrix}\right] \underrightarrow{-}
\left[\begin{matrix} \delta & -\gamma
\\ -\beta & \alpha \end{matrix}\right]$$ and
$$\left[\begin{matrix} \alpha & \gamma
\\ \beta & \delta \end{matrix}\right] \underrightarrow{r_1}
\left[\begin{matrix} -\gamma & \alpha
\\ -\delta & \beta \end{matrix}\right] \underrightarrow{*}
\left[\begin{matrix} -\gamma & -\alpha
\\ \delta & \beta \end{matrix}\right],$$
$$\left[\begin{matrix} \alpha & \gamma
\\ \beta & \delta \end{matrix}\right] \underrightarrow{*}
\left[\begin{matrix} \alpha & -\gamma
\\ -\beta & \delta \end{matrix}\right] \underrightarrow{r_1}
\left[\begin{matrix} \gamma & \alpha
\\ -\delta & -\beta \end{matrix}\right].$$
Hence, we have $-r_1-r_1=r_1-r_1-$ and $-*=*-$ and $r_1*=*r_1$.
Also, $--$ and $r_1r_1$ and $**$ are the identity function from
$PM_{2\times2}$ to $PM_{2\times2}$. We show that the group
generated by the elementary operations on $PM_{2\times2}$ induced
by those on $\textbf{\textit{ST}}$ has the group presentation
$\langle\, x,y,z\, |\, x^2=y^2=z^2=1,\, xyxy=yxyx,\, xz=zx,\,
yz=zy\, \rangle$ which is a Coxeter group.

\begin{thm} The group $G(F)$ generated by the
elementary operations on $PM_{2\times2}$ induced by those on
$\textbf{\textit{ST}}$ has the group presentation $$\langle\,
x,y,z\, |\, x^2=y^2=z^2=1,\, xyxy=yxyx,\, xz=zx,\, yz=zy\,
\rangle.$$ Furthermore, $G(F)$ is isomorphic to the Coxeter group
$C_M$ with the Coxeter matrix $M=\begin{pmatrix} 1 & 4 & 2 \\
4 & 1 & 2 \\ 2 & 2 & 1 \end{pmatrix}$. That is, $$G(F)=\langle\,
x,y,z\, |\,
x^2=y^2=z^2=(xy)^4=(yx)^4=(xz)^2=(zx)^2=(yz)^2=(zy)^2=1\,
\rangle.$$
\end{thm}

\begin{proof} Let $G=\langle\, x,y,z\, |\, x^2=y^2=z^2=1,\,
xyxy=yxyx,\, xz=zx,\, yz=zy\, \rangle$. Suppose that $\phi:G
\rightarrow G(F)$ is the epimorphism such that $\phi(x)=-$,
$\phi(y)=r_1$, $\phi(z)=*$. We claim that ${\rm Ker}\,\phi=\{1\}$.
Let $W(x,y,z)$ be a word in ${\rm Ker}\,\phi$. Then
$\phi(W(x,y,z))=W(-,r_1,*)=Id_{PM_{2\times2}}$. Since
$x^2=y^2=z^2=1$, we may assume that $W(x,y,z)$ has no consecutive
letters and no inverses of letters. Since $xz=zx$ and $yz=zy$ and
$z^2=1$, we have either $W(x,y,z)=W_1(x,y)z$ or
$W(x,y,z)=W_1(x,y)$ for some word $W_1(x,y)$ in $\{x,y\}$. We may
also assume that $W_1(x,y)$ has no consecutive letters and no
inverses of letters. We show that $W(x,y,z) \neq W_1(x,y)z$. If
$W(x,y,z)=W_1(x,y)z$, then $W_1(-,r_1)*=Id_{PM_{2\times2}}$. That
is, $W_1(-,r_1)=*$.

Observe that $$\begin{matrix} - \neq *, & -r_1 \neq *, & -r_1-
\neq *, & -r_1-r_1 \neq *, \\ r_1 \neq *, & r_1- \neq *, & r_1-r_1
\neq *, & r_1-r_1- \neq *. \end{matrix}$$ By $-r_1-r_1=r_1-r_1-$
and $-^2=r_1^2=Id_{PM_{2\times2}}$, we have $W_1(-,r_1) \neq *$.
This is a contradiction. Hence, $W(x,y,z) \neq W_1(x,y)z$.
Therefore, $W(x,y,z)=W_1(x,y)$ and the number of $z$ in $W(x,y,z)$
must be even.

Since $W_1(x,y)$ has no consecutive letters and no inverses of
letters, we have either

there are $k \in \mathbb N \cup \{0\}$ and $R \in
\{1,x,xy,xyx,xyxy,xyxyx,xyxyxy,xyxyxyx\}$

such that $W_1(x,y)=(xy)^{4k}R$ or

there are $k' \in \mathbb N \cup \{0\}$ and $R' \in
\{1,y,yx,yxy,yxyx,yxyxy,yxyxyx,yxyxyxy\}$

such that $W_1(x,y)=(yx)^{4k'}R'$.

Also, since $W_1(x,y)=W(x,y,z) \in {\rm Ker}\,\phi$, we have
either

$Id_{PM_{2\times2}}=W_1(-,r_1)=(-r_1)^{4k}\phi(R)$ or
$Id_{PM_{2\times2}}=W_1(-,r_1)=(r_1-)^{4k'}\phi(R')$.

Similarly, as above, observe that
$$\begin{matrix} - \neq
Id_{PM_{2\times2}}, & -r_1 \neq Id_{PM_{2\times2}}, & -r_1-
\neq Id_{PM_{2\times2}}, & -r_1-r_1 \neq Id_{PM_{2\times2}}, \\
r_1 \neq Id_{PM_{2\times2}}, & r_1- \neq Id_{PM_{2\times2}}, &
r_1-r_1 \neq Id_{PM_{2\times2}}, & r_1-r_1- \neq
Id_{PM_{2\times2}}.
\end{matrix}$$

Also, notice that $$\begin{matrix} -r_1-r_1-=r_1-r_1, &
-r_1-r_1-r_1=r_1-, \\ -r_1-r_1-r_1-=r_1, &
-r_1-r_1-r_1-r_1=Id_{PM_{2\times2}} \end{matrix}$$ and
$$\begin{matrix} r_1-r_1-r_1=-r_1-, & r_1-r_1-r_1-=-r_1, \\
r_1-r_1-r_1-r_1=-, & r_1-r_1-r_1-r_1-=Id_{PM_{2\times2}}.
\end{matrix}$$ Hence, we know that $\phi(R)=Id_{PM_{2\times2}}$ if and only if $R=1$
and $\phi(R')=Id_{PM_{2\times2}}$ if and only if $R'=1$.

Since $(-r_1)^{4k}=Id_{PM_{2\times2}}$ and
$(r_1-)^{4k'}=Id_{PM_{2\times2}}$, we have
$\phi(R)=Id_{PM_{2\times2}}$ and $\phi(R')=Id_{PM_{2\times2}}$.
Hence, $R=1$ and $R'=1$.

Therefore, $W_1(x,y)=(xy)^{4k}$ or $W_1(x,y)=(yx)^{4k}$ for some
$k \in \mathbb N \cup \{0\}$. Since $(xy)^4=1$ and $(yx)^4=1$,
$W_1(x,y)=1$. That is, $W(x,y,z)=1$. We have proved ${\rm
Ker}\,\phi=\{1\}$. Hence, $\phi:G \rightarrow G(F)$ is a group
isomorphism and $G(F)$ has the group presentation $\langle\,
x,y,z\, |\, x^2=y^2=z^2=1,\, xyxy=yxyx,\, xz=zx,\, yz=zy\,
\rangle$.

Now, we show that $G(F)$ is isomorphic to $C_M$. Since
$(xy)^2=(yx)^2$, $(xy)^2(xy)^2=(yx)^2(xy)^2$ and
$(xy)^2(yx)^2=(yx)^2(yx)^2$. Since $x^2=y^2=1$, $(xy)^4=(yx)^4=1$.
Also, since $xz=zx$, $(xz)(xz)=(zx)(xz)$ and $(xz)(zx)=(zx)(zx)$.
Since $x^2=z^2=1$, $(xz)^2=(zx)^2=1$. Similarly, since $yz=zy$,
$(yz)(yz)=(zy)(yz)$ and $(yz)(zy)=(zy)(zy)$. Since $y^2=z^2=1$,
$(yz)^2=(zy)^2=1$. Hence, the consequence of relators of $C_M$ is
contained in that of $G(F)$. Conversely, Since $(xy)^4=1$,
$(xy)^4(yx)^2=(yx)^2$. Since $x^2=y^2=1$, $(xy)^2=(yx)^2$. Also,
since $(xz)^2=1$, $(xz)^2(zx)=zx$. Since $x^2=z^2=1$, $xz=zx$.
Similarly, since $(yz)^2=1$, $(yz)^2(zy)=zy$. Since $y^2=z^2=1$,
$yz=zy$. Hence, the consequence of relators of $G(F)$ is contained
in that of $C_M$. Thus, $G(F)$ is isomorphic to $C_M$.
\end{proof}

\section{Nonsurjectivity of the spherical tangle invariant $F$}

Recall the $\Delta$-move on knot diagrams introduced in
\cite{M-Y}. It is illustrated in Figure 14 (a). If we apply the
Kauffman states to the diagrams involved in the $\Delta$-move, we
get the 5 basis diagrams without closed components as shown in
Figure 14 (b).

\vskip 0.2in

\bigskip
\centerline{\epsfxsize=3.8 in \epsfbox{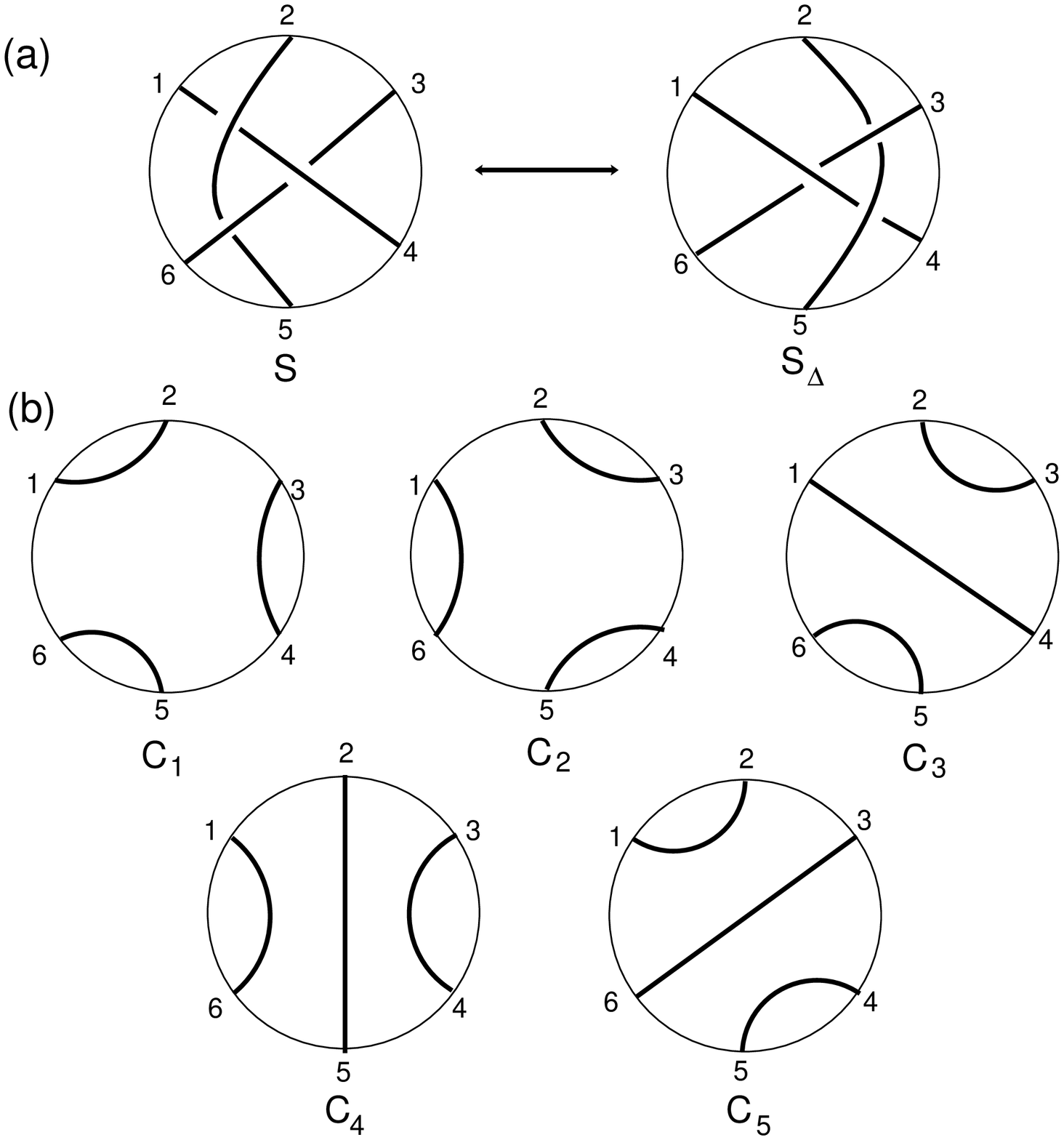}}
\medskip
\centerline{\small Figure 14. (a) A $\Delta$-move on a diagram;
(b) The 5 basis diagrams.}
\bigskip

\begin{df} Link diagrams $D_1$ and $D_2$ are said to be
$\Delta$-equivalent if $D_2$ can be obtained from $D_1$ by a
finite sequence of $\Delta$-moves and Reidemeister moves.
\end{df}

For a spherical tangle $S$, we can get four links
$S_{11},S_{12},S_{21},S_{22}$ by taking closures of $S$ with its
hole filled by fundamental tangles. We say that two spherical
tangles $S$ and $S'$ are $\Delta$-equivalent if each $S'_{ij}$ can
be obtained from $S_{ij}$ by a finite sequence of $\Delta$-moves
and Reidemeister moves.

\begin{lm} Let $S$ and $S'$ be spherical tangles such that $S$ and $S'$ are
$\Delta$-equivalent. If $F(S)=\left[\begin{matrix} \alpha & \gamma
\\ \beta & \delta
\end{matrix}\right]$ and $F(S')=\left[\begin{matrix} \alpha'
& \gamma' \\ \beta' & \delta' \end{matrix}\right]$, then $\alpha
\equiv \epsilon\alpha'$ {\rm mod} 4, $\beta \equiv \epsilon\beta'$
{\rm mod} 4, $\gamma \equiv \epsilon\gamma'$ {\rm mod} 4, $\delta
\equiv \epsilon\delta'$ {\rm mod} 4 for some $\epsilon=\pm 1$.
\end{lm}

\begin{proof} Suppose that $L$ is a link diagram and $L_\Delta$ is a link
diagram obtained from $L$ by a single $\Delta$-move. Then
$$
\begin{aligned}
\langle L \rangle &= A^3\langle C_2 \rangle + A\langle C_4 \rangle
+ A\langle C_3 \rangle + A^{-1}\langle C_1 \rangle + A\langle C_5
\rangle + A^{-1}\langle C_1 \rangle + A^{-1}\langle C_1 \rangle\\
&=3A^{-1}\langle C_1 \rangle + A^3\langle C_2 \rangle + A\langle
C_3 \rangle + A\langle C_4 \rangle + A\langle C_5 \rangle
\end{aligned}$$
and
$$
\begin{aligned}
\langle L_\Delta \rangle &= A^3\langle C_1 \rangle + A\langle C_3
\rangle + A\langle C_5 \rangle + A^{-1}\langle C_2 \rangle +
A\langle C_4 \rangle + A^{-1}\langle C_2 \rangle + A^{-1}\langle
L_2 \rangle\\
&=A^3\langle C_1 \rangle + 3A^{-1}\langle C_2 \rangle + A\langle
C_3 \rangle + A\langle C_4 \rangle + A\langle C_5 \rangle.
\end{aligned}$$
Hence, $$\langle L \rangle - \langle L_\Delta
\rangle=4A^{-1}\langle C_1 \rangle - 4A^{-1}\langle C_2 \rangle.$$

Suppose that $\langle L \rangle=aA^k$, $\langle L_\Delta
\rangle=a'A^{k'}$, $\langle C_1 \rangle=bA^l$, $\langle C_2
\rangle=b'A^{l'}$ for some $a,a',b,b',k,k',l,l' \in \mathbb{Z}$.
Then $aA^k - a'A^{k'}=4bA^{l-1} - 4b'A^{l'-1}$.

\newpage

{\it Case 1}. If $A^k=A^{k'}$, then $(a - a')A^k=4bA^{l-1} -
4b'A^{l'-1}$, hence,

$(a - a')A^k=0; -4b'A^{l'-1}; 4bA^{l-1}; 4(b \pm b')A^{l-1}$ if
$b=0, b'=0$; $b=0, b'\neq0$; $b\neq0, b'=0$; $b\neq0, b'\neq0$,
respectively, so $a - a' \equiv 0$ {\rm mod} 4.

{\it Case 2}. If $A^k=-A^{k'}$, then $(a + a')A^k=4bA^{l-1} -
4b'A^{l'-1}$, hence,

$(a + a')A^k=0; -4b'A^{l'-1}; 4bA^{l-1}; 4(b \pm b')A^{l-1}$ if
$b=0, b'=0$; $b=0, b'\neq0$; $b\neq0, b'=0$; $b\neq0, b'\neq0$,
respectively, so $a + a' \equiv 0$ {\rm mod} 4.

{\it Case 3}. If $A^k \neq \pm A^{k'}$, then $aA^k -
a'A^{k'}=4bA^{l-1} - 4b'A^{l'-1}$, hence,

$aA^k - a'A^{k'}=0; -4b'A^{l'-1}; 4bA^{l-1}; 4bA^{l-1} -
4b'A^{l'-1}$ if $b=0, b'=0$; $b=0, b'\neq0$; $b\neq0, b'=0$;
$b\neq0, b'\neq0$, respectively, so

1) $a=0, a'=0$ if $b=0, b'=0$,

2) $a=0, a'=\pm 4b'$ or $a=\pm 4b', a'=0$ if $b=0, b'\neq0$,

3) $a=0, a'=\pm 4b$ or $a=\pm 4b, a'=0$ if $b\neq0, b'=0$,

4) $a=\pm 4b, a'=\pm 4b'$ or $a=\pm 4b', a'=\pm 4b$ if $b\neq0,
b'\neq0$.

Hence, $a \equiv 0$ {\rm mod} 4 and $a' \equiv 0$ {\rm mod} 4.
Therefore, we always have $a \equiv \epsilon a'$ {\rm mod} 4 with
$\epsilon=\pm1$.

In general, suppose that a link $L'$ can be obtained from $L$ by a
finite sequence of $\Delta$-moves and Reidemeister moves of type
II and III. If $\langle L\rangle =aA^k$ and $\langle L'\rangle
=a'A^{k'}$, then $a\equiv\epsilon a'$ mod 4 and $\epsilon=\pm1$
depends only on the powers $k$ and $k'$.

For the spherical tangle $S$, we need to consider 4 links
$S_{11},S_{12},S_{21},S_{22}$. If $F(S) =\left[\begin{matrix}
\alpha & \gamma \\ \beta & \delta
\end{matrix}\right]$, then $\langle S_{11}\rangle=\alpha A^k$,
$\langle S_{12}\rangle=\gamma A^{k+2}$, $\langle S_{21}\rangle =
\beta A^{k-2}$, and $\langle S_{22}\rangle = \delta A^{k}$. Also,
if $F(S') =\left[\begin{matrix} \alpha' & \gamma' \\ \beta' &
\delta'
\end{matrix}\right]$, then $\langle S'_{11}\rangle=\alpha' A^{k'}$,
$\langle S'_{12}\rangle=\gamma' A^{k'+2}$, $\langle S'_{21}\rangle
= \beta' A^{k'-2}$, and $\langle S'_{22}\rangle = \delta' A^{k'}$.
Notice that since a $\Delta$-move will not change the writhe and
we can postpone all Reidemeister moves of type I in any finite
sequence of diagram moves to the end of that sequence, we do not
need to worry that the Kauffman bracket is only a regular isotopy
invariant.

Thus, for the four corresponding links we obtained from $S'$, the
sign $\epsilon$ is a constant. Thus, we have $\alpha \equiv
\epsilon\alpha'$ {\rm mod} 4, $\beta \equiv \epsilon\beta'$ {\rm
mod} 4, $\gamma \equiv \epsilon\gamma'$ {\rm mod} 4, $\delta
\equiv \epsilon\delta'$ {\rm mod} 4.
\end{proof}

We will make use of the following theorem.

\begin{thm} {\rm (Matveev \cite{Mat} and Murakami-Nakanishi \cite {M-Y})}:
Oriented links $L=L_1 \sqcup \cdots \sqcup L_n$ and $L'=L'_1
\sqcup \cdots \sqcup L'_n$ are $\Delta$-equivalent if and only if
${\rm lk}(L_i,L_j)={\rm lk}(L'_i, L'_j)$ for all $i,j$ with $1
\leq i < j \leq n$.
\end{thm}

Suppose that a spherical tangle diagram $S$ has no circle
components. It has 4 components $K_1,K_2,K_3,K_4$. We will look at
a diagram of $S$ and orient each $K_i$ arbitrarily. We use $-K_i$
to mean to reverse the orientation of $K_i$. We define ${\rm
lk}(K_i,K_j)$ to be a half of the sum of the signs of crossings
between $K_i$ and $K_j$. We have ${\rm lk}(-K_i,K_j)=-{\rm
lk}(K_i, K_j)$.

Note: For each $S_{11},S_{12},S_{21},S_{22}$, we get a link whose
components are unions of some of $K_1,K_2,K_3,K_4$. If components
of $S_{11}$ etc. are oriented, they are unions of some of
$\epsilon K_1$, $ \epsilon K_2$, $\epsilon K_3$, $\epsilon K_4$,
where $\epsilon_i=\pm1$.

Let $S'$ be another spherical tangle with components
$K'_1,K'_2,K'_3,K'_4$. Suppose the end points of $K'_i$ are the
same as the end points of $K_i$. i.e., $\partial K'_1=\partial
K_1$, $\partial K'_2=\partial K_2$, $\partial K'_3=\partial K_3$,
$\partial K'_4=\partial K_4$. So we can orient each $K_i$ and
$K'_i$ consistently. In this case, we can orient the links
$S_{ij}$ and $S'_{ij}$ consistently in the sense that the
corresponding components of $S_{ij}$ and $S'_{ij}$ are the same
union of $\epsilon_i K_i$ and $\epsilon_iK'_i$, respectively.

\begin{lm} The linking numbers of $S_{ij}$ and $S'_{ij}$ are equal
for all $i,j \in\{1,2\}$ if ${\rm lk}(K_i,K_j)={\rm lk}(K'_i,
K'_j)$ for all $i,j \in\{1,2,3,4\}$.
\end{lm}

\begin{proof} This is obvious from the definition of consistent orientations
of $S_{ij}$ and $S'_{ij}$ and
$${\rm lk}(\epsilon_iK_i,\epsilon_jK_j)=\epsilon_i\epsilon_j{\rm lk}(K_i,K_j).$$
\end{proof}

\begin{lm} Let $\textbf{\textit{J}}$ be the spherical tangle shown in Figure
15. Let $p_1,p_2,p_3,p_4$ be the number of half twists inside of
the balls marked by 1,2,3,4, respectively. Then
$${\rm det}\,F(\textbf{\textit{J}})=(p_1p_4-p_2p_3)^2.$$
\end{lm}

\begin{proof} This is by a direct calculation. We have
$$F(\textbf{\textit{J}})=\left[
\begin{matrix}p_1p_2p_3+p_1p_2p_4+p_1p_3p_4+p_2p_3p_4 &
-p_1p_3-p_1p_4-p_2p_3-p_2p_4 \\
p_1p_2+p_1p_4+p_3p_2+p_3p_4 & -p_1-p_2-p_3-p_4
\end{matrix}\right].$$
We calculate the determinant of $F( \textbf{\textit{J}})$ and get
$(p_1p_4-p_2p_3)^2$.
\end{proof}

\bigskip
\centerline{\epsfxsize=3.3 in \epsfbox{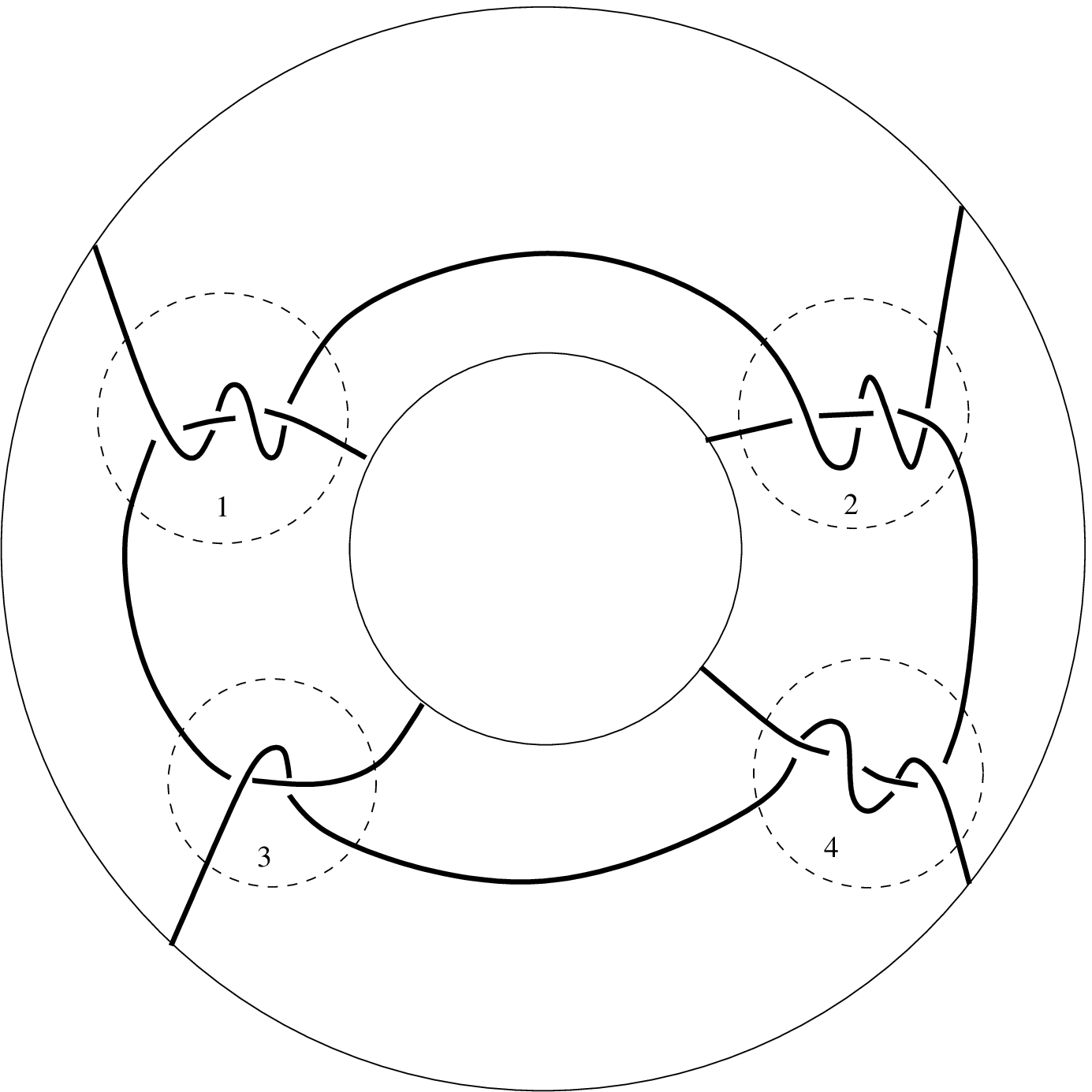}}
\medskip
\centerline{\small Figure 15. The spherical tangle
$\textbf{\textit{J}}$.}
\bigskip

\begin{lm} Let $S$ be a spherical tangle without closed
components. Then either there is a spherical tangle $S'$ which is
either $\textbf{\textit{I}}$-reducible or
$\textbf{\textit{J}}$-reducible, with $\textbf{\textit{J}}$ as
shown in Figure 15 or that $\textbf{\textit{J}}$ after some
possible operations of $(\cdot)^{r_1}$ and/or $(\cdot)^{r_2}$,
such that $\partial K'_i=\partial K_i$ and ${\rm lk}(K_i,K_j)={\rm
lk}(K'_i, K'_j)$ for all $i,j \in\{1,2,3,4\}$.
\end{lm}

\begin{proof} If $S$ has one component whose end points lie on different
boundary components of $S^2\times I$, then we can find such a
spherical tangle $S'$ that is $\textbf{\textit{I}}$-reducible.
Suppose now $S$ has no such components. The linking number of two
components whose end points lie on the same boundary component of
$S^2\times I$ can be realized by adding ball tangles. So we may
assume that there is no linking between such components. Then,
after some possible operations of $(\cdot)^{r_1}$ and/or
$(\cdot)^{r_2}$, we can take $S'$ as $\textbf{\textit{J}}$ with
the number of full twistings equal to the linking numbers of $S$.
\end{proof}

\begin{thm} If $S$ is a spherical tangle diagram without closed
components, then ${\rm det}\,F(S) \equiv n^2$ {\rm mod} 4 for some
integer $n$.
\end{thm}

\begin{proof} Let $S'$ be the spherical tangle in Lemma 4.31. Then $S_{ij}$
and $S'_{ij}$ have the same linking numbers, for each
$ij=11,12,21,22$. By Theorem 4.28, $S$ and $S'$ are
$\Delta$-equivalent. The theorem then follows from Lemma 4.27,
Theorem 4.20, and Lemma 4.30.
\end{proof}

Suppose now that a spherical tangle $S$ has closed components, we
have the following theorem.

\vskip 0.2in

\bigskip
\centerline{\epsfxsize=5 in \epsfbox{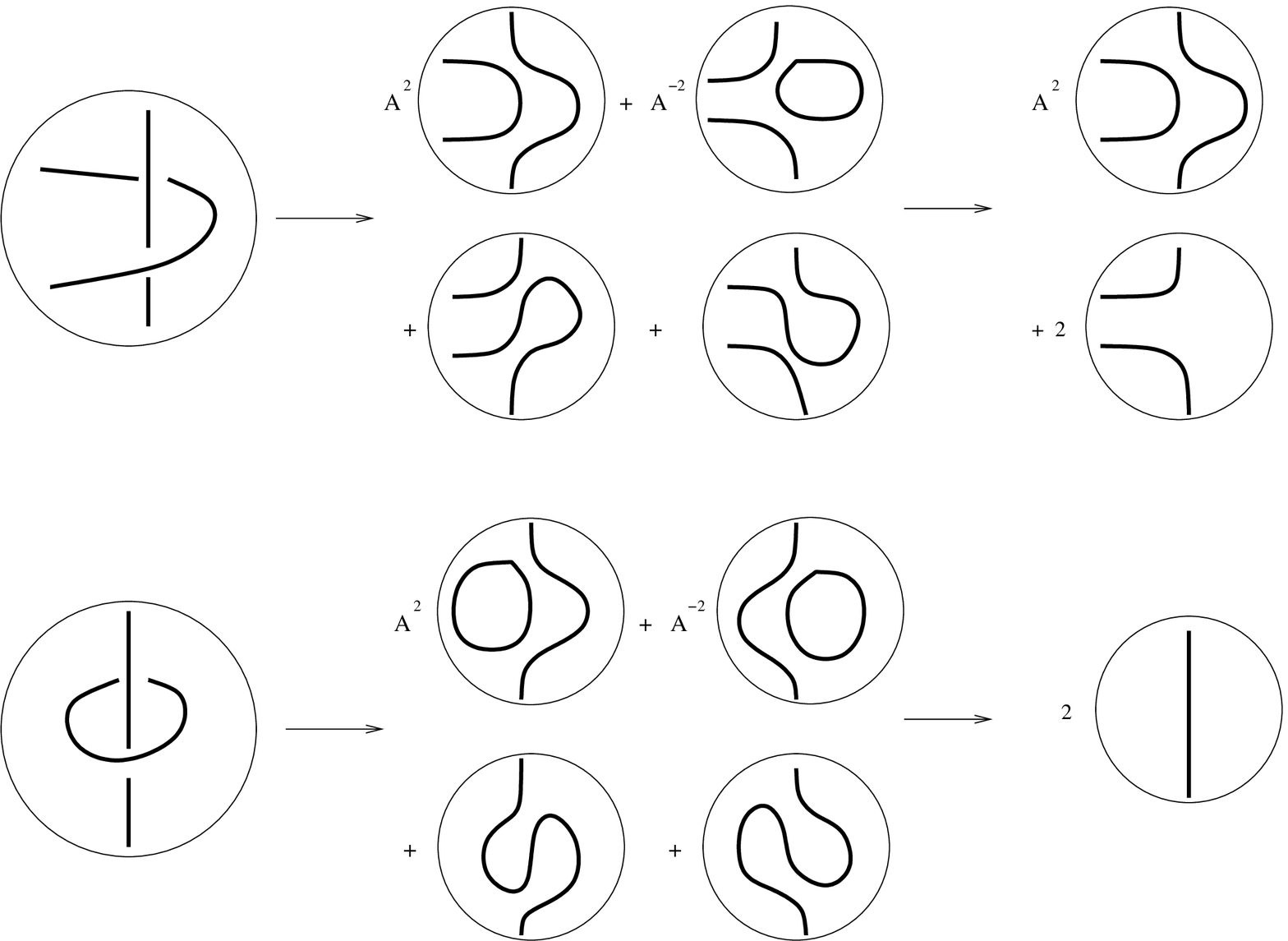}}
\medskip
\centerline{\small Figure 16. The case that $S$ has closed
components.}
\bigskip

\begin{thm} If $S \in \textbf{\textit{ST}}$ with
$F(S)=\left[\begin{matrix} \alpha & \gamma \\ \beta & \delta
\end{matrix}\right]$ and $S$ has closed components, then
$\alpha \equiv 0$ {\rm mod} 2, $\beta \equiv 0$ {\rm mod} 2,
$\gamma \equiv 0$ {\rm mod} 2, $\delta \equiv 0$ {\rm mod} 2.
\end{thm}

\begin{proof} See Figure 16 (top left), where a closed component
of $S$ is hooked with another component of $S$ as shown. Applying
the Kauffman skein relation to the local picture there, we get two
diagrams with coefficients $A^2$ and 2, respectively. The diagram
with coefficient 2 has one less closed component than $S$ and the
diagram with coefficient $A^2$ is obtained from $S$ by unhook the
closed component at that place. We keep performing this unhooking
process until the closed component is hooked with another
component only once, as illustrated in Figure 16 (bottom left).
Applying the Kauffman skein relation again, we can unhook this
closed component entirely and we end up with a diagram having one
less closed component and a coefficient 2. Note that this closed
component may itself being knotted. But this is not important
since a knotted closed component separating from other components
will make no contribution to the Kauffman bracket when $A=e^{i
\pi/4}$.

What we have shown is the fact that when $S$ has a closed
component, then 2 divides $\langle S_{ij}\rangle$ for all
$ij=11,12,21,22$. This proves the theorem.
\end{proof}

Note that $n^2 \equiv 0$ or 1 {\rm mod} 4. So combine Theorem 4.32
and Theorem 4.33, we get the following theorem.

\begin{thm} For every $S \in \textbf{\textit{ST}}$, either
${\rm det}\,F(S) \equiv 0$ {\rm mod} 4 or ${\rm det}\,F(S) \equiv
1$ {\rm mod} 4.
\end{thm}

From this theorem, we can conclude that there is no spherical
tangle $S$ such that
$$F(S)=\left[\begin{matrix} 1 & 0 \\ 0 & -1\end{matrix}\right],$$
since the determinant of the matrix above is not equal to 0 or 1
{\rm mod} 4.

Finally, let us indicate a direct way to compute the invariant
$F(\textbf{\textit{J}})$ of the spherical tangle
$\textbf{\textit{J}}$ in Figure 15, in the special case of
$p_1=p_2=p_4=-4$ and $p_3=2$. Check with the formula in Lemma
4.30. Suppose that $T^5,B^{(1)},B^{(2)},B^{(3)},B^{(4)}$ are the
$5,0,0,0,0$-punctured ball tangle diagrams in Figure 17,
respectively. Then
$$\textbf{\textit{J}}=T^5(B^{(1)},B^{(2)},B^{(3)},B^{(4)},\textbf{\textit{I}}).$$
By Theorem 3.8, $F(\textbf{\textit{J}})=F^5(T^5)[\eta^5]
(F^0(B^{(1)}),F^0(B^{(2)}),F^0(B^{(3)}),F^0(B^{(4)}),F(\textbf{\textit{I}}))$.

We have
$F^0(B^{(1)})=F^0(B^{(2)})=F^0(B^{(4)})=\left[\begin{matrix} -4 \\
1 \end{matrix}\right]$ and $F^0(B^{(3)})=\left[\begin{matrix} 2 \\
1 \end{matrix}\right]$.

\bigskip
\centerline{\epsfxsize=4.8 in \epsfbox{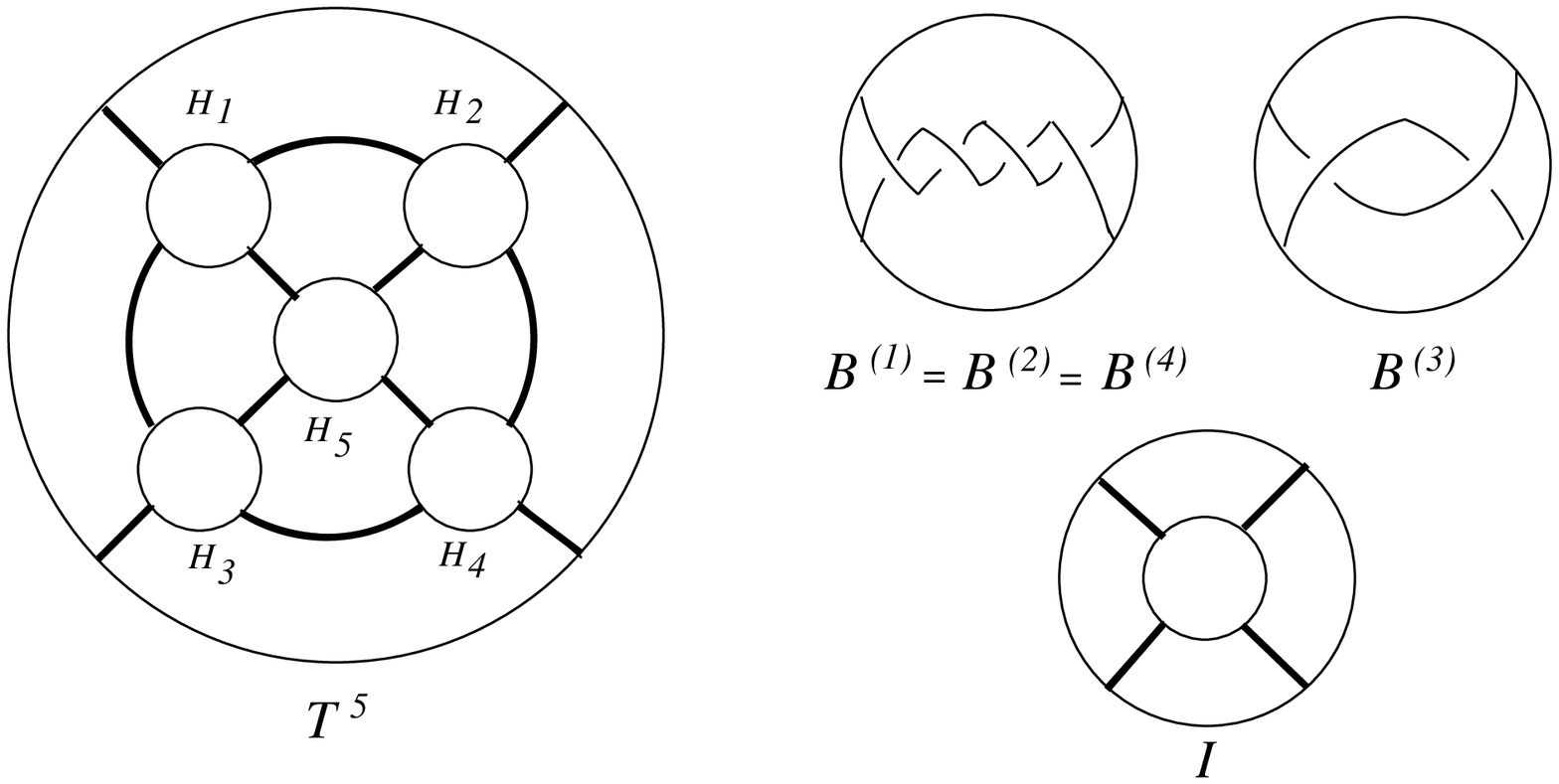}}
\medskip
\centerline{\small Figure 17. A decomposition of the spherical
tangle $\textbf{\textit{J}}$.}
\bigskip

First, let us compute $F^5(T^5)$ as the following steps:

1) The matrix $$\begin{pmatrix} \langle T_{1\alpha_1^5}^5 \rangle
& \cdots & \langle T_{1\alpha_{2^5}^5}^5 \rangle \\ \langle
T_{2\alpha_1^5}^5 \rangle & \cdots & \langle T_{2\alpha_{2^5}^5}^5
\rangle \end{pmatrix}$$ is
$$\begin{pmatrix} 0010 & 1000 & 1001 & 0100 & 1001 & 0100 & 0000 & 0000 \\
0000 & 0010 & 0000 & 1001 & 0010 & 0001 & 1001 & 0100 &
\end{pmatrix}.$$

2) Let $$F^5(T^5)=\left\{\begin{pmatrix} (-i)^{t_1}z\langle
T_{1\alpha_1^5}^5 \rangle & \cdots & (-i)^{t_{2^5}}z\langle T_{1\alpha_{2^5}^5}^5 \rangle \\
(-i)^{t_1}iz\langle T_{2\alpha_1^5}^5 \rangle & \cdots &
(-i)^{t_{2^5}}iz\langle T_{2\alpha_{2^5}^5}^5 \rangle
\end{pmatrix}\,|\,z \in \Phi\right\}
\cap M_{2\times2^5}(\mathbb Z).$$ Then the sequence $(t_k)_{1 \leq
k \leq 2^5}$ of exponents of $-i$ is $$\begin{matrix} 0112 & 1223
& 1223 & 2334 & 1223 & 2334 & 2334 & 3445 \end{matrix}.$$
Therefore, by taking $z=\pm i$, we have the invariant $F^5(T^5)$
as follows.
$$F^5(T^5)=\left[\begin{matrix}
0010 & 1000 & 100-1 & 0-100 & 100-1 & 0-100 & 0000 & 0000 \\
0000 & 0010 & 0000 & 100-1 & 0010 & 000-1 & 100-1 & 0-100 &
\end{matrix}\right].$$

Second, we compute $[\eta^5]
(F^0(B^{(1)}),F^0(B^{(2)}),F^0(B^{(3)}),F^0(B^{(4)}),F(\textbf{\textit{I}}))$
and describe it row-by-row as follows. That is, each pair of the
following means a row of the matrix $[\eta^5]
(F^0(B^{(1)}),F^0(B^{(2)}),F^0(B^{(3)}),F^0(B^{(4)}),F(\textbf{\textit{I}}))$.
$$128\,\,\,0;\, 0\,\,\,128;\, -32\,\,\,0;\, 0\,\,\,-32;\, 64\,\,\,0;\, 0\,\,\,64;\, -16\,\,\,0;\, 0\,\,\,-16;$$
$$-32\,\,\,0;\, 0\,\,\,-32;\, 8\,\,\,0;\, 0\,\,\,8;\, -16\,\,\,0;\, 0\,\,\,-16;\, 4\,\,\,0;\, 0\,\,\,4;$$
$$-32\,\,\,0;\, 0\,\,\,-32;\, 8\,\,\,0;\, 0\,\,\,8;\, -16\,\,\,0;\, 0\,\,\,-16;\, 4\,\,\,0;\, 0\,\,\,4;$$
$$8\,\,\,0;\, 0\,\,\,8;\, -2\,\,\,0;\, 0\,\,\,-2;\, 4\,\,\,0;\, 0\,\,\,4;\, -1\,\,\,0;\, 0\,\,\,-1.$$

Therefore, $$F(\textbf{\textit{J}})=F^5(T^5)[\eta^5]
(F^0(B^{(1)}),F^0(B^{(2)}),F^0(B^{(3)}),F^0(B^{(4)}),F(\textbf{\textit{I}}))$$
$$=\left[\begin{matrix} -32+64-32+0+0-32+0+0 & 0+0+0-8+16+0-8+16 \\
-16-16+0+8+0+8+0+0 & 0+0-4+0-4+0+2-4 \end{matrix}\right]$$
$$=\left[\begin{matrix} -32 & 16 \\ -16 & -10 \end{matrix}\right].$$

Also, we have
$${\rm det}\,F(\textbf{\textit{J}})=(-32)(-10)-16(-16)=576=24^2.$$
Thus, ${\rm det}\,F(\textbf{\textit{J}})$ is a square of integer.

\chapter{Open questions}

Here are some open questions that we are unable to answer at this
moment.

(1) Can one describe exactly the image of the invariant $F$ in
$PM_{2\times2}$?

The following two questions make the above question more specific

(1a) Is it true that ${\rm det}\,F(S)$ is the square of an integer
for any spherical tangle $S$?

(1b) Is it true that if a matrix $[A]$ in $PM_{2\times2}$ has its
determinant equal to the square of an integer, then there is a
spherical tangle $S$ such that $F(S)=[A]$?

The spherical tangle $\textbf{\textit{J}}$ does not look like
$\textbf{\textit{I}}$-reducible. But we do not know how to verify
this observation.

(2) How to show that $\textbf{\textit{J}}$ is not
$\textbf{\textit{I}}$-reducible? In general, given spherical
tangles $S$ and $S'$, how to show that $S$ is not $S'$-reducible?

\end{document}